\documentclass{bmcart}

\usepackage[utf8]{inputenc} 
\usepackage{array}
\usepackage{color}
\usepackage{tabularx}
\usepackage{graphicx}
\usepackage{amsmath}
\usepackage{amssymb}
\usepackage{amsfonts}	
\usepackage{moreverb}
\usepackage{dsfont}
\usepackage{bm}
\usepackage{multirow}
\usepackage{soul}
\usepackage{csquotes}

\newcommand\BibTeX{{\rmfamily B\kern-.05em \textsc{i\kern-.025em b}\kern-.08em
T\kern-.1667em\lower.7ex\hbox{E}\kern-.125emX}}

\newcommand{\bea}{\begin{eqnarray}}
\newcommand{\eea}{\end{eqnarray}}
\newcommand{\x}{\mbf{x}}

\newcommand{\mbf}[1]{\mathbf{#1}}			%

\renewcommand{\u}{\mathbf{Q}}
\newcommand{\w}{\mathbf{w}}
\newcommand{\q}{\mathbf{q}}
\newcommand{\Q}{\mathbf{Q}}
\newcommand{\F}{\mathbf{F}}

\renewcommand{\v}{\mathbf{v}}
\newcommand{\V}{\mathbf{V}}
\newcommand{\B}{\mathbf{B}}
\newcommand{\p}{\mathbf{p}}
\newcommand{\M}{\mathbf{M}}

\newcommand{\halb}{\frac{1}{2}}

\newcommand{\ar}{\phi_1\rho_1}
\newcommand{\arr}{\phi_2\rho_2}

\newcommand{\be}{\begin{equation}}
\newcommand{\ee}{\end{equation}}
\newcommand{\bdm}{\begin{displaymath}}
\newcommand{\edm}{\end{displaymath}}
\newcommand{\xL}{x_{i-\frac{1}{2}}}
\newcommand{\xR}{x_{i+\frac{1}{2}}}
\newcommand{\yL}{y_{j-\frac{1}{2}}}
\newcommand{\yR}{y_{j+\frac{1}{2}}}
\newcommand{\zL}{z_{k-\frac{1}{2}}}
\newcommand{\zR}{z_{k+\frac{1}{2}}}
\newcommand{\tn}{t^n}
\newcommand{\tnext}{t^{n+1}}

\newcommand{\aposteriori}{\textit{a posteriori} }

\startlocaldefs
\endlocaldefs

\begin{document}

\begin{frontmatter}

\begin{fmbox}
\dochead{Research}


\title{Efficient conservative ADER schemes based on WENO reconstruction and space-time predictor in primitive variables} 


\author[
   addressref={aff1},                   
   corref={aff1},                       
   email={olindo.zanotti@unitn.it}   
]{\inits{JE}\fnm{Olindo} \snm{Zanotti}}
\author[
   addressref={aff1},
   email={michael.dumbser@unitn.it}
]{\inits{JRS}\fnm{Michael} \snm{Dumbser}}


\address[id=aff1]{
  \orgname{Laboratory of Applied Mathematics, Department of Civil, Environmental and Mechanical Engineering, University of Trento}, 
  \street{Via Mesiano 77},                     %
  \postcode{38123}                                
  \city{Trento},                              
  \cny{Italy}                                    
}



\end{fmbox}


\begin{abstractbox}

\begin{abstract} 
We present a new version of conservative ADER-WENO finite volume schemes, in which both the high order spatial reconstruction as well as the time evolution 
of the reconstruction polynomials in the local space-time predictor stage are performed in \textit{primitive} variables, rather than in conserved ones. To obtain a 
conservative method, the underlying finite volume scheme is still written in terms of the cell averages of the conserved quantities. Therefore, 
our new approach performs the spatial WENO reconstruction \textit{twice}: the \textit{first} WENO reconstruction is carried out on the known \textit{cell averages} 
of the conservative variables. The WENO polynomials are then used at the cell centers to compute \textit{point values} of the \textit{conserved variables}, which are 
subsequently converted into \textit{point values} of the \textit{primitive variables}. This is the only place where the conversion from 
conservative to primitive variables is needed in the new scheme. Then, a \textit{second} WENO reconstruction is performed on the point values of the primitive 
variables to obtain piecewise high order reconstruction polynomials of the primitive variables. 
The reconstruction polynomials are subsequently evolved in time with a \textit{novel} space-time finite element predictor that is directly applied to the governing PDE 
written in \textit{primitive form}. The resulting space-time polynomials of the primitive variables can then be directly used as input for the numerical fluxes at 
the cell boundaries in the underlying \textit{conservative} finite volume scheme. Hence, the number of necessary conversions from the conserved to the primitive variables is 
reduced to just \textit{one single conversion} at each cell center. 
We have verified the validity of the new approach over a wide range of hyperbolic systems, including the classical Euler equations of gas dynamics, the special relativistic 
hydrodynamics (RHD) and ideal magnetohydrodynamics (RMHD) equations, as well as the Baer-Nunziato model for compressible two-phase flows. 
In all cases we have noticed that the new ADER schemes provide \textit{less oscillatory solutions} when compared to ADER finite volume schemes based on the reconstruction 
in conserved variables, especially for the RMHD and the Baer-Nunziato equations. 
For the RHD and RMHD equations, the overall accuracy is improved and the CPU time is reduced by about $25 \%$. 
Because of its increased accuracy and due to the reduced computational cost, we recommend to use this version of ADER as the standard one in the relativistic framework. 
At the end of the paper, the new approach has also been extended to ADER-DG schemes on space-time adaptive grids (AMR). 
\end{abstract}




\end{abstractbox}
%

\end{frontmatter}



\section{Introduction} \label{sec:introduction}
Since their introduction by Toro and Titarev 
\cite{toro1,toro3,toro4,titarevtoro,Toro:2006a}, ADER (arbitrary high order derivatives) schemes for hyperbolic partial differential equations (PDE) 
have been improved and developed along different directions.
A key feature of these methods is their ability to achieve uniformly high order of accuracy in space and time in 
a single step, without the need of intermediate Runge-Kutta stages \cite{Pareschi2005,Puppo2015}, by exploiting the approximate solution of 
a Generalized Riemann Problem (GRP) at cell boundaries. 
ADER schemes have been first conceived within the finite volume (FV) framework, but they were soon extended also 
to the discontinuous Galerkin (DG) finite element framework \cite{dumbser_jsc,taube_jsc} and to a unified formulation 
of FV and DG schemes, namely the so-called $\mathbb{P}_N\mathbb{P}_M$ approach \cite{Dumbser2008}. 
In the original ADER approach by Toro and Titarev,  the approximate solution of the GRP is obtained through the solution of a 
conventional Riemann problem between the boundary-extrapolated values, and a sequence of linearized Riemann 
problems for the spatial derivatives. The required time derivatives in the GRP are obtained via the so-called 
Cauchy-Kowalevski procedure, which consists in replacing the time derivatives of the Taylor expansion at each interface 
with spatial derivatives of appropriate order, by resorting to the strong differential form of the PDE. Such an approach, 
though formally elegant, becomes prohibitive or even impossible as the complexity of the equations increases, especially for 
multidimensional problems and for relativistic hydrodynamics and magneto-hydrodynamics. 
On the contrary, in the modern reformulation of ADER~\cite{DumbserEnauxToro,Dumbser2008,Balsara2013934}, the approximate solution of 
the GRP is achieved by first evolving the data locally inside each cell through a \emph{local space-time discontinuous Galerkin predictor} (LSDG) 
step that is based on a weak form of the PDE, and, second, by solving a 
sequence of classical Riemann problems along the time axis at each element interface. This approach
has the additional benefit that it can successfully cope with stiff source terms in the equations, a fact which is often encountered in physical 
applications. 
For these reasons, ADER schemes  have been applied to real physical problems mostly in their modern version. Notable examples of applications include  
the study of Navier--Stokes equations, with or without chemical reactions~\cite{HidalgoDumbser,DumbserNSE}, geophysical flows~\cite{ADERNC}, 
complex three-dimensional free surface flows~\cite{Dumbser2013},
relativistic magnetic reconnection~\cite{DumbserZanotti,Zanotti2011b}, and the study of the Richtmyer--Meshkov instability in the relativistic regime~\cite{Zanotti2015b}.
In the last few years,  ADER schemes have  been enriched with several additional properties, reaching a high level of flexibility. 
First of all, ADER schemes have been soon extended to deal with non-conservative systems of hyperbolic PDE ~\cite{Hidalgo2009,ADERNC,AMR3DNC}, by resorting to path-conservative methods~\cite{Pares2004,pares2006}. %
ADER  schemes have also been extended to the Lagrangian framework, in which they are currently applied to the solution of multidimensional problems on unstructured meshes
for various systems of equations, ~\cite{Lagrange2D,LagrangeNC,LagrangeMDRS,LagrangeMHD,Lagrange3D}. 
On another side, ADER schemes have been combined with Adaptive Mesh Refinement (AMR) techniques ~\cite{AMR3DCL,Zanotti2015}, exploiting the
local properties of the discontinuous Galerkin predictor step, which is applied cell-by-cell irrespective of the level of refinement of the neighbour cells.
Moreover, ADER schemes have also been used in combination with Discontinuous Galerkin methods, even in the presence of shock waves and other discontinuities within the flow, 
thanks to a novel \textit{a posteriori} sub-cell finite volume limiter technique based on the MOOD approach \cite{CDL1,CDL2}, that is designed to stabilize the 
discrete solution wherever the DG approach fails and produces spurious oscillations or negative densities and pressures ~\cite{Dumbser2014,Zanotti2015d,Zanotti2015c}. 

The various implementations of ADER schemes mentioned so far differ under several aspects, but they all share the following common features: they apply the local space-time discontinuous Galerkin predictor to the conserved variables, which in turn implies that, if a WENO finite volume scheme is used, the spatial WENO reconstruction is also performed in terms of the conserved variables. 
Although this may be regarded as a reasonable choice, 
it has two fundamental drawbacks. The first one has to do with the fact that, as shown by \cite{Munz1986},
the reconstruction in conserved variables provides the worst shock capturing fidelity when compared to the reconstruction performed either in primitive or in characteristic variables.
The second drawback is instead related to computational performance. Since the computation of the numerical fluxes requires the calculation of integrals via Gaussian quadrature,
the physical fluxes must necessarily be computed at each space-time Gauss--Legendre quadrature point. However, there are systems of equations (e.g. the relativistic hydrodynamics or magnetohydrodynamics equations) for which the physical fluxes can only be written in terms of the primitive variables. As a result, a conversion from the conserved to the primitive variables is necessary for the calculations of the fluxes, and this operation, which is never analytic for such systems of equations, is rather expensive. For these reasons it would be very desirable to have an ADER scheme in which both the reconstruction and the subsequent local space-time discontinuous Galerkin predictor are performed in primitive variables. It is the aim of the present paper to explore this possibility.
It is also worth stressing that in the context of high order finite difference Godunov methods, based on traditional Runge--Kutta discretization in time, the reconstruction 
in primitive variables has been proved to be very successful by \cite{DelZanna2007} in their ECHO general relativistic code (see also \cite{Bucciantini2011,Zanotti2011}). 
In spite of the obvious differences among the numerical schemes adopted, the approach that we propose here and the ECHO-approach share the common feature of requiring a single (per cell)
conversion from the conserved to the primitive variables.

The plan of the paper is the following: in Sect.~\ref{sec:num-approach} we describe the numerical method, with particular emphasis on Sect.~\ref{sec:WENO_reconstruction} and on Sect.\ref{sec:localDG}, 
where the spatial reconstruction strategy and the local space-time discontinuous Galerkin predictor in primitive variable are described. 
The results of our new approach are presented in Sect.~\ref{sec:num-tests}
for a set of four different systems of equations. In Sect.~\ref{sec:extension} we show that the new strategy can also be extended to pure 
Discontinuous Galerkin schemes, even in the presence of space-time adaptive meshes (AMR).
Finally, Sect.~\ref{sec:conclusions} is devoted to the conclusions of the work.

\section{Numerical Method}
\label{sec:num-approach}

We present our new approach for purely regular Cartesian meshes, although there is no conceptual reason preventing the extension to general curvilinear or unstructured
meshes, which may be considered in future studies.

\subsection{Formulation of the equations}
We consider hyperbolic systems of balance laws that contain both conservative and non-conservative terms, i.e.
\begin{equation}	\label{NCsyst}
	\frac{\partial \u}{\partial t}+\nabla\cdot\bf F(\u)+\bf B(\u)\cdot\nabla \u=\bf S(\u)\,,
\end{equation}
where $\u \in \Omega_\Q \subset \mathds{R}^\nu$ is the state vector of the $\nu$ 
{\em conserved variables}, which, for the typical gas dynamics equations, are related to the conservation of mass, momentum and energy. 
${\bf F}(\u)=[{\bf f}^x(\u),{\bf f}^y(\u),{\bf f}^z(\u)]$ is the flux tensor\footnote{
Since we adopt Cartesian coordinates,
${\bf f}^x(\u),{\bf f}^y(\u),{\bf f}^z(\u)$ express the fluxes along the 
$x$, $y$ and $z$ directions, respectively.}
 for the conservative part of the PDE system,
while ${\bf B}(\u)=[\textbf{B}_x(\u),\textbf{B}_y(\u),\textbf{B}_z(\u)]$ represents the non-conservative part of it. 
Finally, $\bf S(\u)$ is the vector of the source terms, which may or may  not be present.
In the follow up of our discussion it is convenient to recast the system (\ref{NCsyst})
in quasilinear form as
\begin{equation}	\label{Csyst}
\frac{\partial \u}{\partial t}+ \bf{A}(\u) \cdot\nabla \u = \textbf{S}(\u)\,,
\end{equation}
where  ${\bf A}(\u)=[{\bf A}_x,{\bf A}_y,{\bf A}_z]=\partial {\bf F}(\u)/\partial \u+{\bf B}(\u)$ accounts for both the conservative and the non-conservative contributions. 
As we shall see below, a proper discretization of Eq.~(\ref{Csyst}) can provide the time evolution of the conserved variables  $\u$, but when the {\em primitive variables} $\V$ are adopted instead,  
Eq.~(\ref{Csyst}) translates into
\begin{equation}	
\label{Primsyst}
\frac{\partial \V}{\partial t}+ {\bf C}(\u) \cdot\nabla \V = \left( \frac{\partial \u}{\partial \V}\right)^{-1}{\bf S}(\u)\,,
\quad \textnormal{ with } \quad 
{\bf C}(\u)=\left( \frac{\partial \u}{\partial \V}\right)^{-1} \bf{A}(\u) \left( \frac{\partial \u}{\partial \V} \right). 
\end{equation}
In the following we suppose that the conserved variables $\u$ can always be written \textit{analytically} in terms of the primitive variables $\V$, i.e. 
the functions 
\be
\label{eq:prim2cons}
\u=\u(\V)
\ee
are supposed to be analytic for all PDE systems under consideration. On the contrary, the conversion from the conserved to the primitive variables, henceforth 
the {\em cons-to-prim conversion}, is not always available in closed form, i.e. the functions 
\be
\label{eq:cons2prim}
\V=\V(\u)
\ee
may \textit{not} be analytic (e.g. for relativistic hydrodynamics and magnetohydrodynamics to be discussed in Sect.~\ref{sec:RMHD}), 
thus requiring an approximate numerical solution.
As a result, the matrix $\left( \frac{\partial \u}{\partial \V}\right)^{-1}$, which in principle could be simply computed as 
\be
\left( \frac{\partial \u}{\partial \V}\right)^{-1}=\left( \frac{\partial \V}{\partial \u}\right)\,,
\ee
in practice it cannot be obtained in this manner, but it must be computed as 
\be
\left( \frac{\partial \u}{\partial \V}\right)^{-1}= \M^{-1} \,,
\ee
where we have introduced the notation
\be
\M = \left( \frac{\partial \u}{\partial \V}\right)\,,
\ee
which will be used repeatedly below. Since $\u(\V)$ is supposed to be analytic, the matrix $\M$ can be easily computed. 
Equation (\ref{NCsyst}) will serve us as the master equation to evolve the cell averages of the conserved variables $\u$ via 
a standard finite volume scheme. However, both the spatial WENO reconstruction and the subsequent LSDG predictor will act on the primitive variables $\V$, 
hence relying on the alternative formulation given by Eq.~(\ref{Primsyst}). The  necessary steps to obtain such a scheme are described in the 
Sections \ref{sec:FV}--\ref{sec:localDG} below.

\subsection{The finite volume scheme}
\label{sec:FV}
In Cartesian coordinates, we discretize
the computational domain $\Omega$ through space-time control volumes 
${\mathcal I}_{ijk}=I_{ijk}\times [t^n,t^n+\Delta t]=[\xL,\xR]\times[\yL,\yR]\times[\zL,\zR]\times [\tn,\tn+\Delta t]$, with $\Delta x_i=\xR-\xL$, $\Delta y_j=\yR-\yL$, $\Delta z_k=\zR-\zL$ and $\Delta t=\tnext-\tn$. 
Integration of Eq.~(\ref{NCsyst}) over  ${\mathcal I}_{ijk}$ yields the usual  finite volume discretization
\begin{eqnarray} \label{FVformula}
\label{eq:finite_vol}
{\bar \u}_{ijk}^{n+1}&=&{\bar \u}_{ijk}^{n}-
\frac{\Delta t}{\Delta x_i}\left[\left({\textbf f}^x_{i+\halb,j,k} \,\, -{\textbf f}^x_{i-\halb,j,k}\right)+\frac{1}{2} \,\,
\left({{D}}^x_{i+\halb,j,k} +{{D}}^x_{i-\halb,j,k}\right)\right]\nonumber\\
&& \hspace{9mm} -\frac{\Delta t}{\Delta y_j}\left[\left({\textbf f}^y_{i,j+\halb,k}-{\textbf f}^y_{i,j-\halb,k}\right)+\frac{1}{2}
\left({{D}}^y_{i,j+\halb,k}+{{D}}^y_{i,j-\halb,k}\right)\right]\nonumber\\
&& \hspace{9mm} -\frac{\Delta t}{\Delta z_k}\left[\left({\textbf f}^z_{i,j,k+\halb}-{\textbf f}^z_{i,j,k-\halb}\right)+\frac{1}{2}
\left({{D}}^z_{i,j,k+\halb}+{{D}}^z_{i,j,k-\halb}\right)\right]
 + \Delta t({\bf \bar{S}}_{ijk}- {\bf \bar{P}}_{ijk})\,, \nonumber \\ 
\end{eqnarray}
where the cell average 
\begin{equation}
{\bar \u}_{ijk}^{n}=\frac{1}{\Delta x_i}\frac{1}{\Delta y_j}\frac{1}{\Delta z_k}\int_{x_{i-\halb}}^{x_{i+\halb}}\int_{y_{j-\halb}}^{y_{j+\halb}}\int_{z_{k-\halb}}^{z_{k+\halb}}{\u}(x,y,z,t^n)dz\,dy\,\,dx
\end{equation}
is the spatial average of the vector of conserved quantities at time $\tn$. In Eq.~(\ref{eq:finite_vol}) we recognize two different sets of terms, namely
 those  due to the conservative part of the system 
(\ref{NCsyst}), and those coming from the non-conservative part of it. In the former set we include the three time-averaged fluxes
\begin{equation}
\label{averF}
{\bf f}^x_{i+\halb,jk}= \frac{1}{\Delta t}\frac{1}{\Delta y_j}\frac{1}{\Delta z_k} \hspace{-1mm}  \int \limits_{t^n}^{t^{n+1}} \! \int \limits_{y_{j-\halb}}^{y_{j+\halb}} \! \int \limits_{z_{k-\halb}}^{z_{k+\halb}} \hspace{-1mm} 
{\bf \tilde f}^x \! \left({\v}_h^-(x_{i+\halb},y,z,t),{\v}_h^+(x_{i+\halb},y,z,t)\right) dz \, dy \, dt, 
\end{equation}
\begin{equation}
\label{averG}
{\bf f}^y_{i,j+\halb,k}=\frac{1}{\Delta t}\frac{1}{\Delta x_i}\frac{1}{\Delta z_k} \hspace{-1mm}  \int \limits_{t^n}^{t^{n+1}} \! \int \limits_{x_{i-\halb}}^{x_{i+\halb}} \! \int \limits_{z_{k-\halb}}^{z_{k+\halb}} \hspace{-1mm} 
{\bf \tilde f}^y \! \left({\v}_h^-(x,y_{j+\halb},z,t),{\v}_h^+(x,y_{j+\halb},z,t)\right) dz\,dx\,dt, \\
\end{equation}
\begin{equation}
\label{averH}
{\bf f}^z_{ij,k+\halb}=\frac{1}{\Delta t}\frac{1}{\Delta x_i}\frac{1}{\Delta y_j} \hspace{-1mm}  \int \limits_{t^n}^{t^{n+1}} \! \int \limits_{x_{i-\halb}}^{x_{i+\halb}} \! \int \limits_{y_{j-\halb}}^{y_{j+\halb}} \hspace{-1mm}  
{\bf \tilde f}^z\! \left({\v}_h^-(x,y,z_{k+\halb},t),{\v}_h^+(x,y,z_{k+\halb},t)\right) dy\,dx\,dt  
\end{equation}
and the space-time averaged source term 
\begin{equation}
\label{source:S}
{\bf \bar{S}}_{ijk}=\frac{1}{\Delta t}\frac{1}{\Delta x_i}\frac{1}{\Delta y_j}\frac{1}{\Delta z_k}\int \limits_{t^n}^{t^{n+1}}\int \limits_{x_{i-\halb}}^{x_{i+\halb}}\int \limits_{y_{j-\halb}}^{y_{j+\halb}}\int \limits_{z_{k-\halb}}^{z_{k+\halb}}{ \bf S}\left(\v_h(x,y,z,t)\right) dz\,dy\,dx\,dt\,.
\end{equation}
We emphasize that
the terms ${\v}_h$ in Eq.~(\ref{averF})--(\ref{source:S}), as well as in the few equations below, 
are piecewise space-time polynomials of degree $M$ in {\em primitive variables}, computed according to a suitable LSDG predictor based on the formulation (\ref{Primsyst}), as we will discuss in Sect.~\ref{sec:localDG}. This marks a striking difference with respect to traditional ADER schemes, in which such polynomials are instead
computed in conserved variables and are denoted as ${\bf q}_h$ (see, e.g. \cite{HidalgoDumbser}).
The integrals over the smooth part of the non-conservative terms in Eq.~(\ref{eq:finite_vol}) yield the following contribution, 
\begin{equation}
	{\bf{\bar{P}}}_{ijk}=\frac{1}{\Delta t}\frac{1}{\Delta x_i}\frac{1}{\Delta y_j}\frac{1}{\Delta z_k}\int \limits_{t^n}^{t^{n+1}}\int \limits_{x_{i-\halb}}^{x_{i+\halb}}\int \limits_{y_{j-\halb}}^{y_{j+\halb}}\int \limits_{z_{k-\halb}}^{z_{k+\halb}}{\bf B}({\v}_h) \M \, \nabla {\v}_h \,dz\,dy\,dx\,dt\,,
\end{equation}
while the {\em jumps} across the element boundaries are treated within the framework of path-conservative schemes 
\cite{Pares2004,pares2006,Munoz2007,Castro2006,Castro2008,NCproblems} based on the Dal Maso--Le Floch--Murat theory \cite{DLMtheory} as 
\begin{equation}
{{D}}^x_{i+\halb,j,k} \!=\!\! \frac{1}{\Delta t}\frac{1}{\Delta y_j}\frac{1}{\Delta z_k} \hspace{-1mm}  \int \limits_{t^n}^{t^{n+1}} \! \int \limits_{y_{j-\halb}}^{y_{j+\halb}} \! \int \limits_{z_{k-\halb}}^{z_{k+\halb}} \hspace{-2mm}  
{\cal{D}}_x  \! \left({\v}_h^-(x_{i+\halb},y,z,t),{\v}_h^+(x_{i+\halb},y,z,t)\right) dz \, dy \, dt, 
\end{equation} 
\begin{equation} 
{{D}}^y_{i,j+\halb,k} \!=\!\! \frac{1}{\Delta t}\frac{1}{\Delta x_i}\frac{1}{\Delta z_k} \hspace{-1mm}  \int \limits_{t^n}^{t^{n+1}} \! \int \limits_{x_{i-\halb}}^{x_{i+\halb}} \! \int \limits_{z_{k-\halb}}^{z_{k+\halb}} \hspace{-2mm}  
{\cal{D}}_y  \! \left({\v}_h^-(x,y_{j+\halb},z,t),{\v}_h^+(x,y_{j+\halb},z,t)\right) dz \, dx \, dt, 
\end{equation} 
\begin{equation} 
{{D}}^z_{i,j,k+\halb} \!=\!\! \frac{1}{\Delta t}\frac{1}{\Delta x_i}\frac{1}{\Delta y_j} \hspace{-1mm}  \int \limits_{t^n}^{t^{n+1}} \! \int \limits_{x_{i-\halb}}^{x_{i+\halb}} \! \int \limits_{y_{j-\halb}}^{y_{j+\halb}} \hspace{-2mm}  
{\cal{D}}_z  \! \left({\v}_h^-(x,y,z_{k+\halb},t),{\v}_h^+(x,y,z_{k+\halb},t)\right) dy \, dx \, dt\,.
\end{equation}
According to this approach, the following path integrals must be prescribed
\begin{equation}
  \label{eqn.dm}
  {\cal{D}}_i({\v}_h^-,{\v}_h^+) = \int_0^1{{\bf B}_i \left(\Psi({\v}_h^-,{\v}_h^+,s)\right) \M\left(\Psi({\v}_h^-,{\v}_h^+,s)\right) \frac{\partial\Psi}{\partial s}ds}, 
	\qquad i \in \left\{ x,y,z \right\}, 
\end{equation} 
where $\Psi(s)$ is a path joining the left and right boundary extrapolated states ${\v}_h^-$ and ${\v}_h^+$ in state space of the primitive variables. 
The simplest option is to use a straight-line segment path 
\begin{equation}	\label{segment}
 \Psi = \Psi({\v}_h^-,{\v}_h^+,s)
={\v}_h^- + s({\v}_h^+ - {\v}_h^-)\,,\qquad 0\leq s \leq 1. 
\end{equation}
Pragmatic as it is\footnote{See \cite{MuellerWB} for more sophisticated paths.}, 
the choice of the path (\ref{segment}) allows to evaluate the terms ${\cal{D}}_i$ in (\ref{eqn.dm}) as 
\begin{equation}
\label{Osher-D}
{\cal{D}}_i({\v}_h^-,{\v}_h^+) = \left( \int_0^1{ {\bf B}_i \left(\Psi({\v}_h^-,{\v}_h^+,s)\right) \M\left(\Psi({\v}_h^-,{\v}_h^+,s)\right) ds} \right) \left( {\v}_h^+ - {\v}_h^- \right)\,,
\end{equation}
that we compute through 
a three-point Gauss-Legendre formula \cite{USFORCE2,OsherNC,OsherUniversal}. 
The computation of the numerical fluxes ${\bf \tilde f}^i$ in Eq.~(\ref{averF}) requires the use of an approximate Riemann solver, see \cite{toro-book}.
In this work we have limited our attention to 
a local Lax-Friedrichs flux (Rusanov flux) and to the Osher-type flux proposed in \cite{OsherUniversal,OsherNC,ApproxOsher}. 
Both of them can be written formally as 
\begin{equation}
  {\bf \tilde f}^i = \frac{1}{2}\left( \mathbf{f}^i(\v_h^-) + \mathbf{f}^i(\v_h^+) \right) - \frac{1}{2}\mathbf{D}_i \, \widetilde \M \, \left( \v_h^+ - \v_h^- \right)\,, \qquad i \in \left\{ x,y,z \right\} 
\label{eqn.numerical.flux} 
\end{equation} 
where $\mathbf{D}_i \geq 0$ is a positive-definite dissipation matrix that depends on the chosen Riemann solver. For the Rusanov flux it simply reads 
\begin{equation}
 \mathbf{D}^{\textnormal{Rusanov}}_i =  |s_{\max}| \mathbf{I} \,,
\label{eqn.rusanov} 
\end{equation}
where $|s_{\max}|$ is the maximum absolute value of the eigenvalues admitted by the PDE and $\mathbf{I}$ is the identity matrix. The matrix $\widetilde\M$ is a \textit{Roe matrix} that 
allows to write the jumps in the conserved variables in terms of the jump in the primitive variables, i.e. 
\be
\q_h^+ - \q_h^-  =  \Q(\v_h^+) - \Q(\v_h^-)  = \widetilde\M \, \left( \v_h^+ - \v_h^- \right). 
\ee
Since $\M = \partial \Q / \partial \V$, the Roe matrix $\widetilde \M$ can be easily defined by a path integral as
\begin{equation}
\Q(\v_h^+) - \Q(\v_h^-)  = \int \limits_0^1 \M (\Psi(\v_h^-,\v_h^+,s)) \frac{\partial\Psi}{\partial s} ds = \widetilde\M \, \left( \v_h^+ - \v_h^- \right), 
\label{eqn.MRoe1} 
\end{equation} 
which in the case of the simple straight-line segment path \eqref{segment} leads to the expression 
\begin{equation}
\widetilde \M = \int \limits_0^1 \M (\Psi(\v_h^-,\v_h^+,s)) ds. 
\label{eqn.MRoe2} 
\end{equation} 
In the case of the Osher-type flux, on the other hand, the dissipation matrix reads 
\begin{eqnarray}
\label{eqn.osher}
 \mathbf{D}^{\textnormal{Osher}}_i =  \int \limits_0^1 |{\bf A}_i(\Psi(\v_h^-,\v_h^+,s))| ds\,,
\end{eqnarray}
with the usual definition of the matrix absolute value operator 
\begin{equation}
|{\bf A}|={\bf R}|{\bf \Lambda}|{\bf R}^{-1}\,,\qquad  |{\bf \Lambda}|={\rm diag}(|\lambda_1|, |\lambda_2|, \ldots, |\lambda_\nu|)\,.
\end{equation}
The path $\Psi$ in Eq. (\ref{eqn.osher}) and \eqref{eqn.MRoe2} is the same segment path adopted in (\ref{segment}) for the computation of the jumps ${\cal{D}}_i$.

\subsection{A novel WENO reconstruction in primitive variables}
\label{sec:WENO_reconstruction}
%
Since we want to compute the time averaged fluxes [c.f. Eq.~(\ref{averF})--(\ref{averH})] and the space-time averaged sources
[c.f. Eq.~(\ref{source:S})] directly from the primitive variables $\V$, it is necessary to reconstruct a WENO polynomial in primitive variables.
However, the underlying finite volume scheme \eqref{FVformula} will still advance in time the cell averages of the conserved variables 
$\bar{\Q}_{ijk}^n$, which are the only known input quantities at the reference time level $t^n$. Hence, the whole procedure is performed through 
the following three simple steps:  
\begin{enumerate}
\item We perform a {\em first} standard spatial WENO reconstruction of the conserved variables starting from the cell averages ${\bar \u}_{ijk}^{n}$. This allows to obtain a reconstructed polynomial $\w_h(x,y,z,t^n)$ in conserved variables valid within each cell.
\item Since $\w_h(x,y,z,t^n)$ is defined at any point inside the cell, we simply \textit{evaluate} it at the cell center in order to obtain the \textit{point value} 
$\Q_{ijk}^{n}= \w_h(x_i,y_j,z_k,t^n)$. This conversion from cell averages ${\bar \u}_{ijk}^{n}$ to point values $\Q_{ijk}^{n}$ is the \textbf{main key idea} of our new method, 
since the simple identity  $\Q_{ijk}^{n} = {\bar \u}_{ijk}^{n}$ is valid only up to second order of accuracy! 
After that, we perform a conversion from the point-values of the conserved variables to the point-values in primitive variables, i.e. we apply Eq.~(\ref{eq:cons2prim}), 
thus obtaining the corresponding primitive variables $\V_{ijk}^{n} = \V(\Q_{ijk}^{n})$ at each cell center. This is the only step in the entire algorithm that needs a conversion 
from the conservative to the primitive variables. 
\item Finally, from the point-values of the  primitive  variables at the cell centers, we perform a {\em second} WENO reconstruction to obtain a reconstruction polynomial in  
\textit{primitive variables}, denoted as $\p_h(x,y,z,t^n)$. This polynomial is then used as the initial condition for the new local space--time DG predictor in primitive
variables described in Sect.~\ref{sec:localDG}.
\end{enumerate}
As for the choice of the spatial WENO reconstruction, we have adopted a dimension-by-dimension reconstruction strategy, discussed in full details in our previous works
(see \cite{AMR3DCL,AMR3DNC,Zanotti2015}). Briefly, we first introduce space-time reference coordinates $\xi,\eta,\zeta,\tau\in[0,1]$, defined by
\begin{equation}
\label{eq:xi}
x = x_{i-\halb} + \xi   \Delta x_i, \quad 
y = y_{j-\halb} + \eta  \Delta y_j, \quad 
z = z_{k-\halb} + \zeta \Delta z_k, \quad
t = t^n + \tau \Delta t\,,
\end{equation} 
and, along each spatial direction, we define a basis of polynomials
$\{\psi_l(\lambda)\}_{l=1}^{M+1}$, each of degree $M$, 
formed by the $M+1$ Lagrange interpolating polynomials, 
that pass through the $M+1$
Gauss-Legendre quadrature nodes  $\{\mu_k\}_{k=1}^{M+1}$. According to the WENO philosophy, 
a number of stencils is introduced such that the final polynomial is a
data-dependent nonlinear combination of the polynomials computed from each stencil. 
Here, we use a fixed number $N_s$ of one-dimensional stencils, namely $N_s=3$ 
for odd order schemes (even polynomials of degree $M$), and $N_s=4$ for 
even order schemes (odd polynomials of degree $M$). For example, 
focusing on the $x$ direction for convenience,
every stencil along $x$ 
is formed by the union of $M+1$ adjacent cells, i.e.
\begin{equation}
\label{eqn.stencildef}  
\mathcal{S}_{ijk}^{s,x} = \bigcup \limits_{e=i-L}^{i+R}
        {I_{ejk}}, \quad 
\end{equation}
where $L=L(M,s)$ and $R=R(M,s)$ are the 
spatial extension of the stencil to the left and to the
right.\footnote{See Appendix A of \cite{Zanotti2015} for a graphical representation.}

Now, an important difference emerges depending on whether we are reconstructing the conserved or the primitive variables.
In the former case, corresponding to the computation of $\w_h(x,y,z,t^n)$ at step $1$ above,
we require that the reconstructed polynomial must preserve the \textit{cell-averages} of the 
\textit{conserved variables} over each element $I_{ijk}$. Since the polynomials reconstructed along the $x$ direction 
can be written as 
\begin{equation}
\label{eqn.recpolydef.x} 
 \w^{s,x}_h(x,t^n) = \sum \limits_{r=0}^M \psi_r(\xi) \hat \w^{n,s}_{ijk,r} := \psi_r(\xi) \hat \w^{n,s}_{ijk,r}\,, 
\end{equation}
the reconstruction equations read 
\begin{equation}
 \frac{1}{\Delta x_e} \int \limits_{x_{e-\halb}}^{x_{e+\halb}} \w_h^x(x,t^n) dx = \frac{1}{\Delta x_e} \int \limits_{x_{e-\halb}}^{x_{e+\halb}} \psi_r(\xi(x)) \hat \w^{n,s}_{ijk,r} \, dx 
= {\bar{\Q}}^n_{ejk}, \qquad \forall {I}_{ejk} \in \mathcal{S}_{ijk}^{s,x}\,.
 \label{eqn.rec.x} 
\end{equation}
Equations~(\ref{eqn.rec.x}) provide a system of $M+1$ linear equations for the
unknown coefficients $\hat \w^{n,s}_{ijk,r}$, which is conveniently solved through linear algebra packages.
Once this operation has been performed for each stencil,
we  construct a data-dependent nonlinear combination of the resulting polynomials, i.e. 
\begin{equation}
\label{eqn.weno} 
 \w_h^x(x,t^n) = \psi_r(\xi) \hat \w^{n}_{ijk,r}, \quad \textnormal{ with } \quad  
 \hat \w^{n}_{ijk,r} = \sum_{s=1}^{N_s} \omega_s \hat \w^{n,s}_{ijk,r}\,.
\end{equation}   
The nonlinear weights $\omega_s$ are computed according to the WENO approach \cite{shu_efficient_weno} and their explicit expression can be found in 
\cite{AMR3DCL,AMR3DNC,Zanotti2015}. The whole procedure must be repeated along the two directions $y$ and $z$. 
Hence, although  each direction is treated separately, the net effect provides a genuine multidimensional reconstruction. 
We now proceed with the \textbf{key step} of the new algorithm presented in this paper and compute the \textit{point values} of the 
conserved quantities at the cell centers, simply by \textit{evaluating} the reconstruction polynomials in the barycenter of each control volume: 
\begin{equation}
  \Q_{ijk}^n = \w_h \left( x_i,y_j,z_k,t^n \right). 
	\label{eqn.pointeval} 
\end{equation} 
These point values of the conserved quantities $\Q_{ijk}^n$ are now converted into point values of the primitive variables $\V_{ijk}^n$, 
which requires only a single {\em cons-to-prim conversion} per cell. In RHD and RMHD, this is one of the most expensive 
and most delicate parts of the entire algorithm:  
\begin{equation}
  \V_{ijk}^n = \V \left( \Q_{ijk}^n \right). 
	\label{eqn.cons2prim} 
\end{equation} 

The reconstruction polynomials in primitive variables are spanned by the same basis functions $\psi_r(\xi)$ used for $\w_h$, hence 
\begin{equation}
\label{eqn.recpolydefprim.x} 
 \p^{s,x}_h(x,t^n) = \sum \limits_{r=0}^M \psi_r(\xi) \hat \p^{n,s}_{ijk,r} := \psi_r(\xi) \hat \p^{n,s}_{ijk,r}\,, 
\end{equation}
According to step $3$ listed above, we now require that the reconstructed polynomial must interpolate the 
\textit{point-values} of the \textit{primitive variables} at the centers of the cells forming each stencil, i.e.  
\begin{equation}
 \p_h^x(x_e,t^n) = \psi_r(\xi(x_e)) \hat \p^{n,s}_{ijk,r} = \V_{ejk}^n\,, \qquad \forall {I}_{ejk} \in \mathcal{S}_{ijk}^{s,x}. 
 \label{eqn.recprim.x} 
\end{equation}

The reconstruction equations~(\ref{eqn.recprim.x}) will also generate a system of $M+1$ linear equations for the
unknown coefficients $\hat \p^{n,s}_{ijk,r}$. The rest of the WENO logic applies in the same way, leading to 
\begin{equation}
\label{eqn.weno.prim} 
 \p_h^x(x,t^n) = \psi_r(\xi) \hat \p^{n}_{ijk,r}, \quad \textnormal{ with } \quad  
 \hat \p^{n}_{ijk,r} = \sum_{s=1}^{N_s} \omega_s \hat \p^{n,s}_{ijk,r}\,.
\end{equation}   
We emphasize that thanks to our polynomial WENO reconstruction (instead of the original point-wise WENO reconstruction 
of Jiang and Shu \cite{shu_efficient_weno}), the point-value of $\w_h(x,y,z,t^n)$ at each cell center, 
which is required at step $2$ above, is promptly available after evaluating the basis functions at the cell center. 
In other words, there is no need to perform any special transformation from cell averages to point-values via Taylor series 
expansions, like in \cite{Buchmuller2014,Buchmuller2015}. On the other hand, since the WENO reconstruction is performed twice, once for the 
conserved variables and once for the primitive variables, we expect that our new approach will become convenient in terms of computational 
efficiency  only for those systems of equations characterized by relations $\V(\Q)$ that cannot be written in 
closed form. 
In such circumstances, in fact, reducing the number of {\em cons-to-prim conversions} from $M (M+1)^{d+1} + d (M+1)^d$ in $d$ space dimensions 
(due to the space-time predictor and the numerical flux computation in the finite volume scheme) to just \textit{one single conversion} per cell 
will compensate for the double WENO reconstruction in space that we must perform. 
On the contrary, for systems of equations, such as the compressible Euler, for which the {\em cons-to-prim conversion} is analytic, 
no benefit will be reported in terms of computational efficiency, but still a significant benefit will be reported in terms of 
numerical accuracy. All these comments will be made quantitative in Sect.~\ref{sec:num-tests}.

\subsection{A local space--time DG predictor in primitive variables}
\label{sec:localDG}
%
\subsubsection{Description of the predictor}
\label{sec:Description_of_the_predictor}
As already remarked, the computation of the fluxes through the integrals
(\ref{averF}--\ref{averH}) is more
conveniently performed if the primitive variables are available 
at each space-time quadrature point. In such a case, in fact, no conversion from the conserved
to the primitive variables is required. According  to the discussion of the previous Section, it is possible to obtain
a polynomial $\p_h(x,y,z,t^n)$ in primitive variables at the reference time $t^n$. This is however not enough for a high accurate computation of the 
numerical fluxes, and $\p_h(x,y,z,t^n)$ must be evolved in time, locally for each cell, in order to obtain a polynomial $\v_h(x,y,z,t)$
approximating the solution at any time in the range $[t^n;t^{n+1}]$.

To this extent, we need an operation, to be performed locally for each
cell, which uses as input the high order polynomial
$\v_h$ obtained from the WENO reconstruction, and gives
as output its evolution in time, namely
\begin{equation}
\label{LSDG}
\p_h(x,y,z,t^n)\xrightarrow{LSDG} \v_h(x,y,z,t)\,,\hspace{1cm}t\in[t^n;t^{n+1}]\,.
\end{equation} 
This can be obtained through an 
element--local space--time Discontinuous Galerkin predictor that 
is based on the \textit{weak} integral form of
Eq.~(\ref{Primsyst}). From a mathematical point of view, Eq.~(\ref{Primsyst}) is a hyperbolic system in non-conservative form.
Therefore, the implementation of the space--time Discontinuous Galerkin predictor follows strictly the strategy 
already outlined in \cite{AMR3DNC} for non-conservative systems. 
Here we recall briefly the main ideas, focusing on the novel aspects implied by the formulation of 
Eq.~(\ref{Primsyst}).
The sought polynomial $\v_h(x,y,z,t)$ is supposed to be expanded in space and time as
\begin{equation}
 \v_h = \v_h(\boldsymbol{\xi},\tau) = \theta_l\left(\boldsymbol{\xi},\tau \right) \hat \v^n_l\,,   
 \label{eqn.st.q} 
\end{equation}
where the degrees of freedom $\hat \v^n_l$ are the unknowns. The space-time basis functions $\theta_l$ are given by a dyadic product of the 
Lagrange interpolation polynomials that pass through the Gauss-Legendre quadrature points, i.e. the tensor-product quadrature points on 
the hypercube $[0,1]^{d+1}$, see \cite{stroud}. 
The system (\ref{Primsyst}) is first rephrased in terms of the reference coordinates $\tau$ and $\boldsymbol{\xi} = (\xi,\eta,\zeta)$, yielding
\begin{equation}
\label{NCsyst_ref}
\frac{\partial{\V}}{\partial \tau} 
+ \mathbf{C}_1^\ast \frac{\partial{\V}}{\partial \xi} + \mathbf{C}_2^\ast \frac{\partial{\V}}{\partial \eta} + \mathbf{C}_3^\ast \frac{\partial{\V}}{\partial \zeta}
 ={\bf S}^\ast \,,
\end{equation}
with
\begin{equation}
{\bf C}_1^\ast= \frac{\Delta t}{\Delta x_i} \, {\bf C}_1, \quad 
{\bf C}_2^\ast= \frac{\Delta t}{\Delta y_j} \, {\bf C}_2, \quad 
{\bf C}_3^\ast= \frac{\Delta t}{\Delta z_k} \, {\bf C}_3, \quad 
{\bf S}^\ast= \Delta t \M^{-1}{\bf S}. 
\end{equation}
Expression (\ref{NCsyst_ref}) is then multiplied by the piecewise space-time 
polynomials $\theta_k(\xi,\eta,\zeta,\tau)$ and integrated over the 
space-time reference control volume, thus providing 
\begin{eqnarray}
\int \limits_{0}^{1} \int \limits_{0}^{1}  \int \limits_{0}^{1}   \int \limits_{0}^{1}   
\theta_k 
  \frac{\partial{\v_h}}{\partial \tau}  d \boldsymbol{\xi}  d\tau = 
\int \limits_{0}^{1} \int \limits_{0}^{1}  \int \limits_{0}^{1}   \int \limits_{0}^{1}   
\theta_k \left( {\bf S}^\ast - \mathbf{C}_1^\ast \frac{\partial \v_h}{\partial \xi} - \mathbf{C}_2^\ast \frac{\partial \v_h}{\partial \eta} - \mathbf{C}_3^\ast \frac{\partial \v_h}{\partial \zeta} \right) d \boldsymbol{\xi}  d\tau\,,
\label{eqn.pde.weak1} 
\end{eqnarray}
where we have replaced $\V$ with its discrete representation $\v_h$.
Integrating the first term by parts in time yields 
\begin{eqnarray}
 && \int \limits_{0}^{1} \int \limits_{0}^{1}  \int \limits_{0}^{1} \theta_k(\boldsymbol{\xi},1) \v_h(\boldsymbol{\xi},1) 
  \, d \boldsymbol{\xi}  
 - \int \limits_{0}^{1} \int \limits_{0}^{1}  \int \limits_{0}^{1}   \int \limits_{0}^{1} \left(\frac{\partial}{\partial \tau} \theta_k \right) \v_h(\boldsymbol{\xi},\tau) 
 \, d \boldsymbol{\xi} d\tau  = 
   \nonumber \\ 
&&  \int \limits_{0}^{1} \int \limits_{0}^{1}  \int \limits_{0}^{1} \theta_k(\boldsymbol{\xi},0) \p_h(\boldsymbol{\xi},t^n) \, d \boldsymbol{\xi} \, 
+ \int \limits_{0}^{1} \int \limits_{0}^{1}  \int \limits_{0}^{1}   \int \limits_{0}^{1}   
\theta_k \left( {\bf S}^\ast - \mathbf{C}_1^\ast \frac{\partial \v_h}{\partial \xi} - \mathbf{C}_2^\ast \frac{\partial \v_h}{\partial \eta} - \mathbf{C}_3^\ast \frac{\partial \v_h}{\partial \zeta} \right) d \boldsymbol{\xi}  d\tau\,. \nonumber \\
&&
\label{eqn.pde.weak3} 
\end{eqnarray}
Eq.~(\ref{eqn.pde.weak3}) is an element-local nonlinear algebraic equation that must be solved locally for each grid-cell in the unknowns 
$\hat \v^n_l$. In practice, we solve the system of Eqs.~(\ref{eqn.pde.weak3}) through a discrete Picard iteration, see \cite{DumbserZanotti,HidalgoDumbser}, 
where additional comments about its solution can be found.

\subsubsection{An efficient initial guess for the predictor}
\label{sec:initial.guess}
A proper choice of the initial guess for each of the space-time degrees of freedom $\hat \v_l$ can improve the convergence of the Picard process. 
The easiest strategy is to set $\v_h(\mathbf{x},t) = \p_h(\mathbf{x},t^n)$ i.e. the reconstruction polynomial is simply extended as a constant in time. 
This is, however, not the best approach. A better strategy for obtaining a good initial guess for the LSDG predictor was presented in 
\cite{HidalgoDumbser}, and it is based on the implementation of a MUSCL scheme for the explicit terms, plus a second-order Crank--Nicholson scheme in 
case stiff source terms are present. 
In the following, we refer to this version of the initial guess for the LSDG predictor as the MUSCL-CN initial guess.
If the source terms are not stiff, however, an even more efficient approach is possible which is based on a space-time extension of multi-level Adams--Bashforth-type
ODE integrators. For that purpose, the space-time polynomial denoted by $\v^{n-1}_h(\mathbf{x},t)$ obtained during the previous time step $[t^{n-1},t^n]$ is simply 
\textit{extrapolated in time} to the new time step $[t^n,t^{n+1}]$ by simple L2 projection: 
\begin{equation}
  \int \limits_{I_{ijk}} \int \limits_{t^n}^{t^{n+1}} \theta_k(\mathbf{x},t) \v^n_h(\mathbf{x},t)     \, dt \, d\mathbf{x} = 
	\int \limits_{I_{ijk}} \int \limits_{t^n}^{t^{n+1}} \theta_k(\mathbf{x},t) \v^{n-1}_h(\mathbf{x},t) \, dt \, d\mathbf{x}.  
\end{equation} 
In terms of the degrees of freedom $\hat \v^n_l$ and $\hat \v^{n-1}_l$ this relation becomes 
\begin{equation}
  \int \limits_0^1 \int \limits_0^1 \int \limits_0^1 \int \limits_0^1 \theta_k(\boldsymbol{\xi},\tau) \theta_l(\boldsymbol{\xi},\tau) \hat \v^n_l     \, dt \, d\boldsymbol{\xi} = 
	\int \limits_0^1 \int \limits_0^1 \int \limits_0^1 \int \limits_0^1 \theta_k(\boldsymbol{\xi},\tau) \theta_l(\boldsymbol{\xi},\tau') \hat \v^{n-1}_l \, dt \, d\boldsymbol{\xi},  
\label{eqn.abig} 
\end{equation} 
with $\tau' = 1 + \tau \frac{\Delta t^n}{\Delta t^{n-1}}$ and $\Delta t^{n-1} = t^n - t^{n-1}$.  

In the following, we refer to this second version of the initial guess for the LSDG predictor as the Adams--Bashforth (AB) initial guess. 
In Tab.~\ref{tab.CPU.initial.guess} we show a comparison among the performances of the LSDG predictor with these two different implementations of the initial guess.

\section{Numerical tests with the new ADER-WENO finite volume scheme in primitive variables}
\label{sec:num-tests}

In the following we explore the properties of the new ADER-WENO finite volume scheme by solving a wide set of test problems belonging to
four different systems of equations: the classical Euler equations, the relativistic hydrodynamics (RHD) and magnetohydrodynamics (RMHD) equations and
the Baer-Nunziato equations for compressible two-phase flows. For the sake of clarity, we introduce the notation
\enquote{ADER-Prim} to refer to the novel approach of this work
for which both the spatial WENO reconstruction and the subsequent
LSDG predictor are performed on the primitive variables. On the contrary, we denote  
the traditional ADER implementation, for which
both the spatial WENO reconstruction and the LSDG predictor are performed on the conserved variables,
as \enquote{ADER-Cons}. In a few circumstances, we have also compared with the \enquote{ADER-Char} scheme, namely a traditional ADER scheme in which, however, the spatial reconstruction is performed
on the characteristic variables. 
In this Section we focus our attention on finite volume schemes, which,
according to the notation introduced in \cite{Dumbser2008}, 
are denoted as $\mathbb{P}_0\mathbb{P}_M$ methods, where $M$ is the degree of the approximating polynomial.
In Sect.~\ref{sec:extension} a brief account is given to Discontinuous Galerkin methods, referred to as $\mathbb{P}_N\mathbb{P}_N$ methods, for which
an ADER-Prim version is also possible.

\subsection{Euler equations}
\label{sec:Euler}

First of all we consider the solution of the classical Euler equations 
of compressible gas dynamics, for which the vectors of the conserved variables 
$\Q$ and of the fluxes ${\bf f}^x$, ${\bf f}^y$ and ${\bf f}^z$ are given respectively by
\begin{equation}
{\Q}=\left(\begin{array}{c}
\rho \\ \rho v_x \\ \rho v_y \\ \rho v_z \\ E 
\end{array}\right) \, ,  \,\,\,\, 
{\bf f}^x=
\left(\begin{array}{c}
\rho v_x \\
\rho v_x^2 + p \\
\rho v_xv_y \\
\rho v_xv_z \\
v_x(E+p)
\end{array}\right)\, ,  \,\,\,\, 
{\bf f}^y=
\left(\begin{array}{c}
\rho v_y \\
\rho v_xv_y \\
\rho v_y^2 + p \\
\rho v_yv_z \\
v_y(E+p)
\end{array}\right)\, ,  \,\,\,\, 
{\bf f}^z=
\left(\begin{array}{c}
\rho v_z \\
\rho v_xv_z \\
\rho v_yv_z \\
\rho v_z^2 + p \\
v_z(E+p)
\end{array}\right)\,.
\label{eq:Euler-system}
\end{equation}
Here $v_x$, $v_y$ and $v_z$ are the velocity components, $p$ is the pressure, $\rho$ is the mass density,
$E=p/(\gamma-1)+\rho (v_x^2+v_y^2+v_z^2)/2$ is the total energy density, while $\gamma$ is the adiabatic index of the supposed ideal gas equation of state, which
is of the kind $p=\rho\epsilon(\gamma-1)$, $\epsilon$ being the specific internal energy. 

\subsubsection{2D isentropic vortex}
\label{sec:isentropic}
 \begin{table}[!t] 
 \centering
 \begin{tabular}{|c|c||cc|cc|cc|c|}
   \hline
   \multicolumn{9}{|c|}{\textbf{2D isentropic vortex problem }} \\
   \hline
	\hline
    &       &  \multicolumn{2}{c|}{ ADER-Prim } &  \multicolumn{2}{c|}{ ADER-Cons }  & \multicolumn{2}{c|}{ ADER-Char } &    \\
   \hline
    & $N_x$ &  $L_2$ error &   $L_2$ order &  $L_2$ error &   $L_2$ order  & $L_2$ error &   $L_2$ order &   Theor. \\
   \hline
   \hline
   \multirow{5}{*}{\rotatebox{0}{{$\mathbb{P}_0\mathbb{P}_2$}}}
& 100	&  4.060E-03	& ---	   & 5.028E-03	& ---	   & 5.010E-03	& ---	&\multirow{5}{*}{3}\\
& 120	&  2.359E-03	&  2.98	 & 2.974E-03	&  2.88	 & 2.968E-03	&  2.87	 &\\
& 140	&  1.489E-03	&  2.98	 & 1.897E-03	&  2.92	 & 1.893E-03	&  2.92	 &\\
& 160	&  9.985E-04	&  2.99  & 1.281E-03	&  2.94	 & 1.279E-03	&  2.94	 &\\
& 200	&  5.118E-04	&  2.99  & 6.612E-04	&  2.96	 & 6.607E-04	&  2.96	 &\\		
   \hline
   \multirow{5}{*}{\rotatebox{0}{{$\mathbb{P}_0\mathbb{P}_3$}}}
& 50	&  2.173E-03	&  ---	 & 4.427E-03	& ---	   & 5.217E-03	&  ---	& \multirow{5}{*}{4}\\
& 60	&  8.831E-04	&  4.93	 & 1.721E-03	&  5.18	 & 2.232E-03	&  4.65	 &\\
& 70	&  4.177E-04	&  4.85  & 8.138E-04	&  4.85	 & 1.082E-03	&  4.69	 &\\
& 80	&  2.194E-04	&  4.82  & 4.418E-04	&  4.57	 & 5.746E-04	&  4.74	 &\\
& 100	&  7.537E-05	&  4.79  & 1.605E-04	&  4.53	 & 1.938E-04	&  4.87	 &\\
   \hline
	  \multirow{5}{*}{\rotatebox{0}{{$\mathbb{P}_0\mathbb{P}_4$}}}
& 50	&  2.165E-03	& ---	   & 3.438E-03	& ---	   & 3.416E-03	& ---	&\multirow{5}{*}{5}\\
& 60	&  6.944E-04	&  6.23	 & 1.507E-03	&  4.52	 & 1.559E-03	&  4.30	 &\\
& 70	&  3.292E-04	&  4.84  & 7.615E-04	&  4.43	 & 7.615E-04	&  4.65	 &\\
& 80	&  1.724E-04	&  4.84  & 4.149E-04	&  4.55	 & 4.148E-04	&  4.55	 &\\
& 100	&  5.884E-05	&  4.82  & 1.449E-04	&  4.71	 & 1.448E-04	&  4.72	 &\\
   \hline
  \end{tabular}
 \caption{ \label{tab:Vortex_Error}  $L_2$  errors of the mass density and corresponding convergence rates for the 
   2D isentropic vortex problem. A comparison is shown among the reconstruction in primitive variables (ADER-Prim), in conserved variables (ADER-Cons) and in characteristic variables (ADER-Char). The Osher-type numerical flux has been used.}
	\label{Table:convergence}
 \end{table}
%
It is important to assess the convergence properties of the new scheme, in particular comparing with the traditional ADER scheme in conserved and in characteristic variables. 
To this extent, we have studied
the two-dimensional isentropic vortex,  see e.g. \cite{HuShuTri}.  
The initial conditions are given by a uniform mean flow, to which
a perturbation is added, such that
\begin{equation}
\left( \rho,v_x,v_y,v_z,p \right) =(1+\delta\rho, 1+\delta v_x, 1+\delta v_y, 0, 1+\delta p)\,,
\end{equation} 
with 
\begin{equation}
\left(\begin{array}{c}
\delta \rho \\ \delta v_x \\ \delta v_y \\ \delta p 
\end{array}\right)
=
\left(\begin{array}{c}
(1+\delta T)^{1/(\gamma-1)}-1 \\
-(y-5)\epsilon/2\pi \exp{[0.5(1-r^2)]} \\
\phantom{-}(x-5)\epsilon/2\pi \exp{[0.5(1-r^2)]} \\
(1+\delta T)^{\gamma/(\gamma-1)}-1
\end{array}\right).~~~
\label{eq:pert}
\end{equation}
Whatever the perturbation $\delta T$ in the temperature is, it is easy to verify that there is not any variation in the specific entropy $s=p/\rho^\gamma$, and the flow
is advected smoothly and isentropically with velocity $v=(1,1,0)$. We have solved this test over the computational domain $\Omega=[0;10]\times[0;10]$, assuming
\begin{equation}
\delta T=-\frac{\epsilon^2(\gamma-1)}{8\gamma\pi^2}~\exp{(1-r^2)}\,,
\end{equation}
with $r^2=(x-5)^2+(y-5)^2$, vortex strength $\epsilon=5$ and adiabatic index $\gamma=1.4$. Table~\ref{Table:convergence} contains the results of
our calculation, in which we have compared the convergence properties of three different finite volume ADER schemes: ADER-Prim, ADER-Cons and ADER-Char, obtained with the Osher-type Riemann solver, 
see \cite{OsherUniversal}.
While all the schemes converge to the nominal order, 
it is interesting to note that the smallest $L_2$ error is obtained for the \textit{new} ADER finite volume scheme in \textit{primitive variables}, 
and that the difference with respect to the other two reconstructions increases with the order of the method.  

In addition to the convergence properties, we have compared the performances of the Adams--Bashforth version of the initial guess for the LSDG predictor with the traditional version based on the MUSCL-CN algorithm. The comparison has been performed over a $100\times 100$ uniform grid.
The results are shown in Tab.~\ref{tab.CPU.initial.guess}, from which we conclude that the Adams--Bashforth initial guess is indeed computationally more
efficient in terms of CPU time. However, we have also experienced that it is typically less robust, and in some of the most challenging numerical tests discussed in the rest of the paper we had to use the more traditional MUSCL-CN initial guess.
\begin{table}[!t]
\vspace{0.5cm}
\renewcommand{\arraystretch}{1.0}
\begin{center}
\begin{tabular}{ccc}
\hline
\hline
             & \footnotesize{MUSCL-CN}   & \footnotesize{Adams--Bashforth}   \\
\hline
\hline
$\mathbb{P}_0\mathbb{P}_2$         & 1.0                                  & 0.64                             \\
\hline
$\mathbb{P}_0\mathbb{P}_3$         & 1.0                                  & 0.75                             \\
\hline
$\mathbb{P}_0\mathbb{P}_4$         & 1.0                                  & 0.72                             \\
\hline
\hline
\end{tabular}
\end{center}
\caption{CPU time comparison among different versions of the initial guesses for the LSDG predictor. The comparison has been performed for the isentropic vortex solution and 
the numbers have been normalized to the value obtained with the traditional MUSCL-CN initial guess (see Sect.~\ref{sec:initial.guess} for more details).} 
\label{tab.CPU.initial.guess}
\end{table}
%
%
\subsubsection{Sod's Riemann problem}
%
\begin{table}[!b]
\vspace{0.5cm}
\renewcommand{\arraystretch}{1.0}
\begin{center}
\begin{tabular}{cccc}
\hline
\hline
             & \footnotesize{ADER-Prim}   & \footnotesize{ADER-Cons}   &  \footnotesize{ADER-Char} \\
\hline
$\mathbb{P}_0\mathbb{P}_2$         & 1.0                                  & 0.74                                & 0.81  \\
\hline
$\mathbb{P}_0\mathbb{P}_3$         & 1.0                                  & 0.74                                & 0.80  \\
\hline
$\mathbb{P}_0\mathbb{P}_4$         & 1.0                                  & 0.77                                & 0.81  \\
\hline
\end{tabular}
\end{center}
\caption{CPU time comparison among different ADER implementations for the Sod Riemann problem. The 
numbers have been normalized to the value obtained with ADER-Prim.}
\label{tab.CPU.Sod}
\end{table}
\begin{figure}
\begin{center}
\begin{tabular}{cc} 
{\includegraphics[angle=0,width=7.3cm,height=7.3cm]{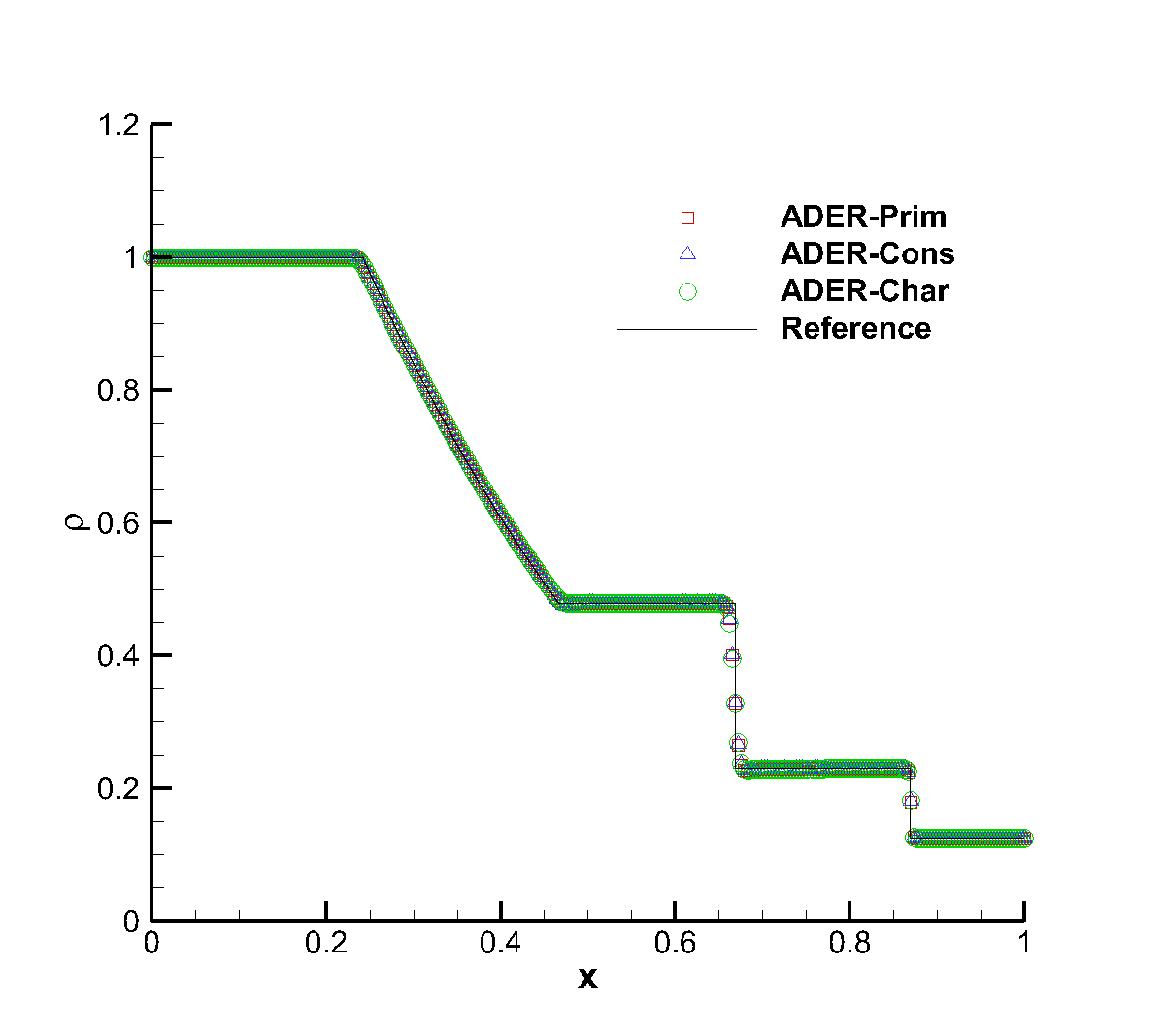}} & 
{\includegraphics[angle=0,width=7.3cm,height=7.3cm]{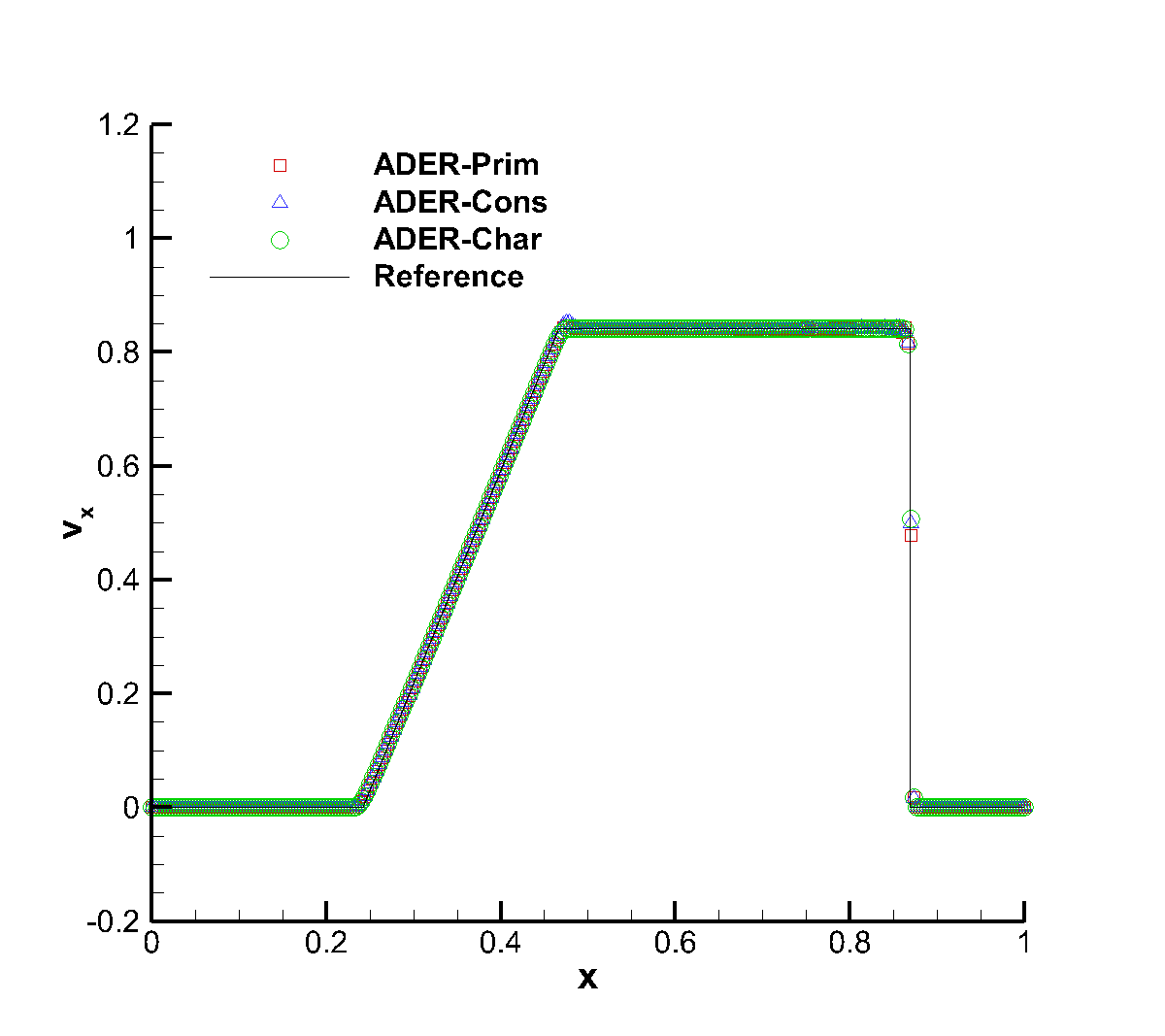}} \\
{\includegraphics[angle=0,width=7.3cm,height=7.3cm]{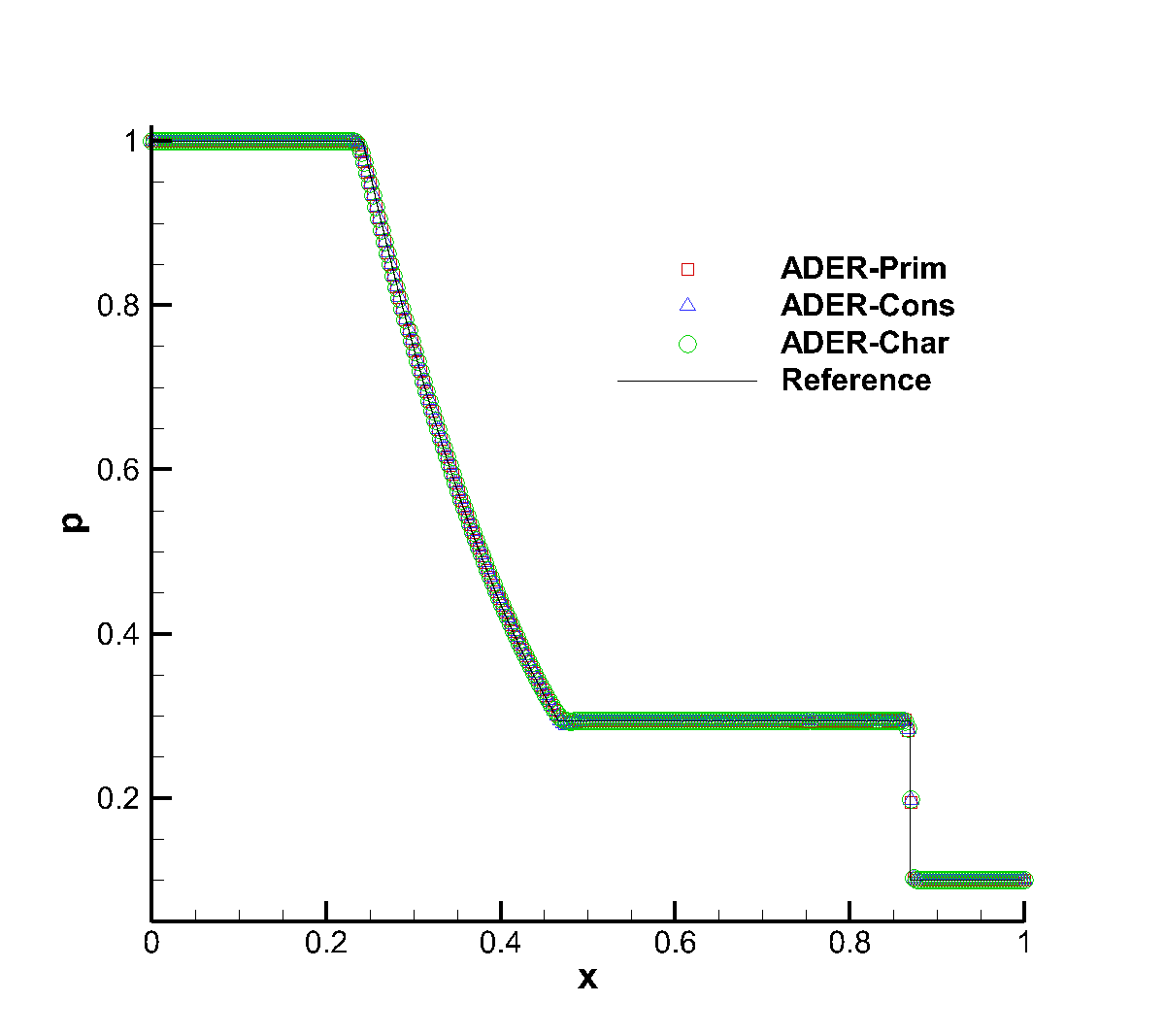}} & 
{\includegraphics[angle=0,width=7.3cm,height=7.3cm]{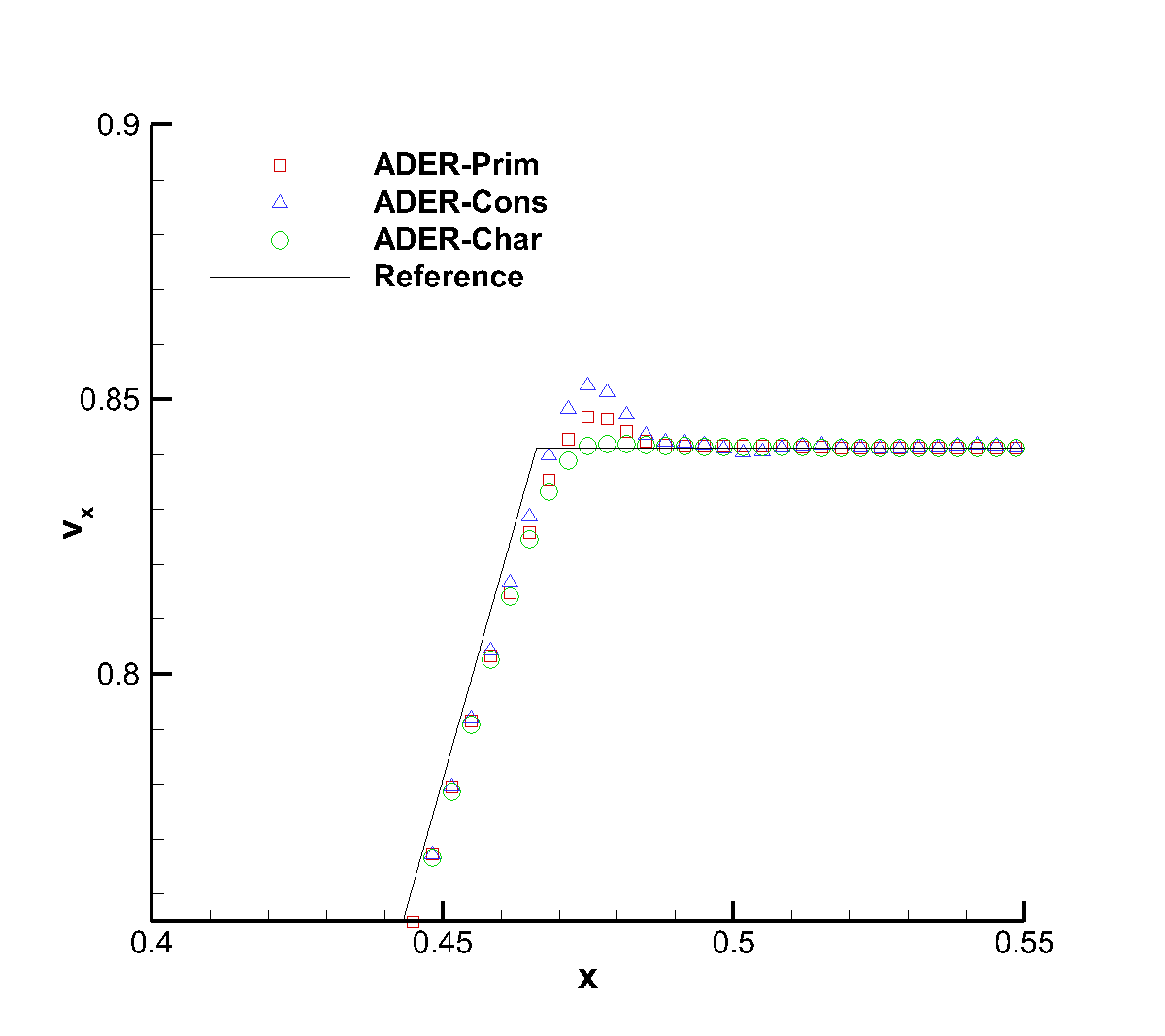}}
\end{tabular} 
\caption{Solution of Sod's Riemann problem with the fourth order ADER-WENO scheme at time $t=0.2$. 
The bottom right panel shows a magnification of the velocity at the tail of the rarefaction.} 
\label{fig:shock-tube-Sod}
\end{center}
\end{figure}
We have then solved the classical Riemann problem named after Sod~\cite{Sod1978}, assuming an adiabatic index $\gamma=1.4$, and evolved until $t_{\rm final}=0.2$.
In spite of the fact that this is a one-dimensional test, we have evolved this problem in two spatial dimensions over the domain $[0,1]\times[-0.2,0.2]$,
using periodic boundary conditions along the passive $y$ direction.
In Fig.~\ref{fig:shock-tube-Sod} we show the comparison among the solutions obtained with ADER-Prim, ADER-Cons and ADER-Char, together with the exact solution
provided in \cite{toro-book}. 
We have adopted the finite volume scheme at the fourth order of accuracy, namely the $\mathbb{P}_0\mathbb{P}_3$ scheme,
in combination with the Rusanov numerical flux and using $400$ cells along the $x$-direction.
Although all of the ADER implementations show a very good agreement with the exact solution, a closer look at the tail of the rarefaction,
highlighted in the bottom right panel, reveals that the ADER-Cons scheme is actually the worst one, while the solution obtained with ADER-Prim is more similar to 
the reconstruction in characteristic variables. 
On the contrary, in terms of CPU-time, ADER-Prim is not convenient for this system of equations 
because the price paid for performing the double WENO reconstruction in space is not significantly compensated by the reduced number of conversions that are needed 
from the conserved to the primitive variables. Table~\ref{tab.CPU.Sod} reports the CPU times, normalized with respect to the ADER-Prim implementation, for different 
orders of accuracy, showing that the ADER-Prim scheme is $\sim25\%$ slower than the traditional ADER-Cons scheme.
As we will see in Tab.~\ref{tab.CPU.RHD} of  Sect.~\ref{sec:RMHD}, the comparison will change in favor of ADER-Prim schemes, when the relativistic equations are solved instead.
%
%

%

\subsubsection{Interacting blast waves}
%
\begin{figure}
\begin{center}
\begin{tabular}{cc} 
{\includegraphics[angle=0,width=6.5cm,height=6.5cm]{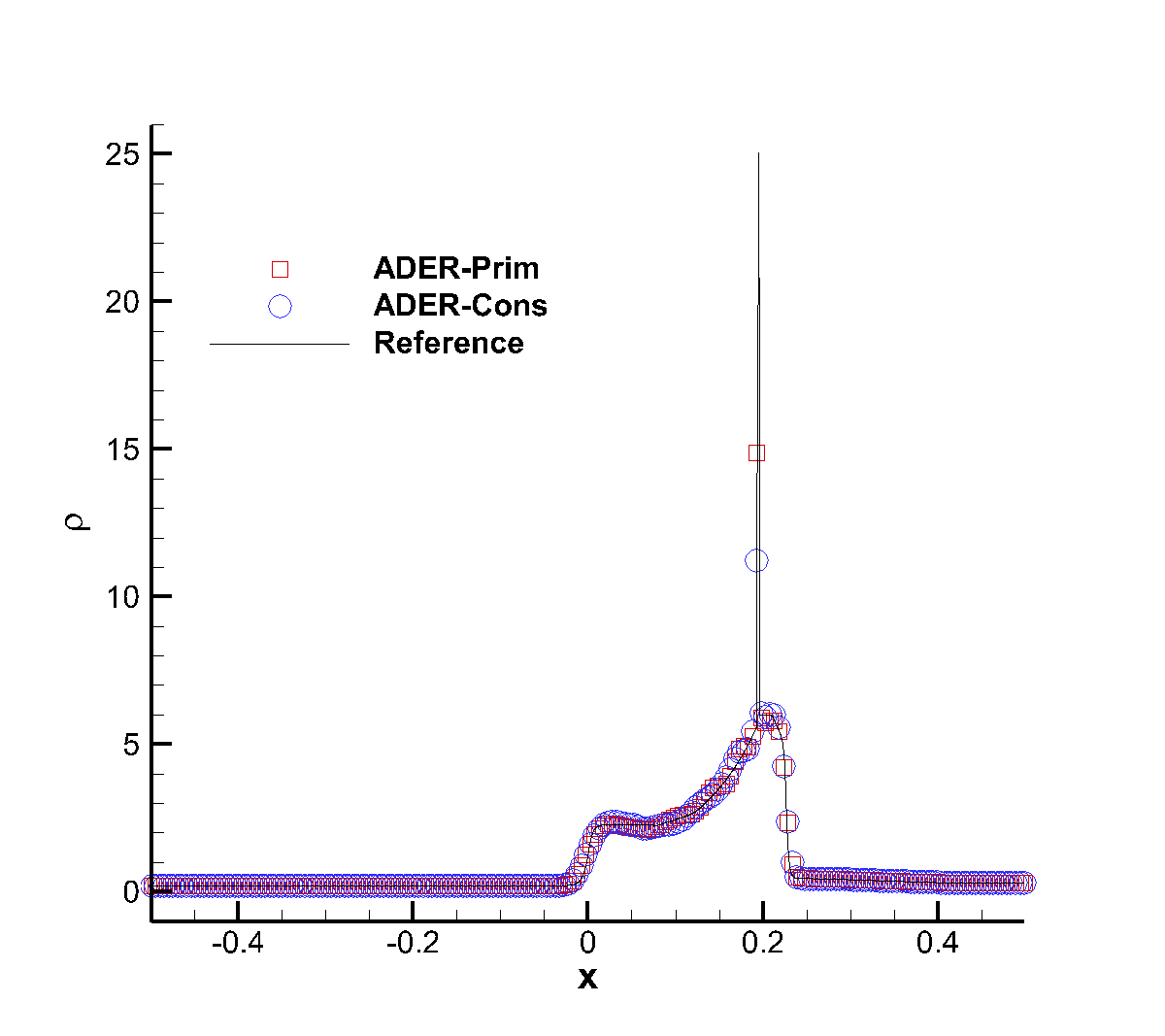}} & 
{\includegraphics[angle=0,width=6.5cm,height=6.5cm]{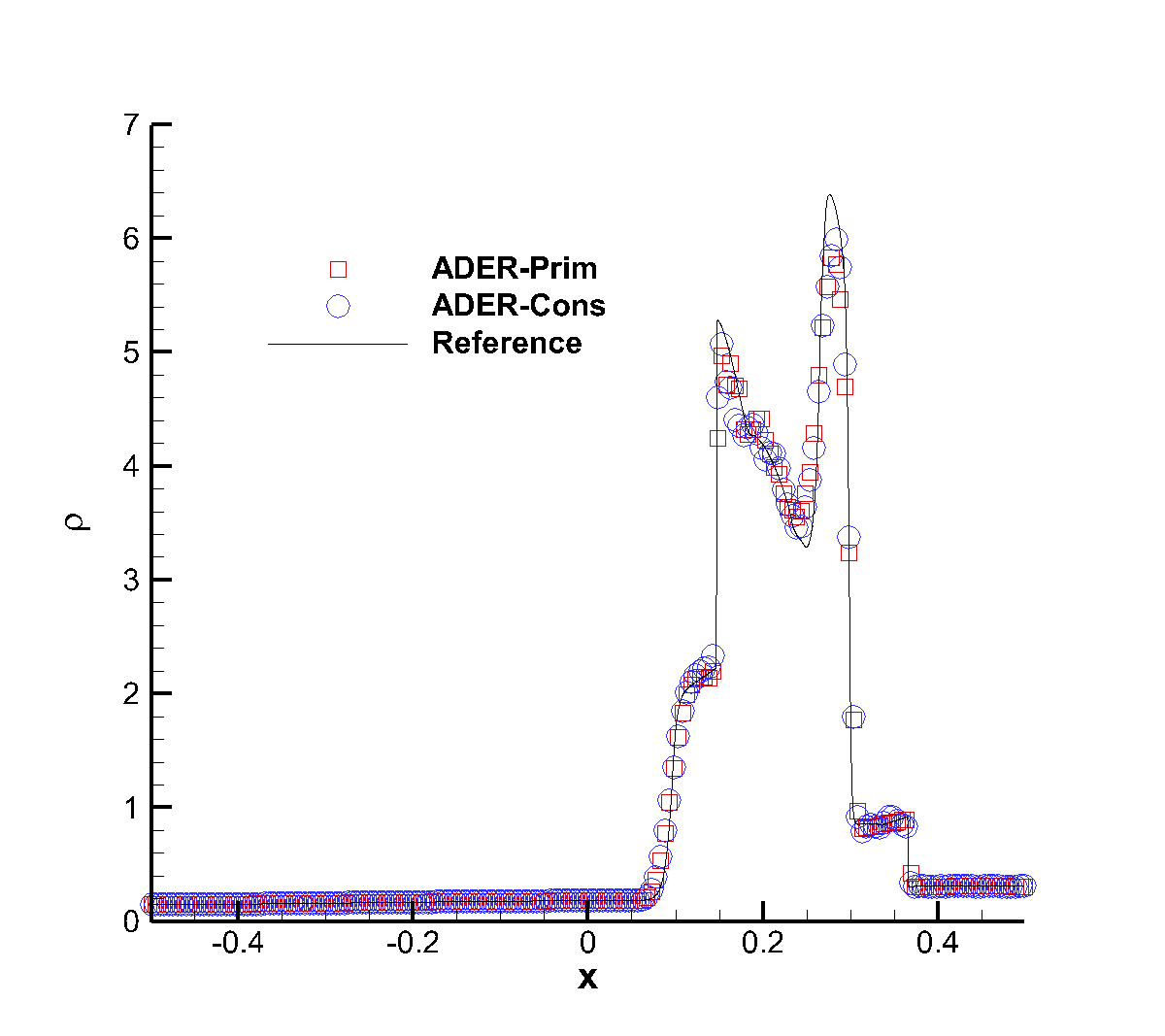}}
\end{tabular} 
\caption{Solution of the interacting Blast-Wave problem at time $t=0.028$ (left panel) and at time $t=0.038$ (right panel)
obtained with the fourth order ADER-WENO scheme. The computation has been performed over a uniform grid of 500 cells.}
\label{fig:BlastWavesEuler}
\end{center}
\end{figure}
%
The interaction between two blast waves was first proposed by \cite{woodwardcol84} and it is now a standard test for
computational fluid dynamics. The initial conditions are given by
\begin{equation}
\label{blast-wave}
(\rho,v_x,p)= \left\{
\begin{array}{llll}
(1.0,0.0,10^3) &   {\rm if} & -0.5 < x < -0.4 \,, \\
 
(1.0,0.0,10^{-2}) &   {\rm if} & -0.4 < x < 0.4 \,, \\
 
(1.0,0.0,10^2) &   {\rm if} & \phantom{-} 0.4 < x < 0.5 \,,
\end{array} \right.
\end{equation}
where the adiabatic index is $\gamma=1.4$. 
We have evolved this problem in two spatial dimensions over the domain $[-0.6,0.6]\times[-0.5,0.5]$,
using reflecting boundary conditions in $x$ direction and periodic boundary conditions along the $y$ direction. The results of our calculations,
obtained with the $\mathbb{P}_0\mathbb{P}_3$ scheme, are reported in Fig.~\ref{fig:BlastWavesEuler}, where only the one-dimensional cuts are shown.
The number of cells chosen along the x-direction, namely $N_x=500$, is not particularly large, at least for this kind of challenging problem. This has been intentionally done to better highlight 
potential differences among the two alternative ADER-Prim and ADER-Cons schemes. As it turns out from the figure, the two methods are very similar in terms of accuracy: the sharp peak 
in the density at time $t=0.028$ (left panel) is somewhat better resolved through the ADER-Prim, while the opposite is true for the highest peak at time $t=0.038$ (right panel).
On the overall, however, the two schemes perform equally well for this test. 

\begin{figure}
\begin{center}
{\includegraphics[angle=0,width=15.0cm,height=5.0cm]{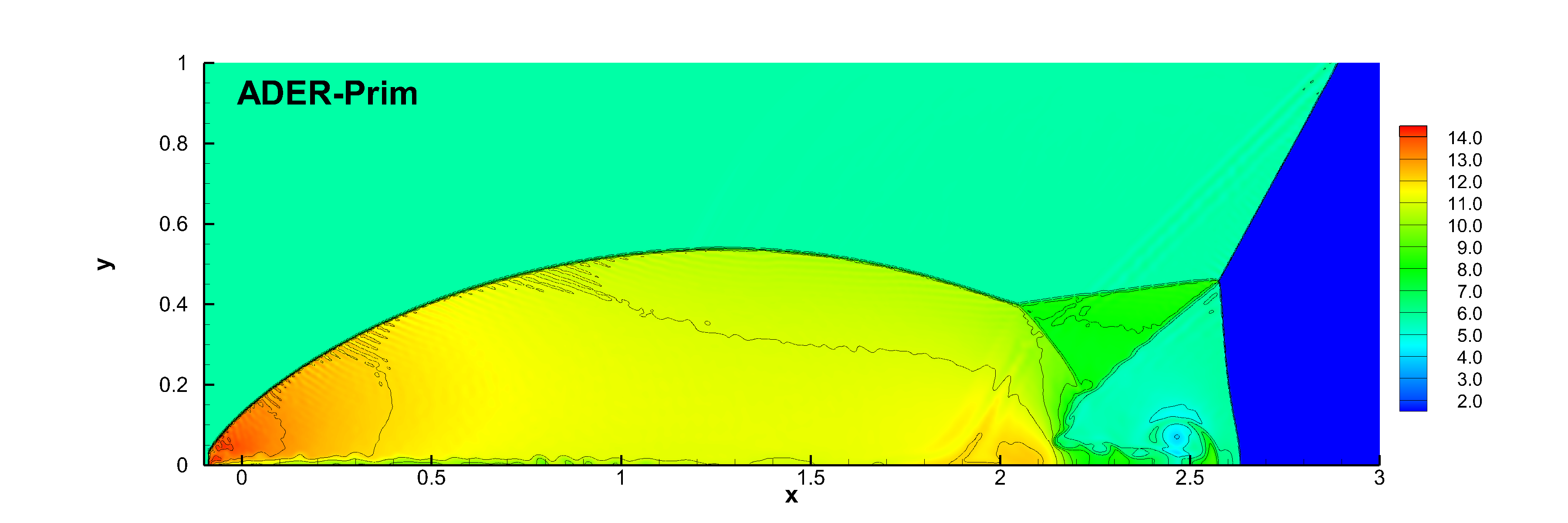}}
{\includegraphics[angle=0,width=15.0cm,height=5.0cm]{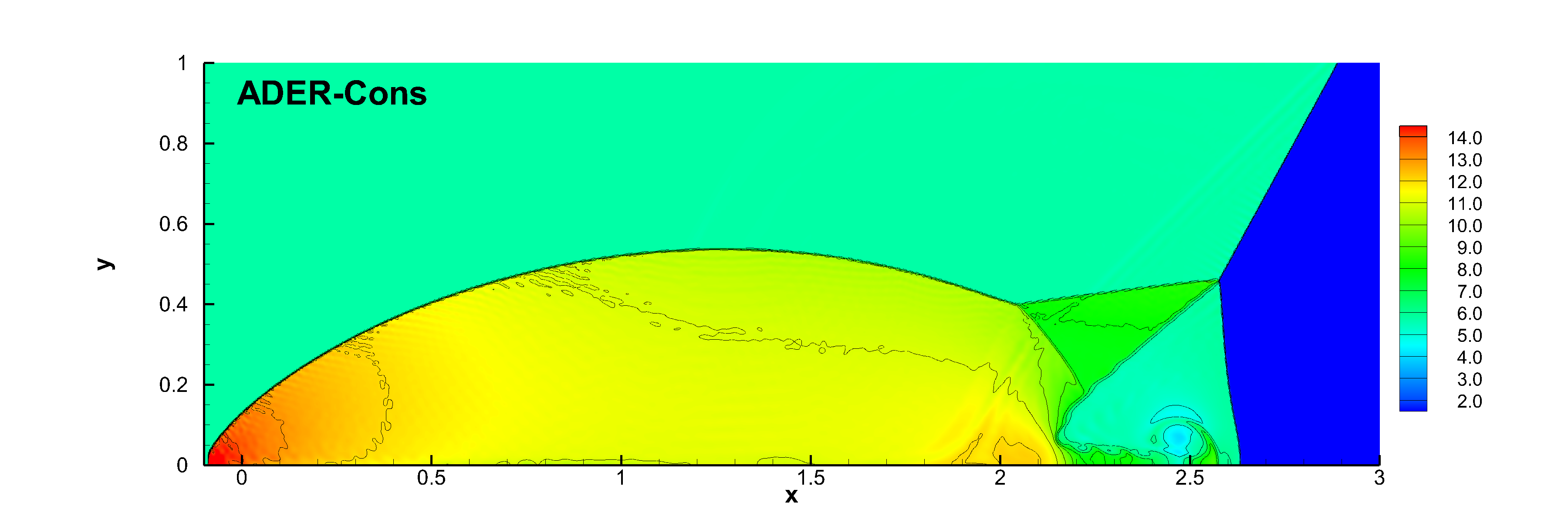}}
\caption{Double Mach reflection problem at time $t=0.2$ 
obtained with the fourth order ADER-WENO scheme and the Rusanov Riemann solver.
The computation has been performed over a uniform grid of $1200\times 300$ cells. 
Top panel: mass density distribution obtained with ADER-Prim.
Bottom panel: mass density distribution obtained with ADER-Cons.}
\label{fig:DMR}
\end{center}
\end{figure}
%

%
\subsubsection{Double Mach reflection problem} 

As a representative test for the Euler equations in two space dimensions, we have considered
the {\em double Mach reflection problem}, which implies the interaction of several waves. 
The dynamics of this problem is triggered by a shock wave propagating towards the right with a Mach number $M=10$, and intersecting
the $x-$ axis at  
$x=1/6$ with an inclination angle of $\alpha=60^{\circ}$. The initial  states ahead and behind the shock are fixed after  solving the Rankine--Hugoniot conditions, obtaining
\begin{eqnarray}
(\rho, u, v, p)( \x,t=0) =
\left\{
\begin{array}{cll}
  \frac{1}{\gamma}(8.0, 8.25, 0.0, 116.5), \quad & \text{ if } & \quad x'<0.1, \\
  (1.0, 0.0, 0.0, \frac{1}{\gamma}),       \quad & \text{ if } & \quad x'\geq 0.1, 
\end{array}
\right.
\end{eqnarray}
where $x' = (x - 1/6) \cos\alpha - y \sin\alpha$.   The adiabatic index is $\gamma=1.4$. 
We fix inflow and outflow boundary conditions 
on the left side and on the right of the numerical domain, respectively, 
while on the bottom we have used reflecting boundary conditions. At the top
we must impose the exact solution of an isolated moving oblique shock wave with the same shock Mach 
number $M_s=10$. 
We have solved the test 
over the  rectangle $\Omega = [0;3.0] \times [0;1]$, covered by a  uniform grid composed of
$1200\times300$ cells, using the Rusanov Riemann solver and a fourth order finite volume scheme.
The two panels of Fig.~\ref{fig:DMR} show the comparison of the solution at time $t=0.2$ obtained with the ADER-Prim (top panel) and with the ADER-Cons (bottom panel) scheme.
The results are very similar in the two cases.

As a tentative conclusion about the performances of ADER-Prim for the Euler equations, we may say that, although it is the most accurate on smooth solutions 
(see Tab.~\ref{Table:convergence}), and comparable to a traditional ADER with reconstruction in characteristic variables, it is computationally more expensive than ADER-Cons 
and ADER-Char. Hence, ADER-Prim will
rarely become the preferred choice in standard applications for the Euler equations.

%
%
\subsection{Relativistic hydrodynamics and magnetohydrodynamics}
\label{sec:RMHD}
From a formal point of view, the equations of special relativistic hydrodynamics  and magnetohydrodynamics can be written in conservative form like the 
classical Euler equations (see, however, the comments below), namely as in Eq.~(\ref{NCsyst}), with the vectors of the conserved variables and of the corresponding fluxes given by
\be
{\Q}=\left[\begin{array}{c}
D \\ S_j \\ U \\ B^j 
\end{array}\right],~~~
{\bf f}^i=\left[\begin{array}{c}
 v^i D \\
 W^i_j \\
 S^i \\
\epsilon^{jik}E^k 
\end{array}\right]\,,\hspace{1cm}i=x,y,z\,.
\label{eq:RMHDfluxes}
\ee
where the  conserved
variables $(D,S_j,U,B_j)$ can be expressed as\footnote{We note that, since the spacetime is flat and we are using Cartesian coordinates, the covariant and the contravariant components of spatial vectors can be used interchangeably, namely $A_i=A^i$, for the generic vector $\vec A$.}
\bea
\label{eq:cons1}
&&D   = \rho W ,\\
\label{eq:cons2}
&&S_i = \rho h W^2 v_i + \epsilon_{ijk}E_j B_k, \\
\label{eq:cons3}
&&U   = \rho h W^2 - p + \frac{1}{2}(E^2 + B^2)\,,
\eea
while the spatial projection of the energy-momentum
tensor of the fluid is \cite{DelZanna2007}  
\be
W_{ij} \equiv \rho h W^2 v_i v_j - E_i E_j - B_i B_j + \left[p +\frac{1}{2}(E^2+B^2)\right]\delta_{ij}\,.
\label{eq:W} 
\ee
Here $\epsilon_{ijk}$ is the Levi--Civita tensor and $\delta_{ij}$ is the Kronecker symbol.  
We have used the symbol $h=1+\epsilon+p/\rho$ to denote the specific enthalpy of the plasma and in all our calculations
the usual ideal gas equation of state has been assumed.
 
The components of the electric and of the magnetic field
in the laboratory frame are denoted by $E_i$ and $B_i$, while the Lorentz factor of the fluid
with respect to this reference frame is $W=(1-v^2)^{-1/2}$. 
We emphasize that the electric field
does not need to be evolved in time under the assumption of infinite electrical conductivity, since it can always be computed in terms of the velocity and of the magnetic field as 
$\vec E = - \vec v \times \vec B$.

Although formally very similar to the classical gas dynamics equations, their relativistic counterpart 
present two fundamental differences. The first one is that, while the physical fluxes ${\bf f}^i$ of the 
classical gas dynamics equations can be written analytically in terms of the conserved variables, i.e. ${\bf f}^i={\bf f}^i(\u)$, 
those of the relativistic hydrodynamics (or magnetohydrodynamics) equations  need the knowledge of the primitive variables, i.e. ${\bf f}^i={\bf f}^i(\V)$ for RMHD. 
The second difference is that, in the relativistic case,
the conversion from the conserved to the primitive variables, i.e. the operation
$(D,S_j,U,B_j)\longrightarrow (\rho,v_i,p,B_i)$, is not analytic, and it must be performed numerically through some appropriate iterative procedure.
Since in an ADER scheme such a conversion must be performed in each space-time degree of freedom of the space-time DG predictor and at each Gaussian quadrature 
point for the computation of the fluxes in the finite volume scheme, we may expect a significant computational  advantage by performing  
the WENO reconstruction and the LSDG predictor directly on the primitive variables. In this way, in fact,
the conversion $(D,S_j,U,B_j)\longrightarrow (\rho,v_i,p,B_i)$ is required only once at the cell center (see Sect.~\ref{sec:WENO_reconstruction}), and not in each 
space-time degree of freedom of the predictor and at each Gaussian point for the quadrature of the numerical fluxes.
We emphasize that the choice of the variables to reconstruct for the relativistic velocity is still a matter of debate. The velocity  $v_i$ may seem the most natural one, 
but, as first noticed by \cite{Komissarov1999}, reconstructing $W v_i$ can 
increase the robustness of the scheme. However, this is not always the case (see Sect.~\ref{sec:RMHD-Rotor-Problem} below)
and in our tests we have favored either the first or the second choice according to convenience.
Concerning the specific strategy adopted to recover the primitive variables, 
in our numerical code we have used the third method reported in Sect. 3.2 of 
\cite{DelZanna2007}. Alternative methods can be found in \cite{NGMD2006,Rezzolla_book:2013}.

Finally, there is an important formal change in the transition from purely hydrodynamics systems to genuinely magnetohydrodynamics systems. 
As already noticed by \cite{Londrillo2000},  the RMHD equations should not be regarded as a mere extension 
of the RHD ones, with just a larger number of variables to evolve. Rather, their formal structure is better described in terms of a coupled system of conservation laws (the five equations for the dynamics of the plasma) and a set of Hamilton-Jacobi equations, those for the evolution of the vector potential of the magnetic field \cite{ShiJin1998}.
The different mathematical structure of the RMHD equations reflects the existence of the divergence-free property of the magnetic field, which must be ensured at all times
during the evolution. Numerically, we have adopted a simplified and well known approach, which consists of 
augmenting the system (\ref{NCsyst}) with an additional equation for a scalar 
field $\Phi$, aimed at propagating away the deviations from $\vec \nabla\cdot\vec B=0$. We therefore need to solve 
\be
\label{eq:divB}
\partial_t \Phi + \partial_i B^i = -\kappa \Phi\,,
\ee
while the fluxes for the evolution of the magnetic field are also changed, namely ${\bf f}^i(B^j)\rightarrow \epsilon^{jik}E^k + \Phi \delta^{ij}$,
where $\kappa\in[1;10]$ in most of our calculations.
Originally introduced by \cite{Dedner:2002} for the classical MHD equations, this approach has been extended to the relativistic regime by \cite{Palenzuela:2008sf}.
More information about the mathematical structure of the RMHD equations can be found in \cite{Anile_book,BalsaraRMHD,Komissarov1999,DelZanna2007,Anton2010}.

In the following, we first  limit our attention to 
a few physical systems for which $B_i=E_i=0$, hence to relativistic hydrodynamics, and then we 
consider truly magnetohydrodynamics tests with $B_i\neq0$.

\subsubsection{RHD Riemann Problems}
\begin{figure}
\begin{center}
\begin{tabular}{cc} 
{\includegraphics[angle=0,width=7.3cm,height=7.3cm]{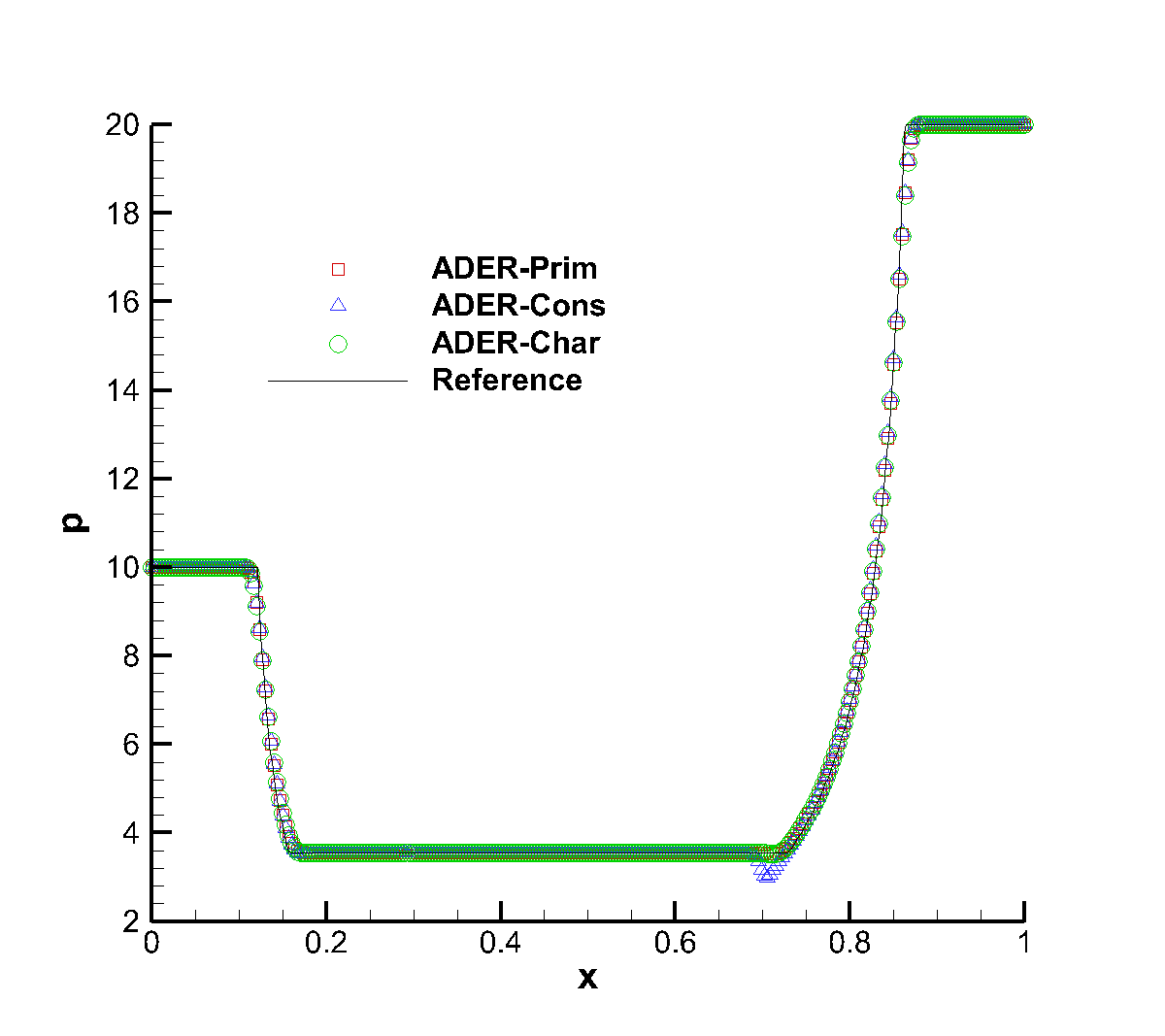}} & 
{\includegraphics[angle=0,width=7.3cm,height=7.3cm]{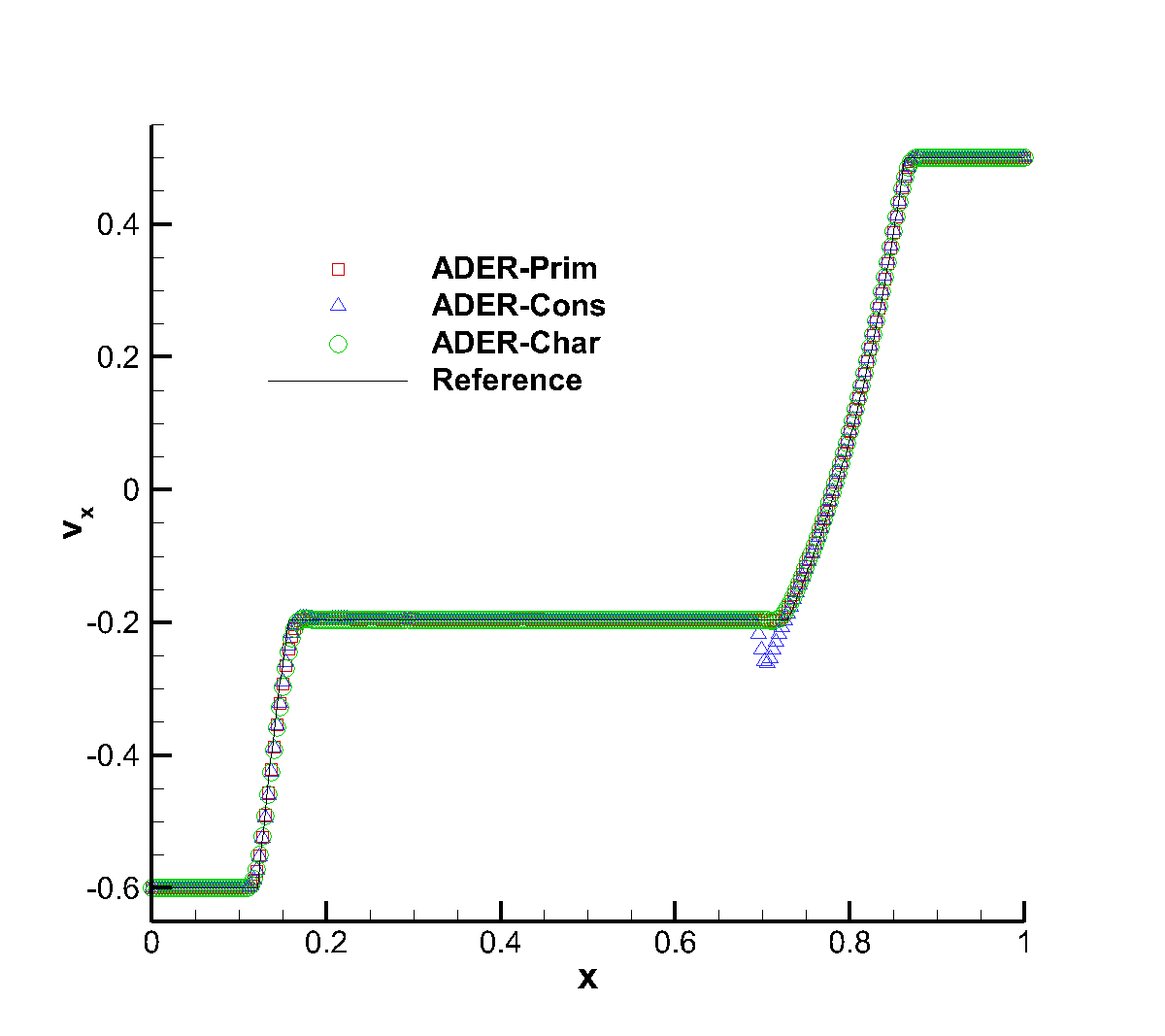}} \\ 
{\includegraphics[angle=0,width=7.3cm,height=7.3cm]{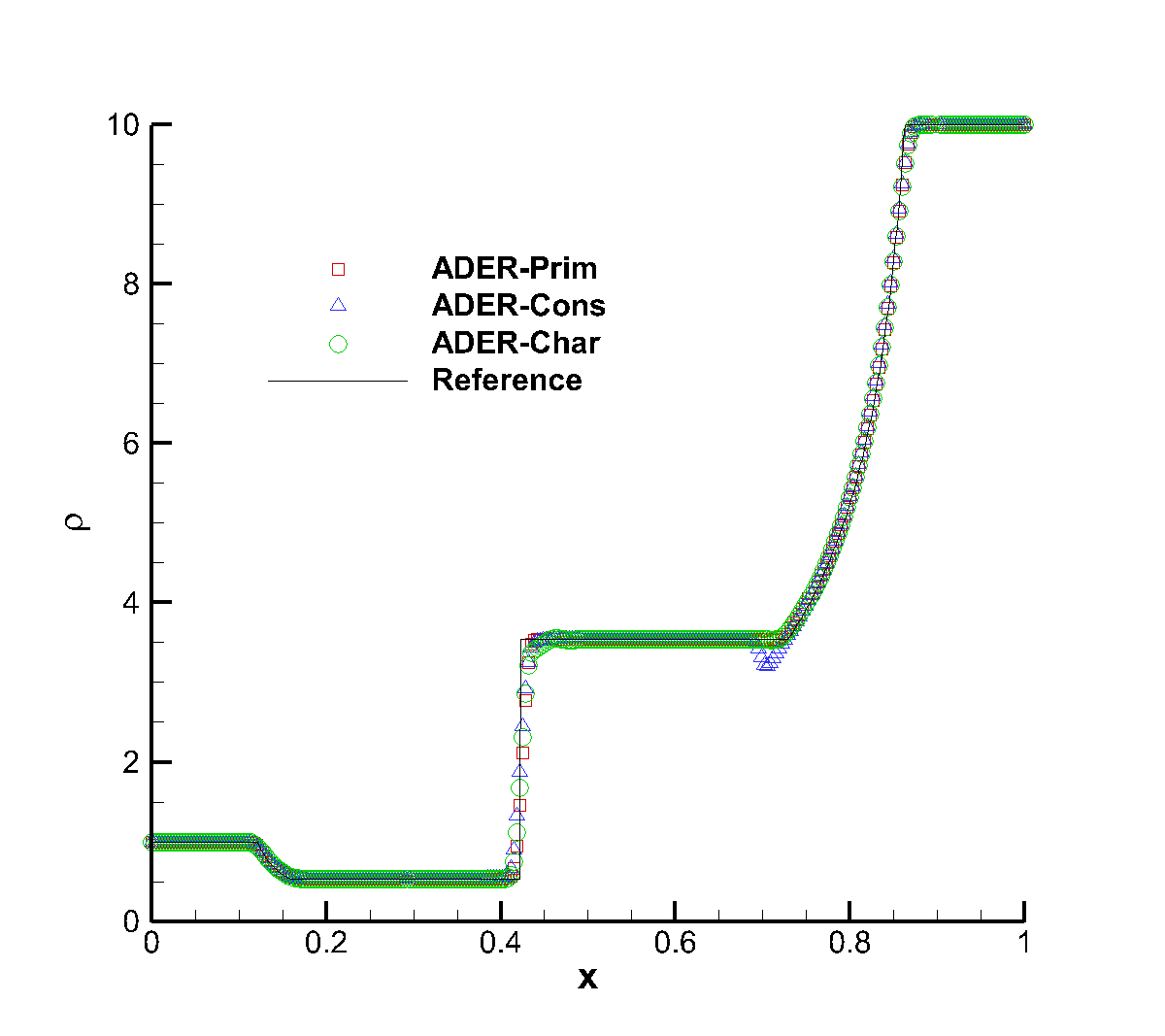}} & 
{\includegraphics[angle=0,width=7.3cm,height=7.3cm]{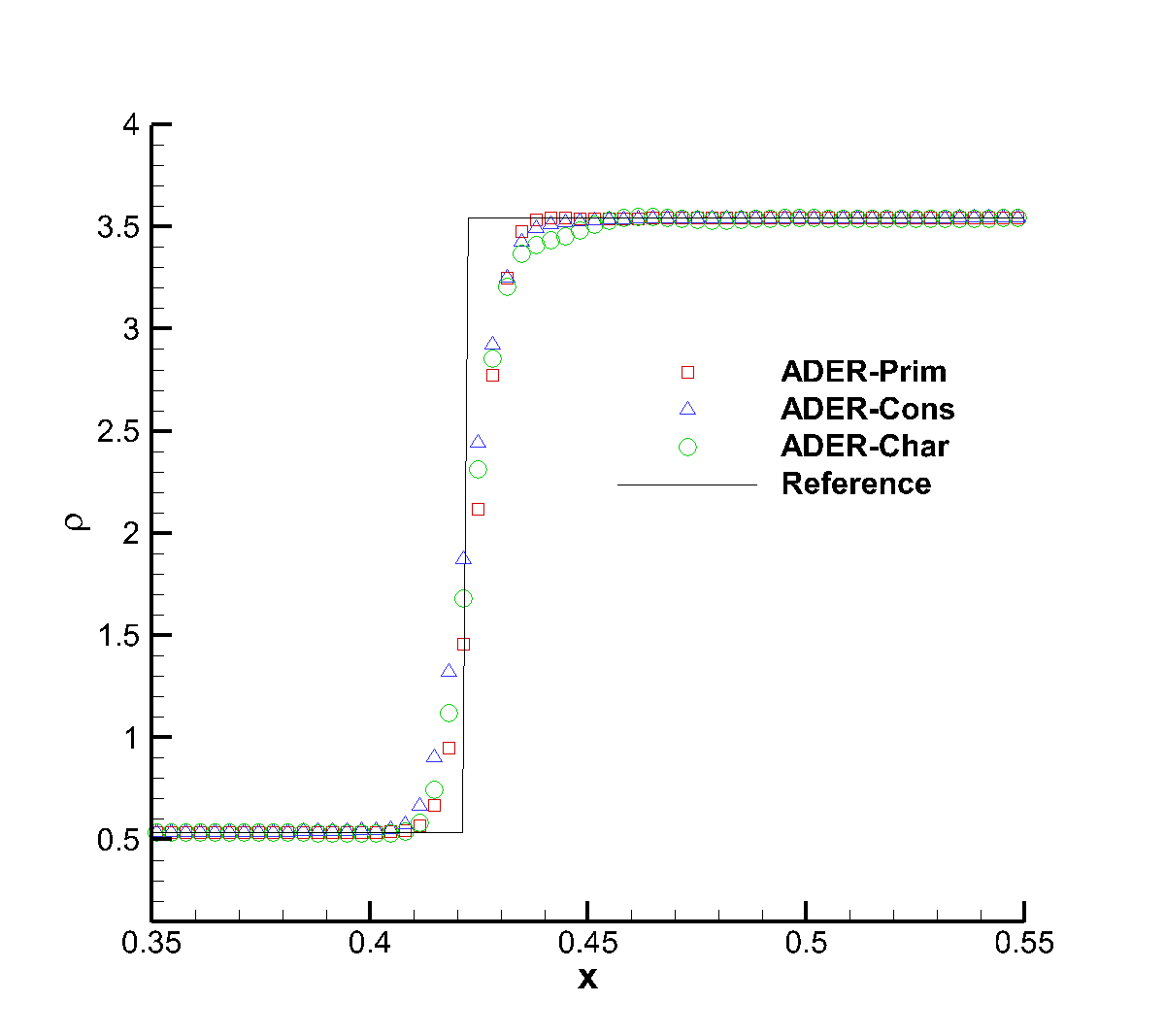}}
\end{tabular} 
\caption{Solution of RHD-RP1 (see Table~\ref{tab.RP.ic}) with the
  fourth order ADER-WENO scheme at time $t=0.4$. The bottom right panel shows a magnification around the contact discontinuity.}
\label{fig:shock-tube-2R}
\end{center}
\end{figure}

\begin{figure}
\begin{center}
\begin{tabular}{cc} 
{\includegraphics[angle=0,width=7.3cm,height=7.3cm]{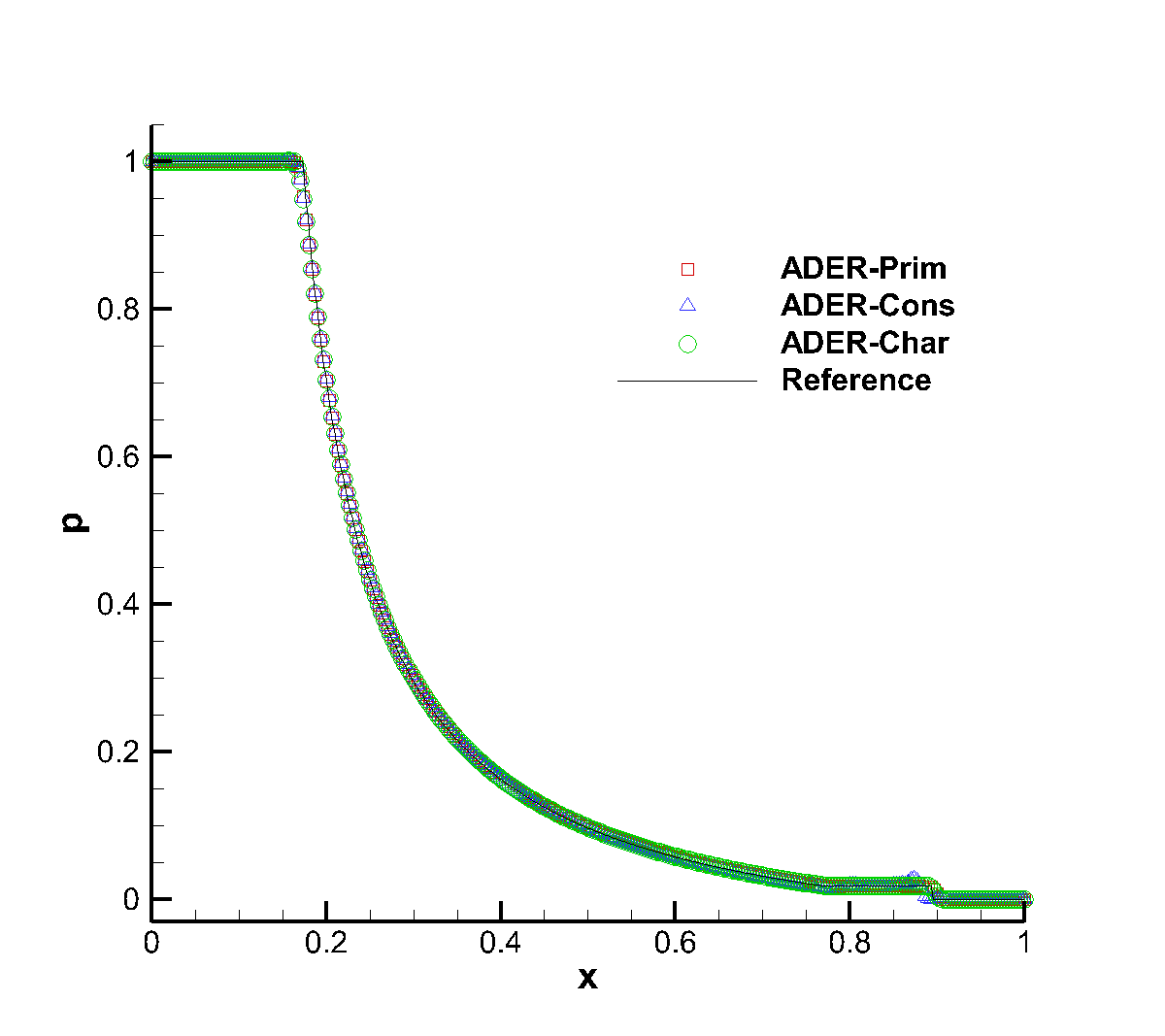}} & 
{\includegraphics[angle=0,width=7.3cm,height=7.3cm]{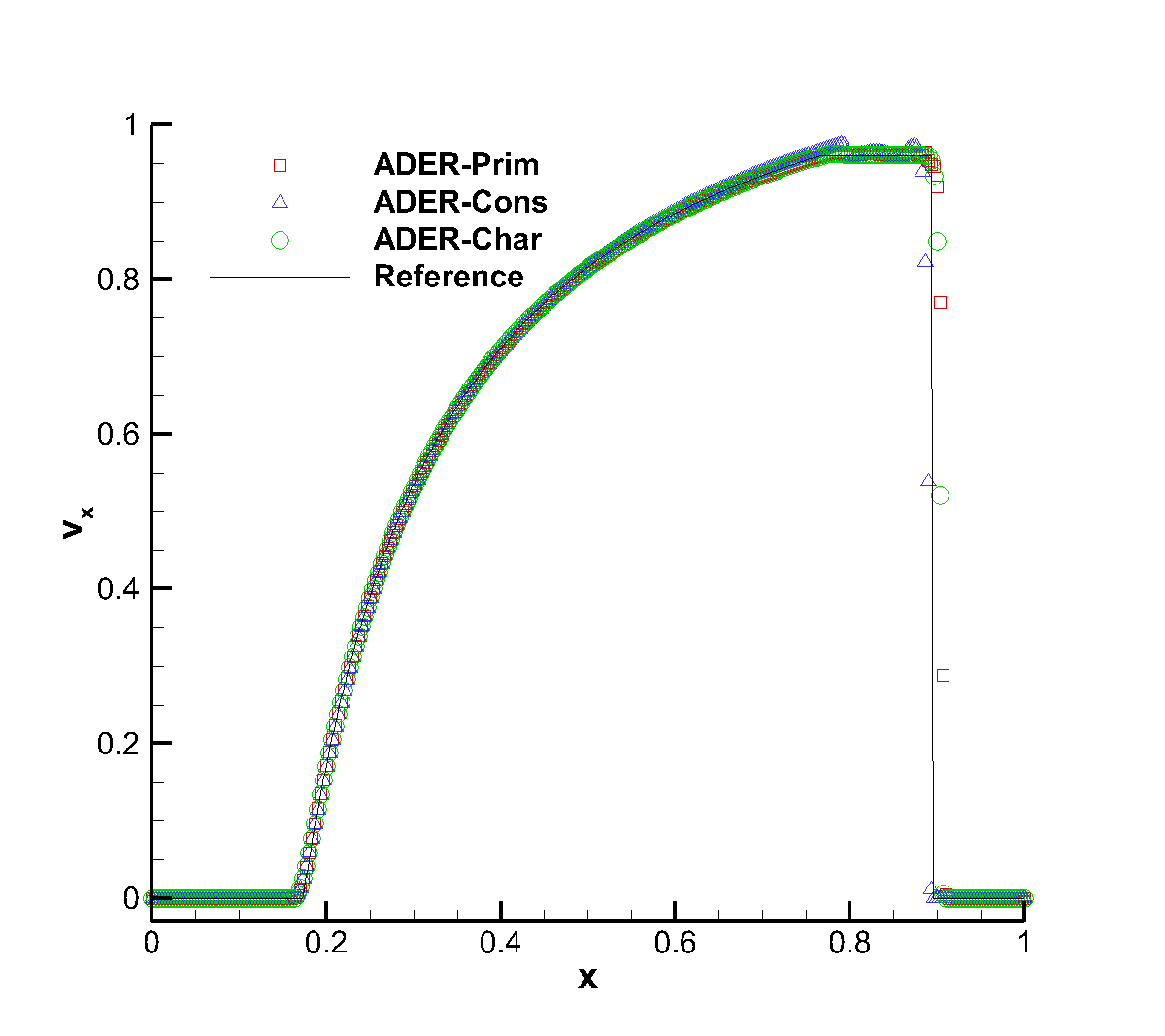}} \\ 
{\includegraphics[angle=0,width=7.3cm,height=7.3cm]{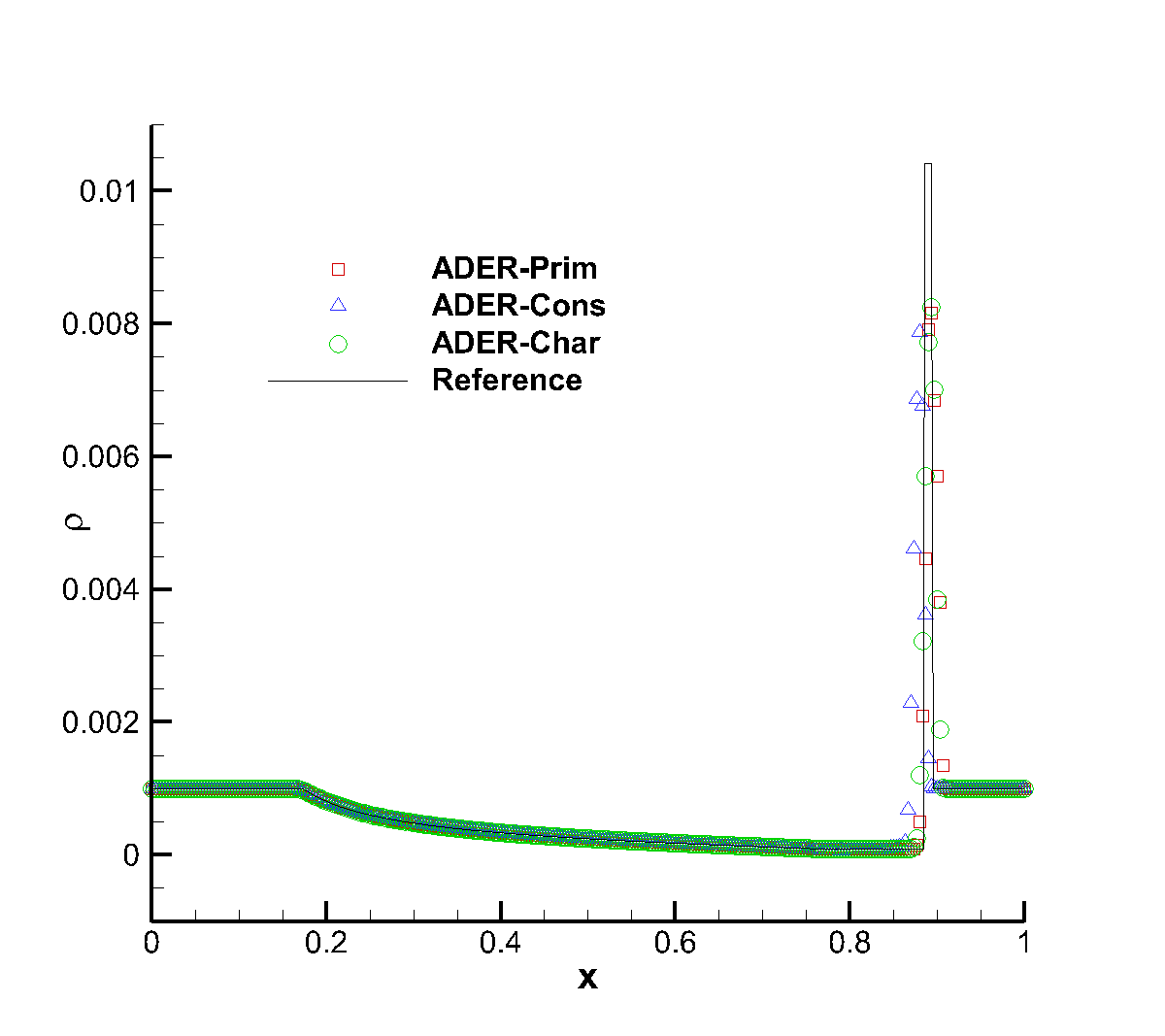}} & 
{\includegraphics[angle=0,width=7.3cm,height=7.3cm]{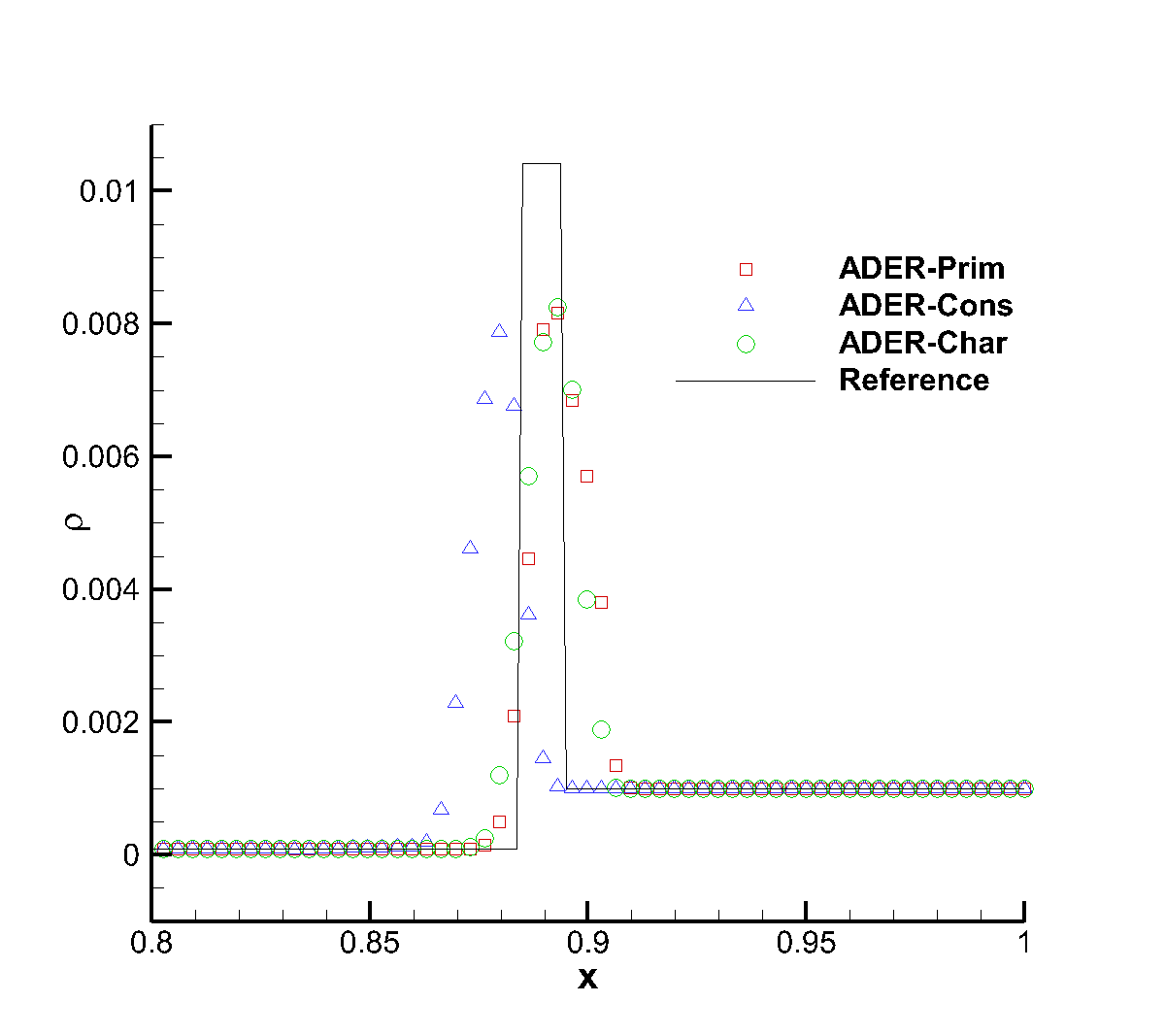}}
\end{tabular} 
\caption{Solution of RHD-RP2 (see Table~\ref{tab.RP.ic}) with the
  third order ADER-WENO scheme at time $t=0.4$. The bottom right panel shows a magnification around the right propagating shock.}
\label{fig:shock-tube-RS}
\end{center}
\end{figure}

%
\begin{table} 
\begin{center} 
\begin{tabular}{c|c||c|ccc|c} 
\hline
\hline
Problem  &   & $\gamma$   & $\rho$ &$v_x$ & $p$ &  $t_{\text{f}}$  \\
\hline
\multirow{2}{*}{\rotatebox{0}{\textbf{RHD-RP1}} } 
&$x > 0$     &\multirow{2}{*}{$\left.5\middle/ 3\right.$} & 1     &  -0.6    & 10   & \multirow{2}{*}{0.4}\\ 
&$x \leq 0$  &                                            & 10    &   0.5    & 20   &    \\ 
\hline
\multirow{2}{*}{\rotatebox{0}{\textbf{RHD-RP2}}}
&$x > 0$     &\multirow{2}{*}{ $\left.5\middle/ 3\right.$} & $10^{-3}$  &  0.0  &  1         & \multirow{2}{*}{0.4}\\ 
&$x \leq 0$  &                                             & $10^{-3}$  &  0.0  &  $10^{-5}$ &  \\ 
\hline
\end{tabular} 
\caption{
\label{tab.RP.ic}
Left and right states of the one--dimensional RHD Riemann problems.} 
\end{center}
\end{table} 
Table~\ref{tab.RP.ic} reports the initial conditions of the two one-dimensional Riemann problems that we have considered, 
and whose wave-patterns at the final time $t_f=0.4$ are shown in Fig.~\ref{fig:shock-tube-2R} and Fig.~\ref{fig:shock-tube-RS}, respectively.
In order to appreciate the differences among the available ADER implementations, we have again solved each problem 
with the three alternative schemes: ADER-Prim, ADER-Cons and ADER-Char.  The reference solution, 
computed as in \cite{Rezzolla01}, is shown too. 
%
%
%

In the first Riemann problem, which was also analyzed by \cite{Mignone2005}, 
two rarefaction waves are produced, separated by a contact discontinuity.
It has been solved through a fourth order $\mathbb{P}_0\mathbb{P}_3$  scheme, using the Rusanov Riemann solver over a uniform grid with $300$ cells. 
As it is clear from Fig.~\ref{fig:shock-tube-2R}, the ADER-Prim scheme performs significantly better than the ADER-Cons. In particular, the overshoot and undershoot 
at the tail of the right rarefaction is absent. In general, the results obtained with ADER-Prim
are essentially equivalent to those of ADER-Char, namely when the reconstruction in characteristic variables is adopted. This is manifest after looking at the
bottom right panel of Fig.~\ref{fig:shock-tube-2R}, where a magnification of the rest mass density at the contact discontinuity is shown.
Additional interesting comparisons can be made about the second Riemann problem, which can be found in \cite{Radice2012a}, and which is displayed in 
Fig.~\ref{fig:shock-tube-RS}. In this case a third order $\mathbb{P}_0\mathbb{P}_2$ scheme has been used, again with 
the Rusanov Riemann solver over a uniform grid with $500$ cells. The right propagating shock has a strong jump in the rest mass density, as it is visible from the 
bottom right panel of the figure, and the position of the shock front is better captured by the two schemes ADER-Prim and ADER-Char.

\begin{table}[!t]
\vspace{0.5cm}
\renewcommand{\arraystretch}{1.0}
\begin{center}
\begin{tabular}{cccc}
\hline
\hline
             & \footnotesize{ADER-Prim}   & \footnotesize{ADER-Cons}   &  \footnotesize{ADER-Char} \\
\hline
\hline
$\mathbb{P}_0\mathbb{P}_2$         & 1.0                                  & 1.26                               & 1.40   \\
\hline
$\mathbb{P}_0\mathbb{P}_3$         & 1.0                                  & 1.13                               & 1.24 \\
\hline
$\mathbb{P}_0\mathbb{P}_4$         & 1.0                                  & 1.04                               & 1.06 \\
\hline
\hline
\end{tabular}
\end{center}
\caption{CPU time comparison among different ADER implementations for the RHD-RP1 problem. 
The numbers have been normalized to the value obtained with ADER-Prim.} 
\label{tab.CPU.RHD}
\end{table}
It is particularly interesting to address the issue of CPU time comparison among different implementations of ADER, as already done for the Euler equations.
The result of such a comparison, performed for the RHD-RP1 problem, are reported in Tab.~\ref{tab.CPU.RHD}, which should be read in synopsis with Tab.~\ref{tab.CPU.Sod}.
Clearly, ADER-Prim is not only more accurate than ADER-Cons, but it is also more efficient. As anticipated, this is in agreement with our expectations, since
in the ADER-Prim implementation a single {\em cons-to-prim} operation is needed within the cell, rather than at each Gaussian quadrature point and at each space-time degree
of freedom.  
For other tests, see for instance Sect.~\ref{sec:RHD-KH}, the CPU time reduction implied by ADER-Prim is even more evident, but the numbers shown in Tab.~\ref{tab.CPU.RHD} 
describe with good fidelity the relative performances of the different ADER in a large number of relativistic tests.

\subsubsection{RHD Kelvin--Helmholtz instability}
\label{sec:RHD-KH}
%
\begin{figure}
\begin{center}
\begin{tabular}{ccc} 
{\includegraphics[angle=0,width=4.0cm,height=8.0cm]{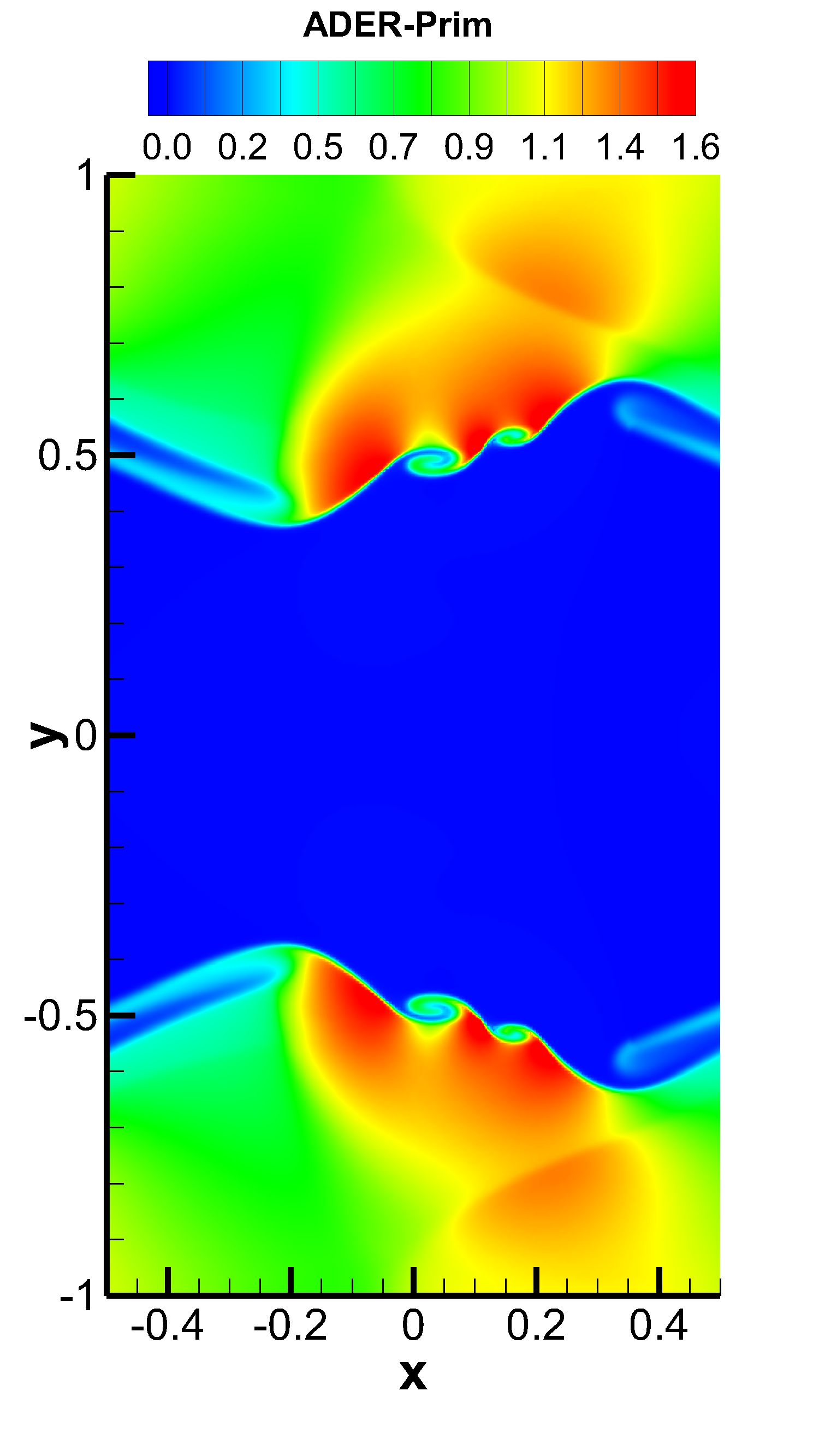}} & 
{\includegraphics[angle=0,width=4.0cm,height=8.0cm]{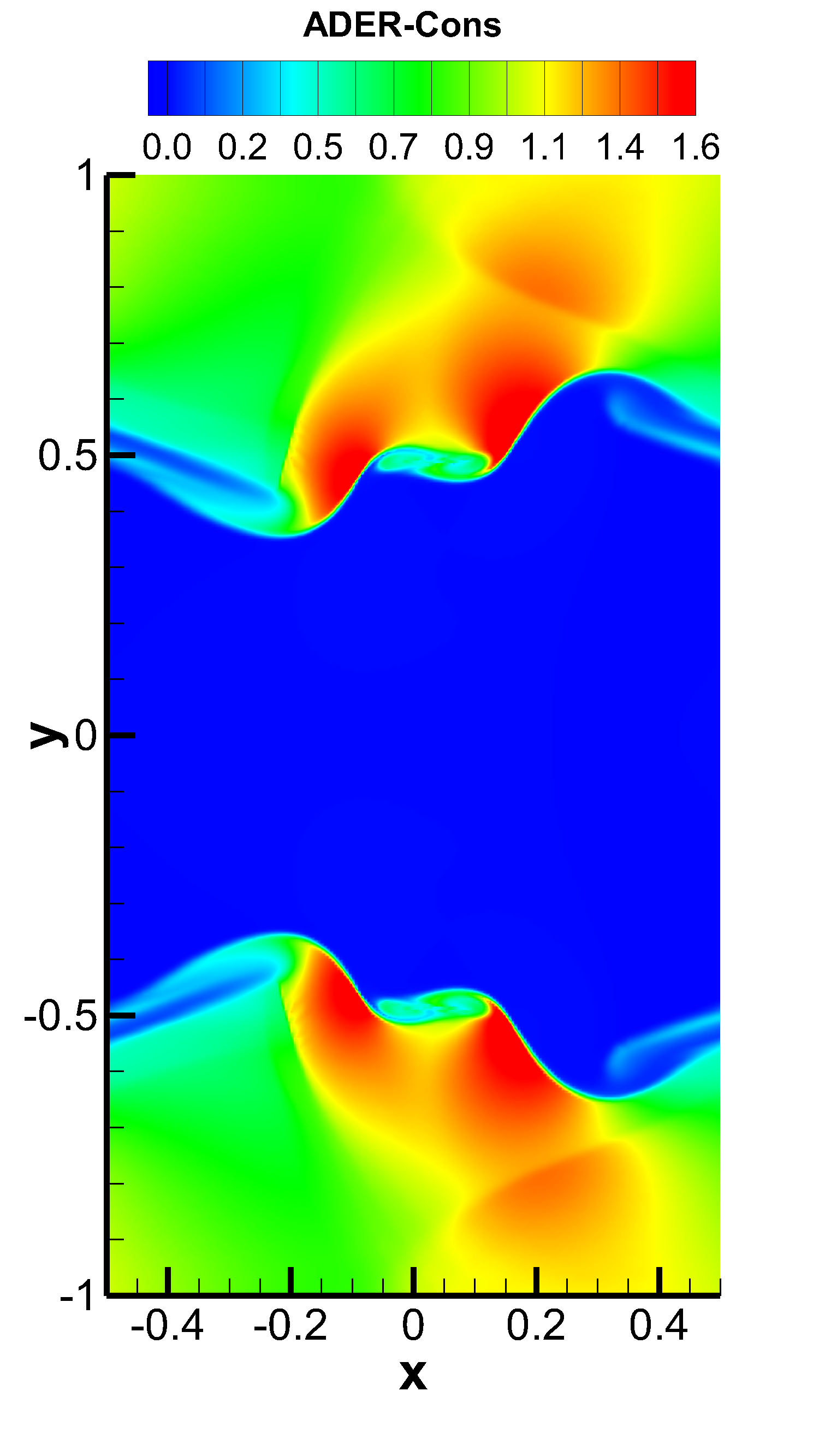}} &
{\includegraphics[angle=0,width=4.0cm,height=8.0cm]{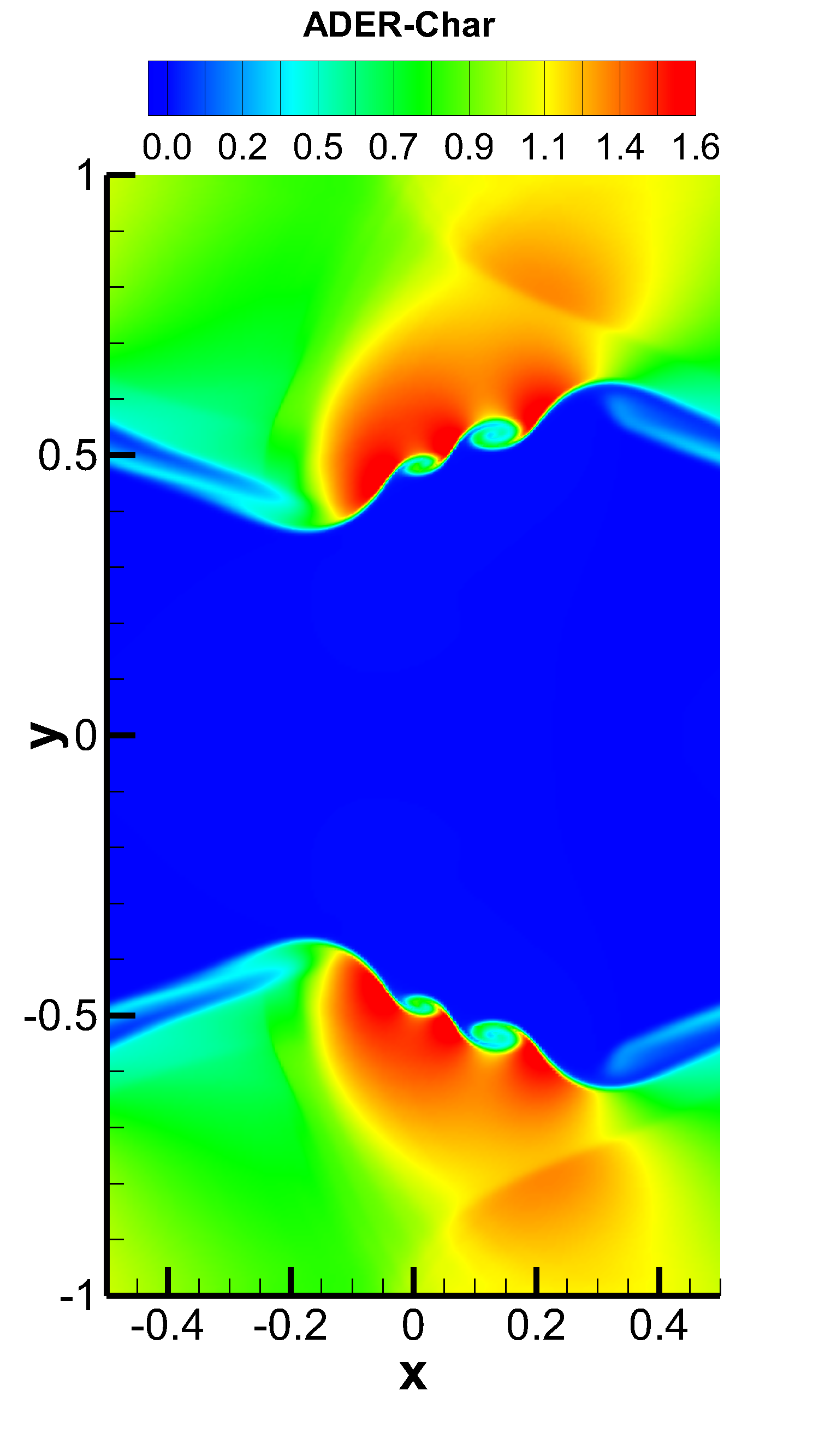}} \\ 
{\includegraphics[angle=0,width=4.0cm,height=8.0cm]{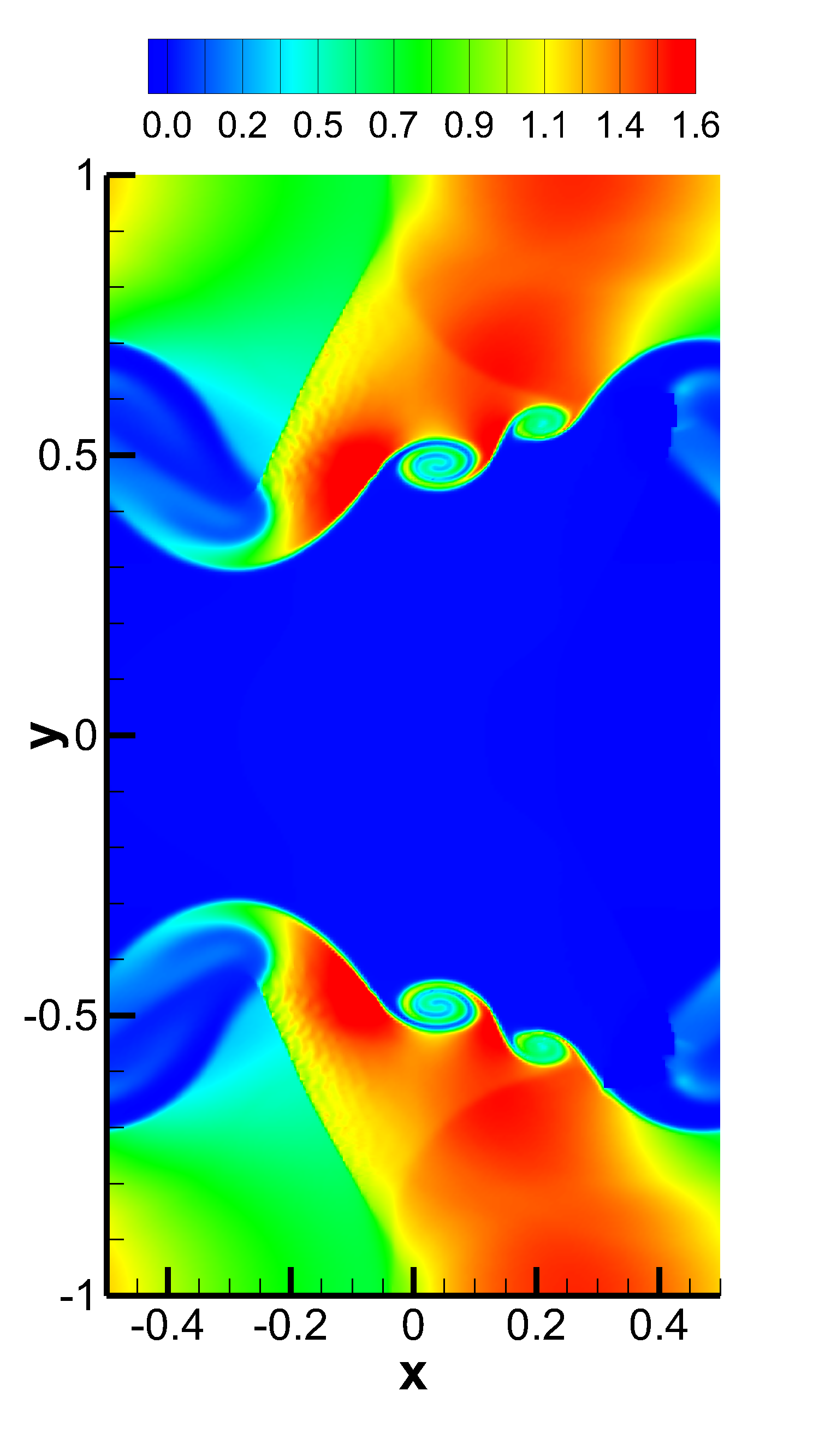}} & 
{\includegraphics[angle=0,width=4.0cm,height=8.0cm]{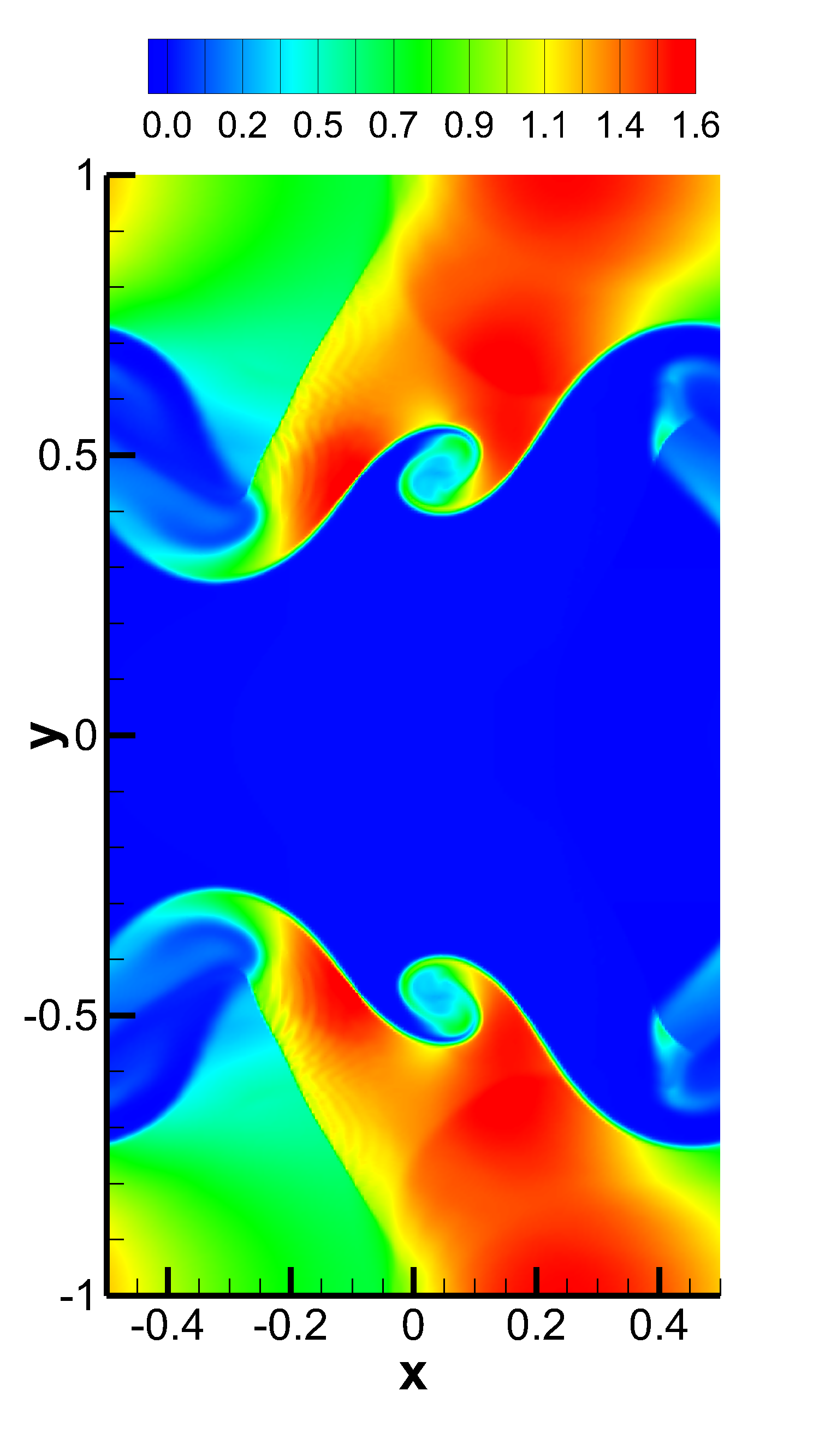}} &
{\includegraphics[angle=0,width=4.0cm,height=8.0cm]{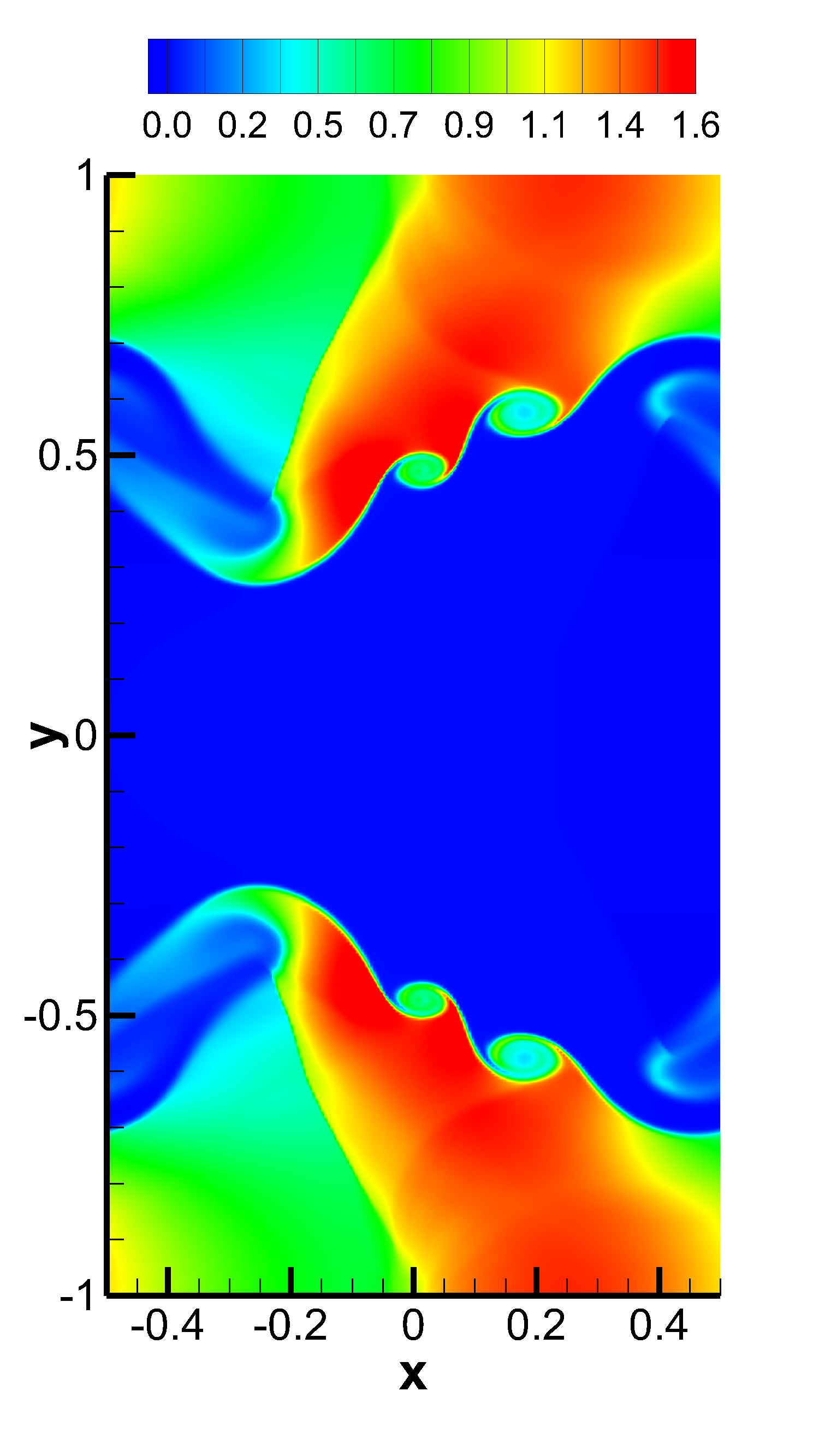}}
\end{tabular} 
\end{center}
\caption{
Two-dimensional Kelvin-Helmholtz instability obtained with the $\mathbb{P}_0\mathbb{P}_3$ scheme and with the Osher flux.
Left panels: solution with ADER-Prim. Central panels: solution  with ADER-Cons.
Right panels: solution  with ADER-Char.
Top panels: solution at $t=2.0$.
Bottom panels: solution at $t=2.5$.
}
\label{fig:KH-RHD}
\end{figure}
In the relativistic regime, the Kelvin--Helmholtz (KH) instability is likely to be responsible for a variety of physical effects, which are encountered in the
dynamics of extragalactic relativistic jets \cite{Bodo2004,Perucho2006,Perucho2007}. As an academic test, we simulate  
the linear growth phase of the KH instability in two spatial dimensions, taking the initial conditions from 
\cite{Mignone2009} (see also \cite{Beckwith2011} and \cite{Radice2012a}). In particular, 
the rest-mass density is chosen as
\begin{equation}\label{KHI-rho}
  \rho = \left\{\begin{array}{ll}
  \rho_0 + \rho_1 \tanh{[(y-0.5)/a]} & \quad y > 0\,, \\
 \noalign{\medskip}
 \rho_0 - \rho_1 \tanh{[(y+0.5)/a]}  & \quad y \leq 0    \,, 
 \end{array}\right.
\end{equation}
with $\rho_0=0.505$ and $\rho_1=0.495$. 
Assuming that the shear layer has a velocity $v_s=0.5$ and a characteristic size $a=0.01$, the velocity along the x-direction is modulated as
\begin{equation}\label{KHI-vx}
  v_x = \left\{\begin{array}{ll}
  v_s \tanh{[(y-0.5)/a]} & \quad y > 0\,, \\
 \noalign{\medskip}
 -v_s \tanh{[(y+0.5)/a]}  & \quad y \leq 0    \,. 
 \end{array}\right.
\end{equation}
It is convenient to  add  a perturbation in the transverse velocity, i.e.
\begin{equation}\label{KHI-vy}
  v_y = \left\{\begin{array}{ll}
  \eta_0 v_s \sin{(2\pi x)} \exp{[-(y-0.5)^2/\sigma]} & \quad y > 0\,, \\
 \noalign{\medskip}
 -\eta_0 v_s \sin{(2\pi x)} \exp{[-(y+0.5)^2/\sigma]}  & \quad y \leq 0    \,, 
 \end{array}\right.
\end{equation}
where $\eta_0=0.1$ is the amplitude of the perturbation, while $\sigma=0.1$ is its length scale. 
The adiabatic index is $\gamma=4/3$ and the pressure is uniform, $p=1$.
The problem has been solved over the computational domain
$[-0.5,0.5]\times[-1,1]$, covered by a uniform mesh with $200\times400$ cells, using the $\mathbb{P}_0\mathbb{P}_3$ scheme and the Osher-type numerical flux.
Periodic boundary conditions are fixed both in $x$ and in $y$ directions. 
Fig.~\ref{fig:KH-RHD} shows the results of the calculations: 
in the left, in the central and in the right panels we have reported the solution obtained with the ADER-Prim, with the ADER-Cons and with the ADER-Char scheme, respectively, 
while the top and the bottom panels correspond to two different times during the evolution, namely $t=2.0$ and $t=2.5$. 
Interestingly, 
two secondary vortices are visible when the reconstruction is performed in primitive and characteristic variables (see left the right panels),
but only one is present in the simulation using the reconstruction in conserved variables.
In \cite{Zanotti2015} we have already commented about the elusive character of these
details in the solution, which depend both on the resolution and on the Riemann solver adopted. Based on our results, we infer that the 
ADER-Cons scheme is the most diffusive, while ADER-Prim and ADER-Char seem to produce the same level of accuracy in the solution.
However, if we look at the CPU times in the two cases, 
we find that ADER-Prim is a factor 2.5 faster than ADER-Cons and a factor 3 faster than ADER-Char, and therefore should be preferred in all relevant applications of RHD.
 \begin{table}[t] 
 \centering
 \begin{tabular}{|c|c||cc|cc|c|}
   \hline
   \multicolumn{7}{|c|}{\textbf{2D circularly polarized Alfv\'en wave }} \\
   \hline
	\hline
      \hline
    & $N_x$ &  $L_1$ error &   $L_1$ order &  $L_2$ error &   $L_2$ order  &    Theor. \\
   \hline
   \hline
   \multirow{5}{*}{\rotatebox{0}{{$\mathbb{P}_0\mathbb{P}_2$}}}
& 50	&  5.387E-02	& ---	   & 9.527E-03	& ---	   & \multirow{5}{*}{3}\\
& 60	&  3.123E-02	&  2.99	 & 5.523E-03	&  2.99	 & \\
& 70	&  1.969E-02	&  2.99	 & 3.481E-03	&  2.99	 & \\
& 80	&  1.320E-02	&  2.99  & 2.334E-03	&  2.99	 & \\
& 100	&  6.764E-03	&  3.00  & 1.196E-03	&  3.00	 & \\		
   \hline
   \multirow{5}{*}{\rotatebox{0}{{$\mathbb{P}_0\mathbb{P}_3$}}}
& 50	&  2.734E-04	&  ---	 & 4.888E-05	& ---	   & \multirow{5}{*}{4}\\
& 60	&  1.153E-04	&  4.73	 & 2.061E-05	&  4.74	 & \\
& 70	&  5.622E-05	&  4.66  & 1.004E-05	&  4.66	 & \\
& 80	&  3.043E-05	&  4.60  & 5.422E-06	&  4.61	 & \\
& 100	&  1.108E-05	&  4.53  & 1.968E-06	&  4.54	 & \\
   \hline
	  \multirow{5}{*}{\rotatebox{0}{{$\mathbb{P}_0\mathbb{P}_4$}}}
& 30	&  2.043E-03	& ---	   & 3.611E-04	& ---	   & \multirow{5}{*}{5}\\
& 40	&  4.873E-04	&  4.98	 & 8.615E-05	&  4.98	 & \\
& 50	&  1.603E-04	&  4.98  & 2.846E-05	&  4.96	 & \\
& 60	&  6.491E-05	&  4.96  & 1.168E-05	&  4.88	 & \\
& 70	&  3.173E-05	&  4.64  & 6.147E-06	&  4.16	 & \\
   \hline
  \end{tabular}
 \caption{   $L_1$ and $L_2$  errors analysis for the 
   2D Alfv\'en wave problem. The errors have been computed with respect to the magnetic field $B^y$.}
	\label{tab:Alfven_Error}
 \end{table}
%

\subsubsection{RMHD Alfv\'en Wave}
\label{sec:RMHD_Alfven_Wave}

In Tab.~\ref{Table:convergence} of Sect.~\ref{sec:isentropic} we have reported the comparison of the convergence rates among three different implementations of ADER for the Euler equations. We believe it is important to verify the
convergence of the new ADER-Prim scheme also for the RMHD equations, which indeed admits an exact, smooth unsteady solution, namely  
the propagation of a circularly polarized Alfv\'en wave (see \cite{Komissarov1997,DelZanna2007} for a full account). The wave is assumed to propagate along 
the $x$ direction in a constant density and constant pressure background, say $\rho=p=1$. The magnetic field, on the other hand, 
is given by
\begin{eqnarray}
B_x&=&B_0 \\
B_y&=&\eta B_0\cos[k(x-v_A t)]\\
B_z&=&\eta B_0\sin[k(x-v_A t)]\,,
\end{eqnarray}
where $\eta=1$ is the amplitude of the wave,
$B_0=1$ is the uniform magnetic field, $k$ is the wave number, while $v_A$ is speed of propagation of the wave.
We have solved this problem over the computational domain $\Omega=[0; 2\pi]\times[0; 2\pi]$, 
using periodic boundary conditions,
the Rusanov Riemann solver and the Adams--Bashforth version for the initial guess of the LSDG predictor.
We have compared the numerical solution with the analytic one
after one period $T=L/v_A=2\pi/v_A$.
Tab.~\ref{tab:Alfven_Error} contains the results of our analysis, showing 
the $L_1$ and the $L_2$ norms of the error of $B^y$. As apparent from the table, the nominal order of convergence of the new ADER-Prim scheme is recovered with very good
accuracy.
%

\subsubsection{RMHD Riemann Problems}

\begin{figure}
\begin{center}
\begin{tabular}{cc} 
{\includegraphics[angle=0,width=7.3cm,height=7.3cm]{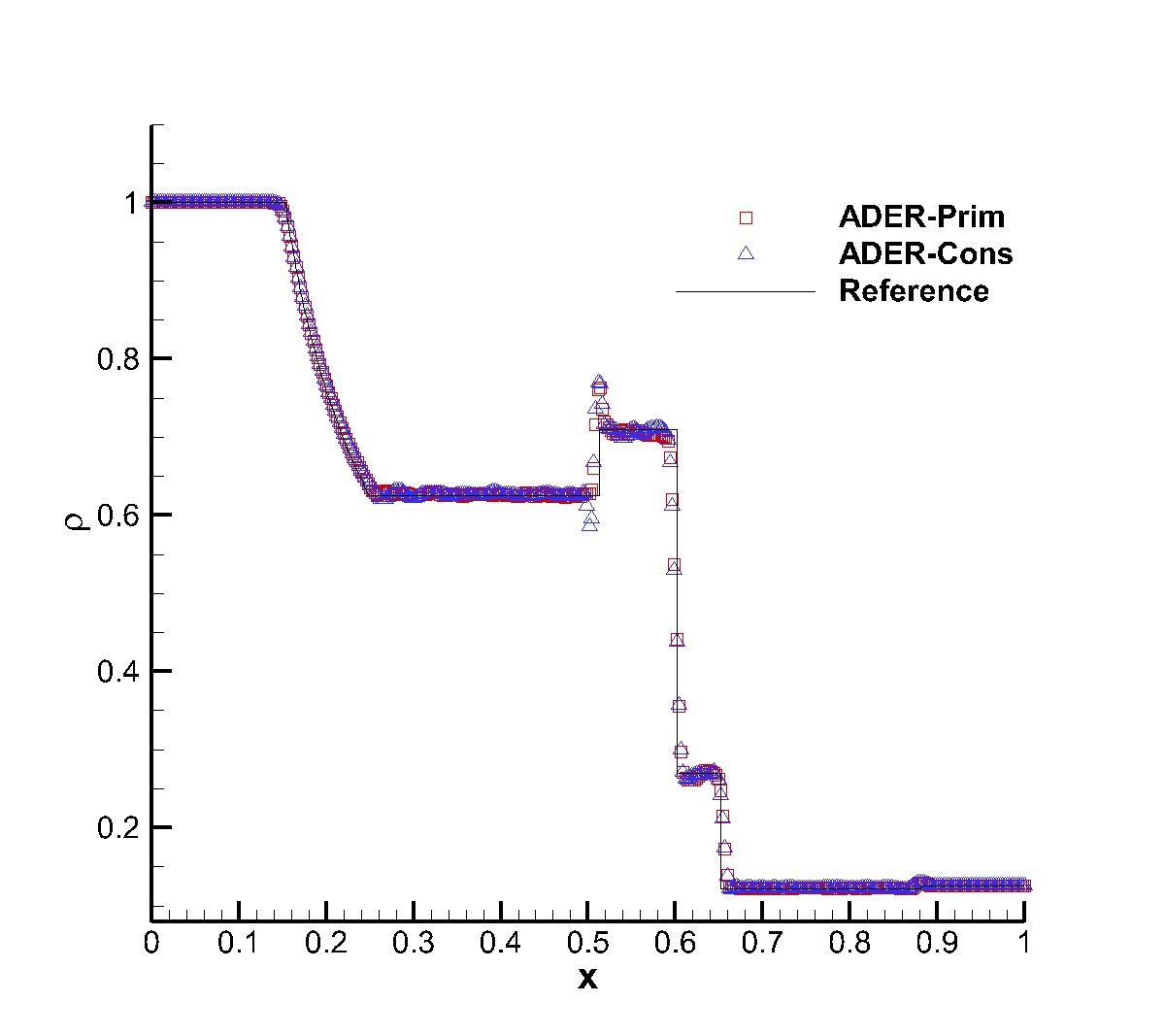}} & 
{\includegraphics[angle=0,width=7.3cm,height=7.3cm]{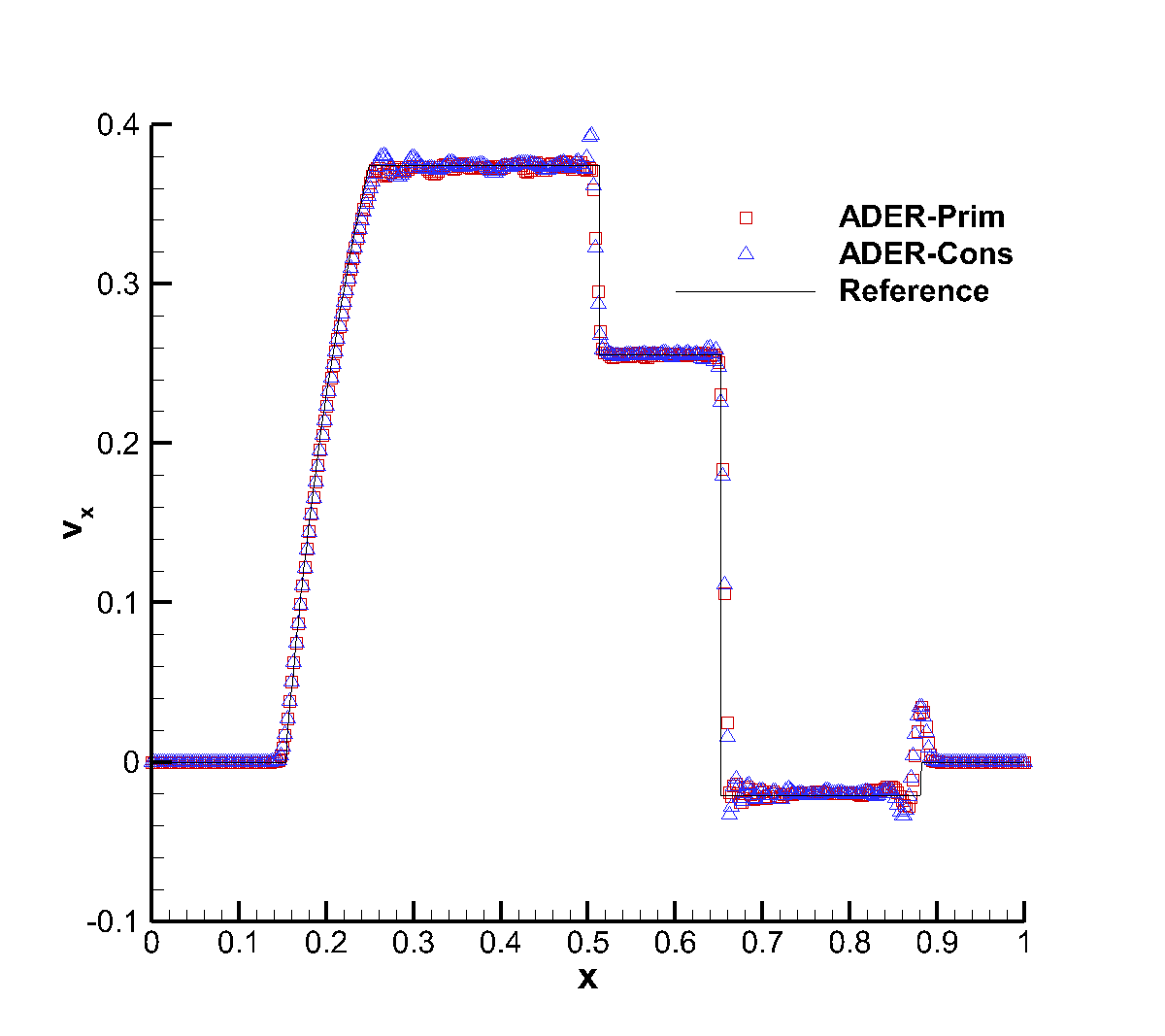}}  \\
{\includegraphics[angle=0,width=7.3cm,height=7.3cm]{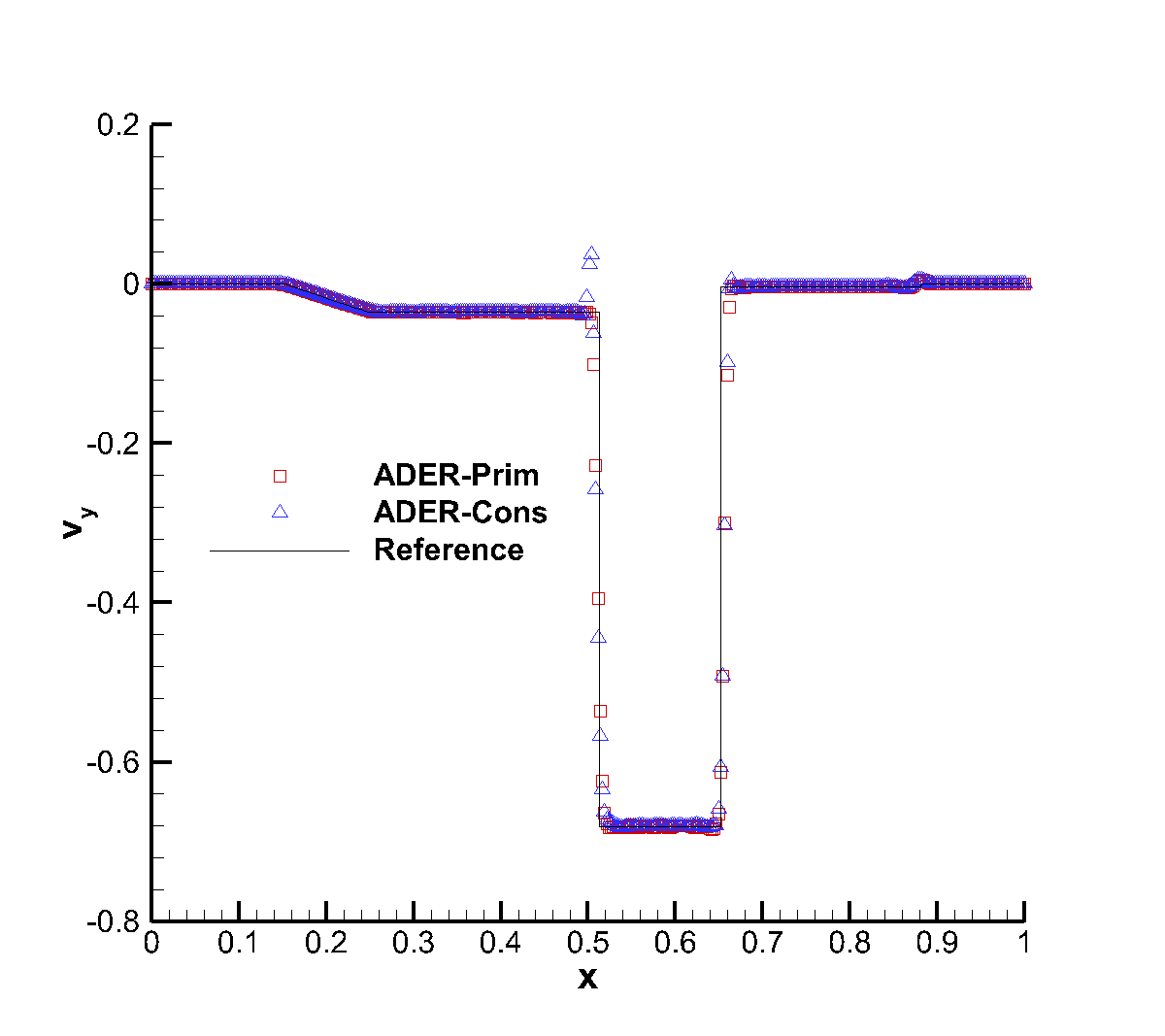}}  & 
{\includegraphics[angle=0,width=7.3cm,height=7.3cm]{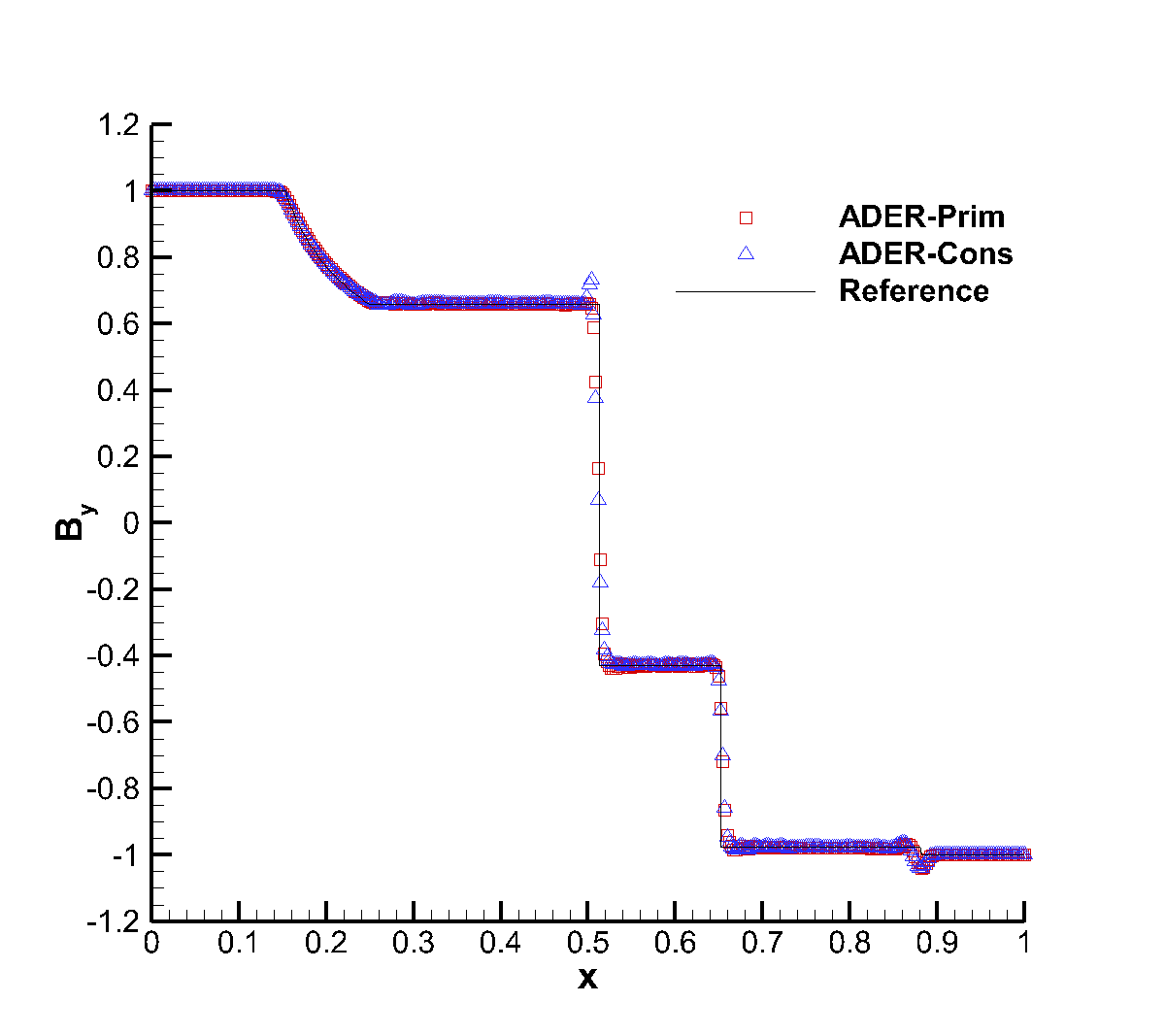}}
\end{tabular} 
\caption{Solution of RMHD-RP1 (see Tab.~\ref{tab:RMHD-RP}) with the
  fourth order ADER-WENO scheme at time $t=0.4$. The Rusanov Riemann solver has been used over a $400$ cells uniform grid.}
\label{fig:Balsara1}
\end{center}
\end{figure}

\begin{figure}
\begin{center}
\begin{tabular}{cc} 
{\includegraphics[angle=0,width=7.3cm,height=7.3cm]{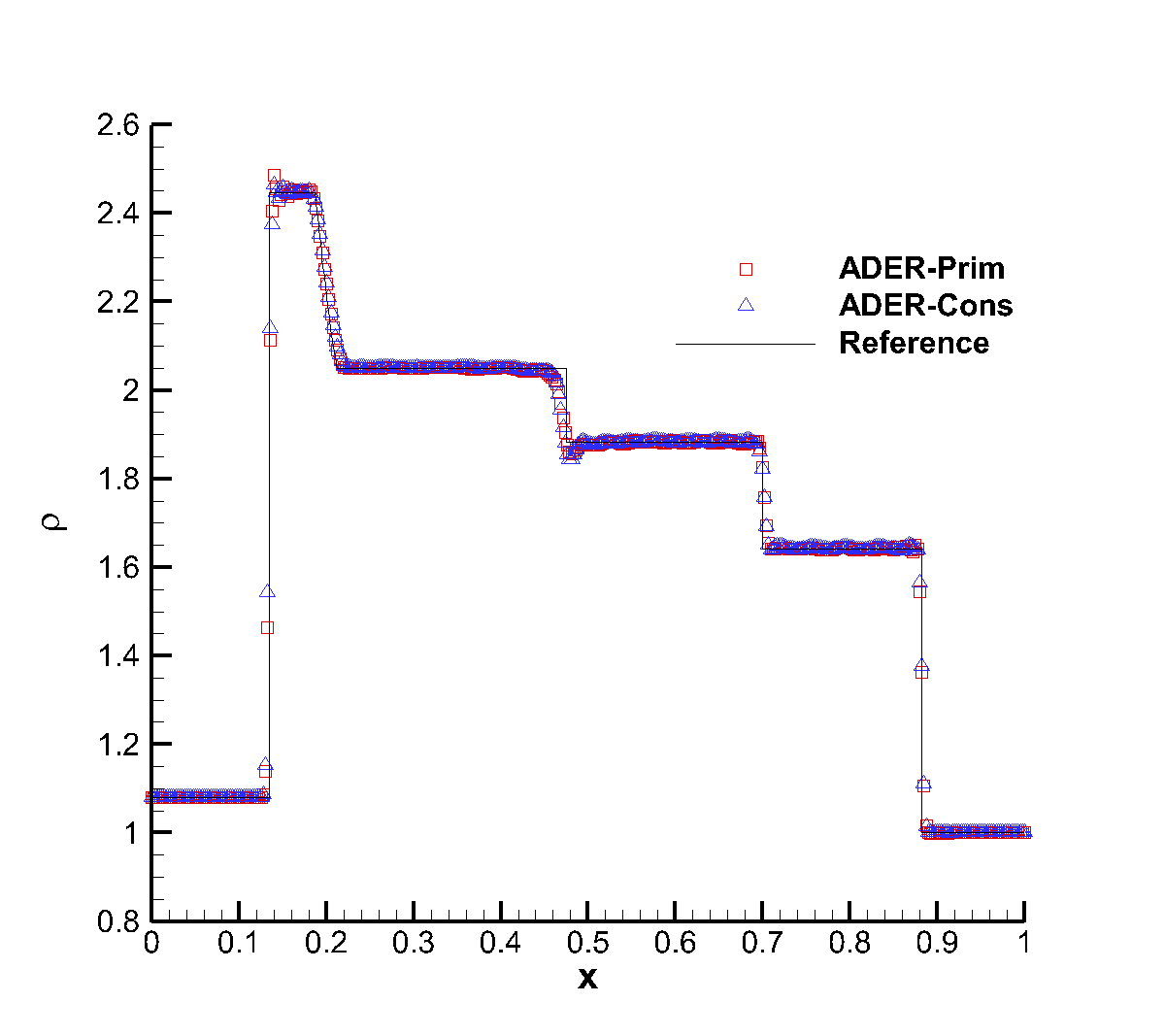}} & 
{\includegraphics[angle=0,width=7.3cm,height=7.3cm]{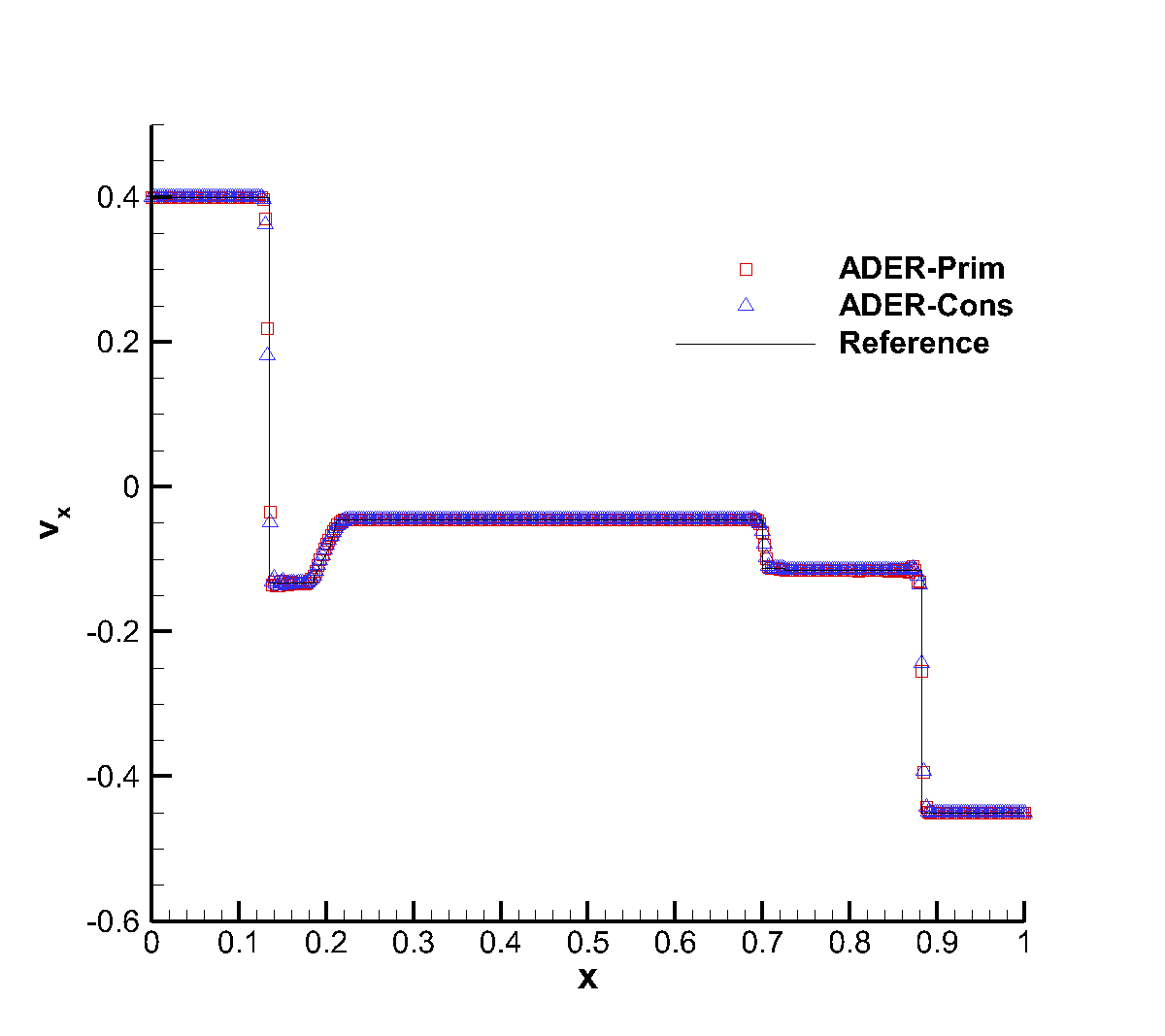}} \\
{\includegraphics[angle=0,width=7.3cm,height=7.3cm]{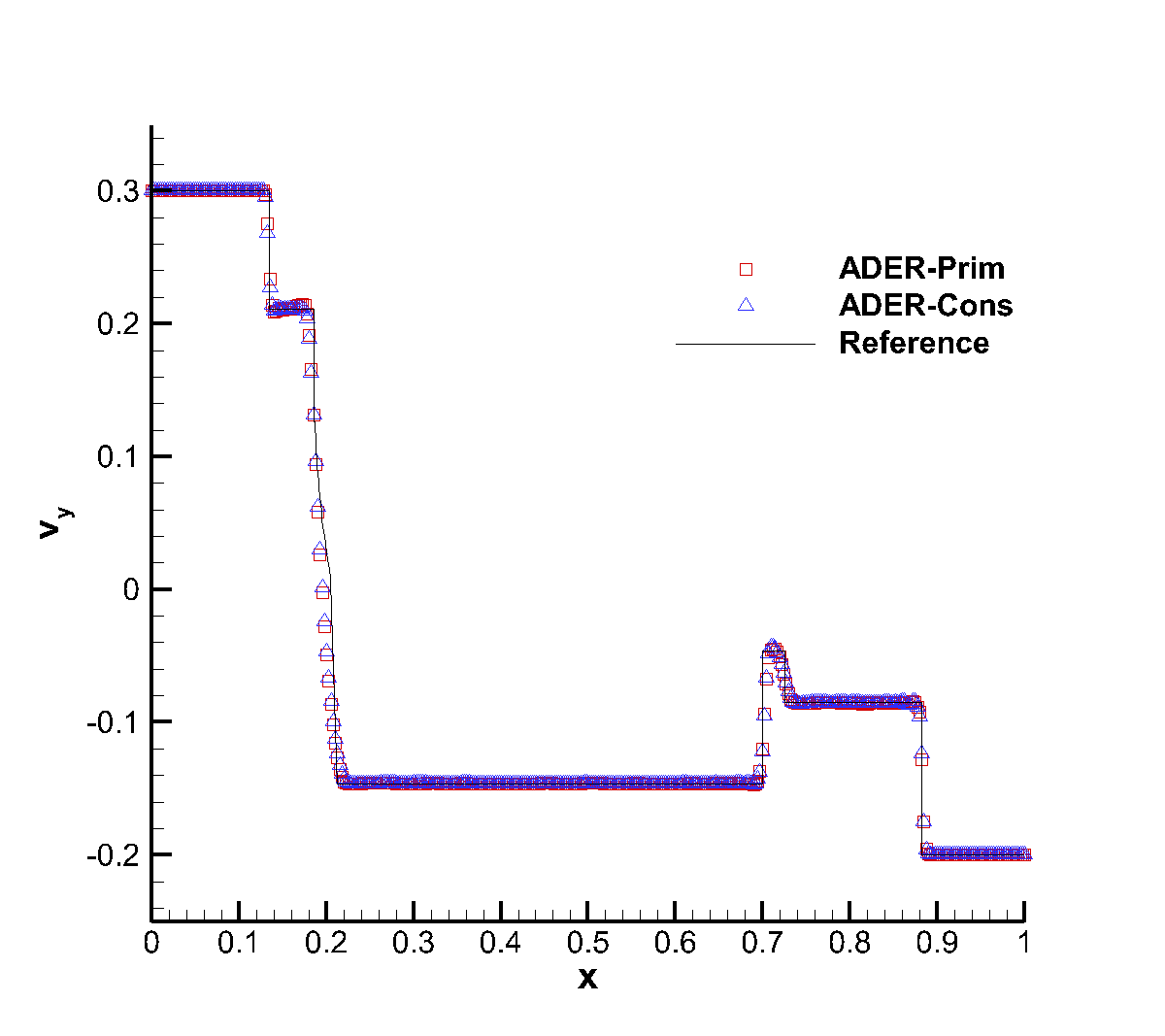}} & 
{\includegraphics[angle=0,width=7.3cm,height=7.3cm]{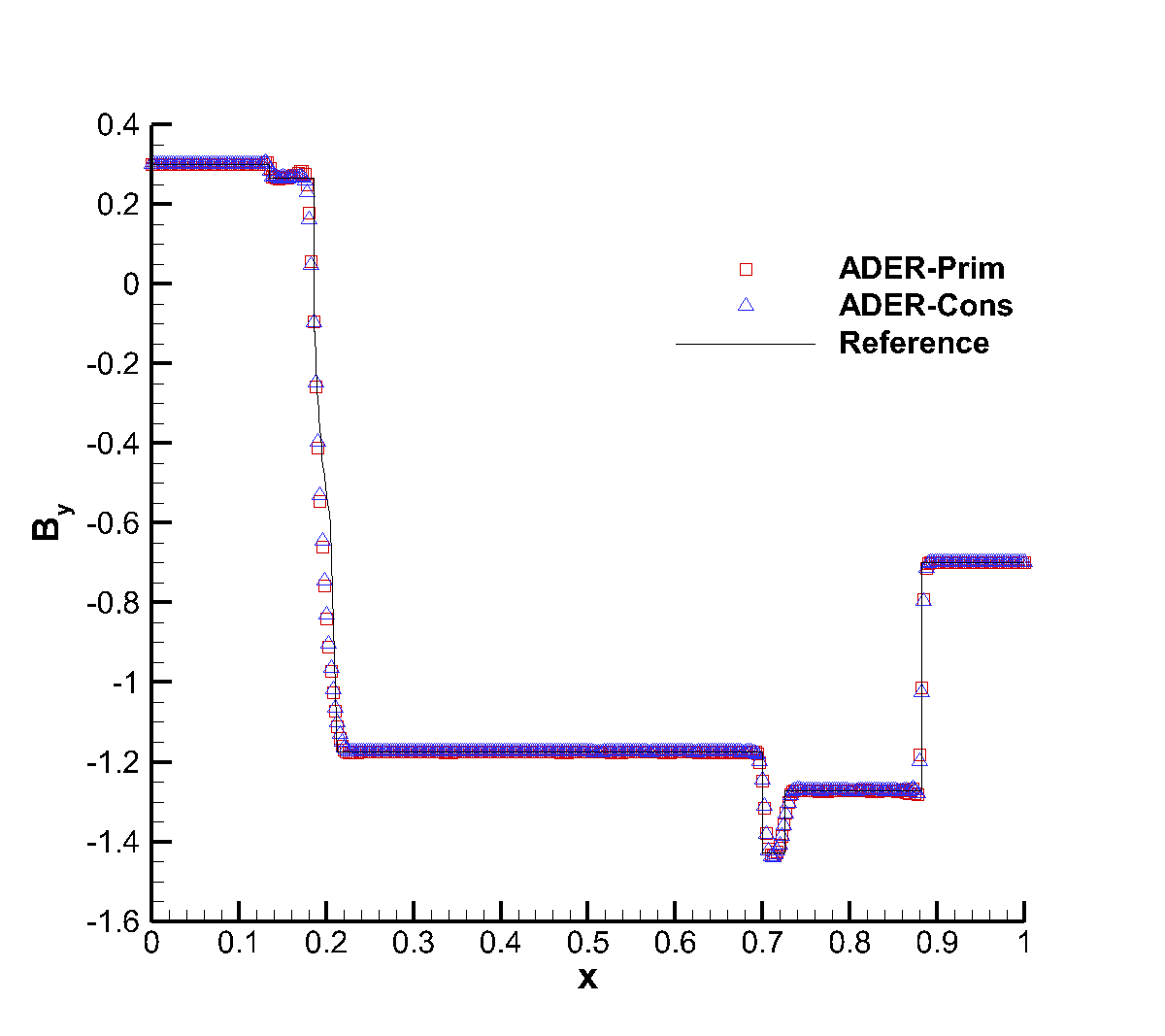}}
\end{tabular} 
\caption{Solution of RMHD-RP2 (see Tab.~\ref{tab:RMHD-RP})  with the
  fourth order ADER-WENO scheme at time $t=0.55$. The Rusanov Riemann solver has been used over a $400$ cells uniform grid.}
\label{fig:Balsara5}
\end{center}
\end{figure}

Riemann problems are very relevant also in RMHD, admitting a larger number of waves than in hydrodynamics. The exact solution was provided by 
\cite{Giacomazzo:2005jy} already ten years ago, making them very popular as a precise tool to validate numerical codes. We have selected 
Test 1 and Test 5 in Table 1 of \cite{BalsaraRMHD}, with initial left and right states that are reported in
Tab.~\ref{tab:RMHD-RP}.
%
\begin{table} 
\begin{center} 
\begin{tabular}{c|c||c|cccccccc|c} 
\hline
\hline
Problem  &   & $\gamma$   & $\rho$ &$(v_x$&$v_y$&$v_z)$ & $p$ & $(B_x$&$B_y$&$B_z)$ & $t_{\text{f}}$  \\
\hline
\multirow{2}{*}{\rotatebox{0}{\textbf{RMHD-RP1}} } 
&$x > 0$     &\multirow{2}{*}{2.0} & 0.125   &  0.0   & 0.0 &0.0 & 0.1  & 0.5&-1.0&0.0 & \multirow{2}{*}{0.4}\\ 
&$x \leq 0$  &                      & 1.0     &  0.0   & 0.0 &0.0 & 1.0  & 0.5& 1.0&0.0 &   \\ 
\hline
\multirow{2}{*}{\rotatebox{0}{\textbf{RMHD-RP2}}}
&$x > 0$     &\multirow{2}{*}{ $\left.5\middle/ 3\right.$} & 1.0   &  -0.45 & -0.2 & 0.2 & 1.0  & 2.0&-0.7&0.5 &\multirow{2}{*}{0.55}\\ 
&$x \leq 0$  &                                             & 1.08  &   0.4  & 0.3  & 0.2 & 0.95 & 2.0& 0.3&0.3 &  \\ 
\hline
\end{tabular} 
\caption{
\label{tab:RMHD-RP}
Left and right states of the one--dimensional RMHD Riemann problems.} 
\end{center}
\end{table} 
Both the tests have been solved using a fourth order ADER-WENO scheme, the Rusanov Riemann solver and over a uniform grid composed of $400$ cells.
The damping factor for the divergence-cleaning procedure is set to $\kappa=10$.
Fig.~\ref{fig:Balsara1} and Fig.~\ref{fig:Balsara5} allow to compare the exact solution 
with the results obtained through the ADER-Prim and the ADER-Cons schemes.
Especially for RMHD-RP1, the solution obtained with the traditional ADER-Cons scheme is significantly more oscillatory than that produced by ADER-Prim. This
is particularly evident in the rest-mass density and in the velocity $v_x$.
We have here a good indication that the ADER-Prim  scheme behaves better than the ADER-Cons  scheme when applied to the equations  of special relativistic magnetohydrodynamics.

%

%

\subsubsection{RMHD Rotor Problem}
\label{sec:RMHD-Rotor-Problem}
\begin{figure}
\begin{center}
\begin{tabular}{cc} 
{\includegraphics[angle=0,width=7.3cm,height=7.3cm]{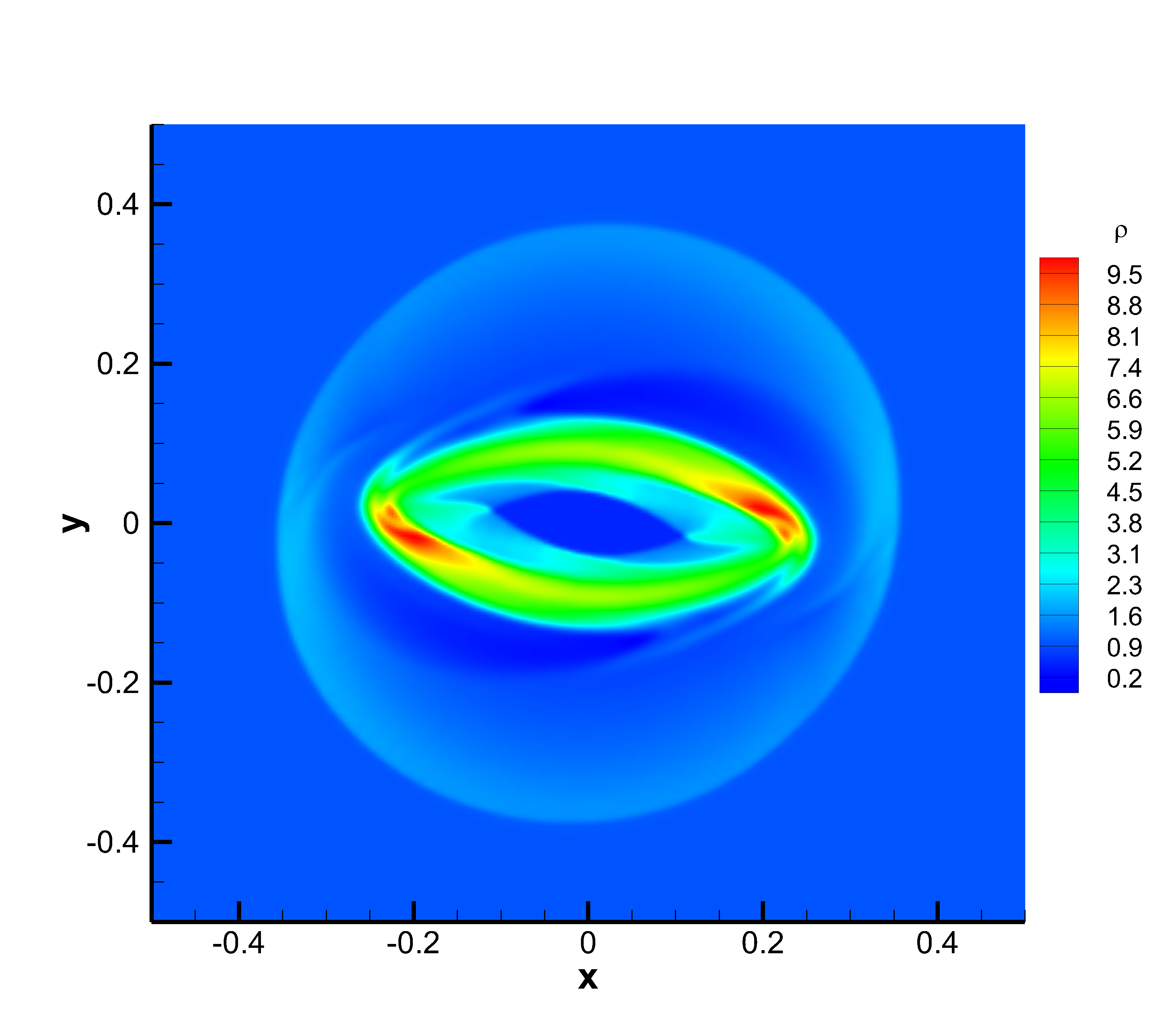}} & 
{\includegraphics[angle=0,width=7.3cm,height=7.3cm]{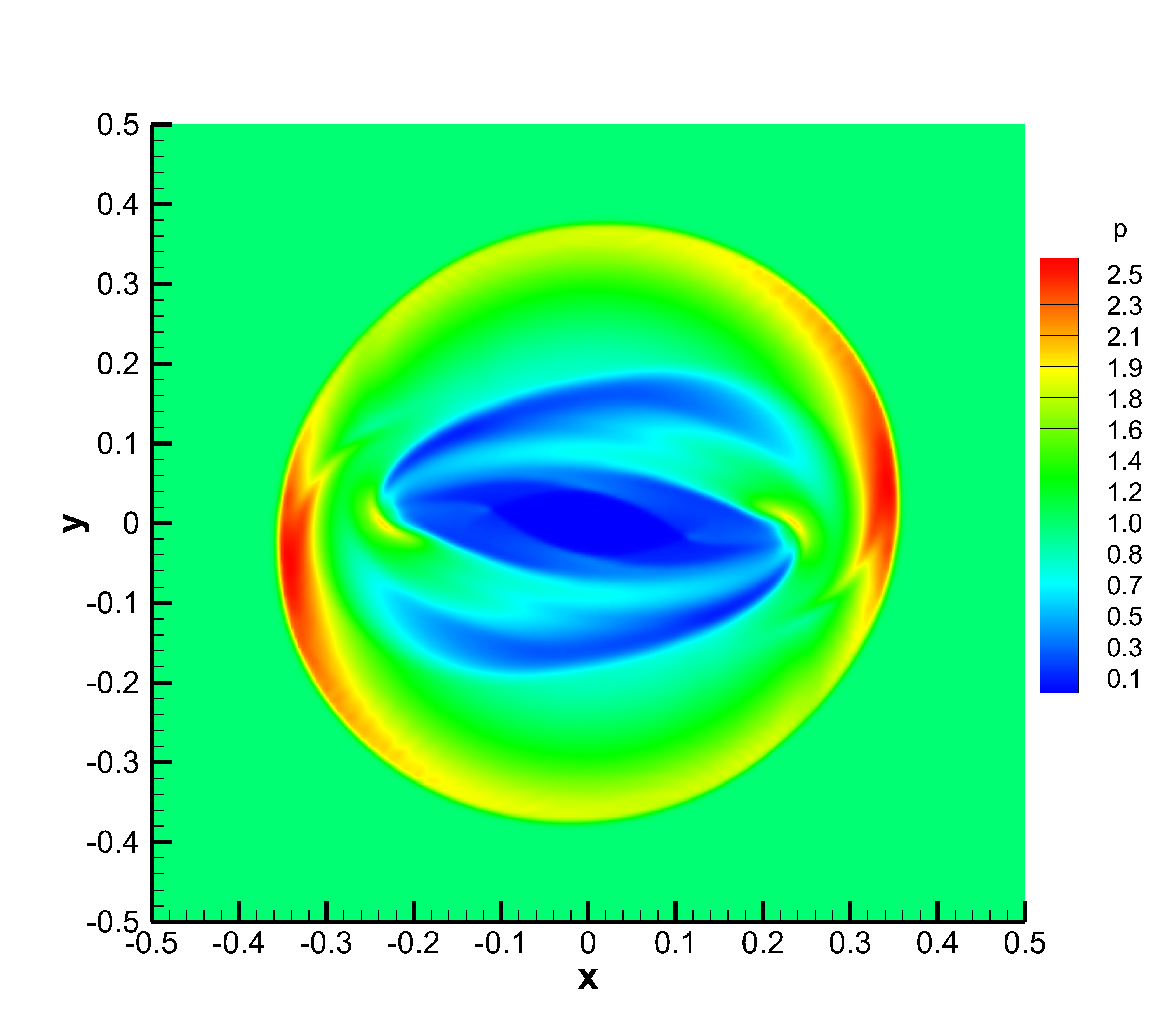}} \\ 
{\includegraphics[angle=0,width=7.3cm,height=7.3cm]{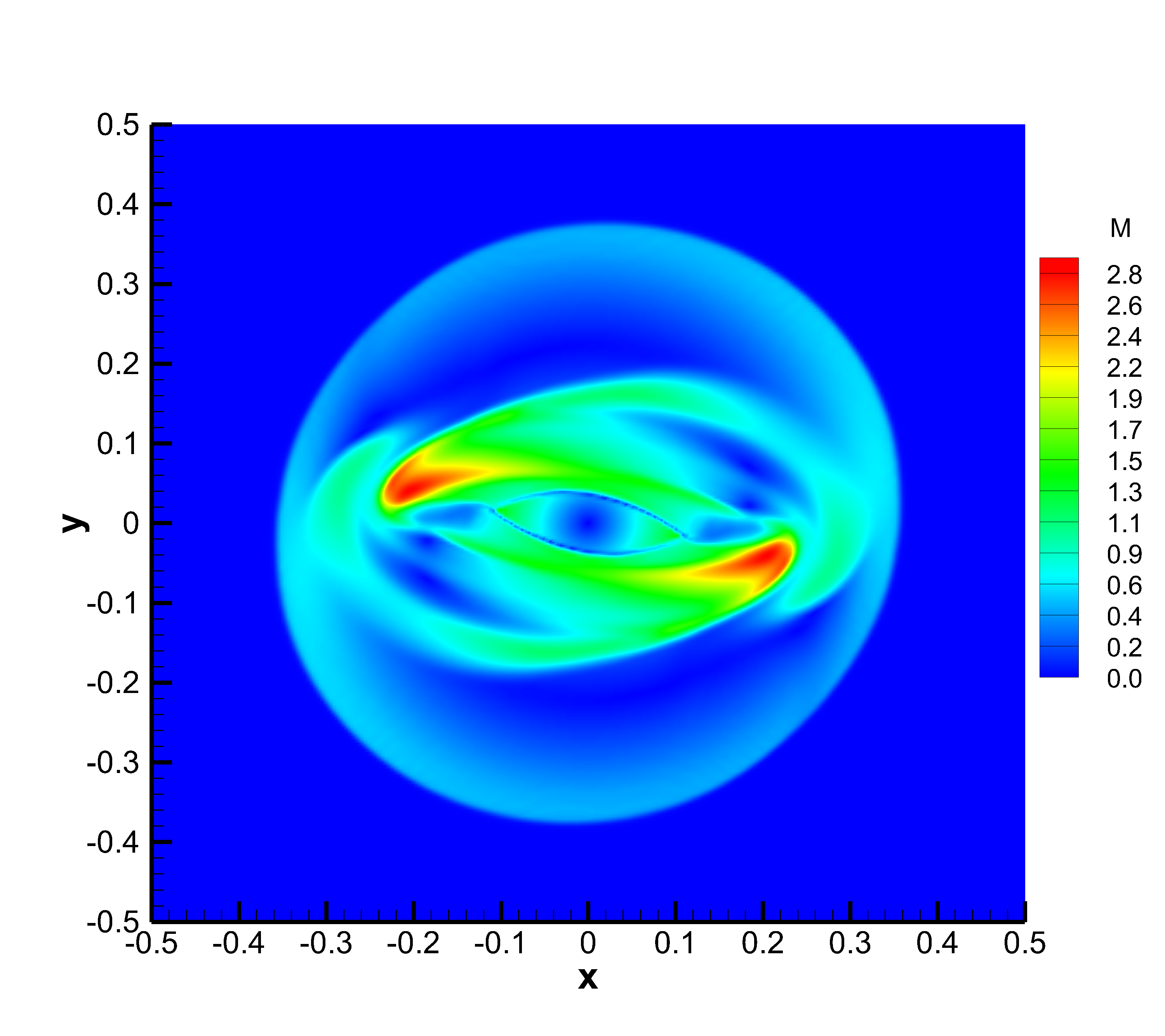}} & 
{\includegraphics[angle=0,width=7.3cm,height=7.3cm]{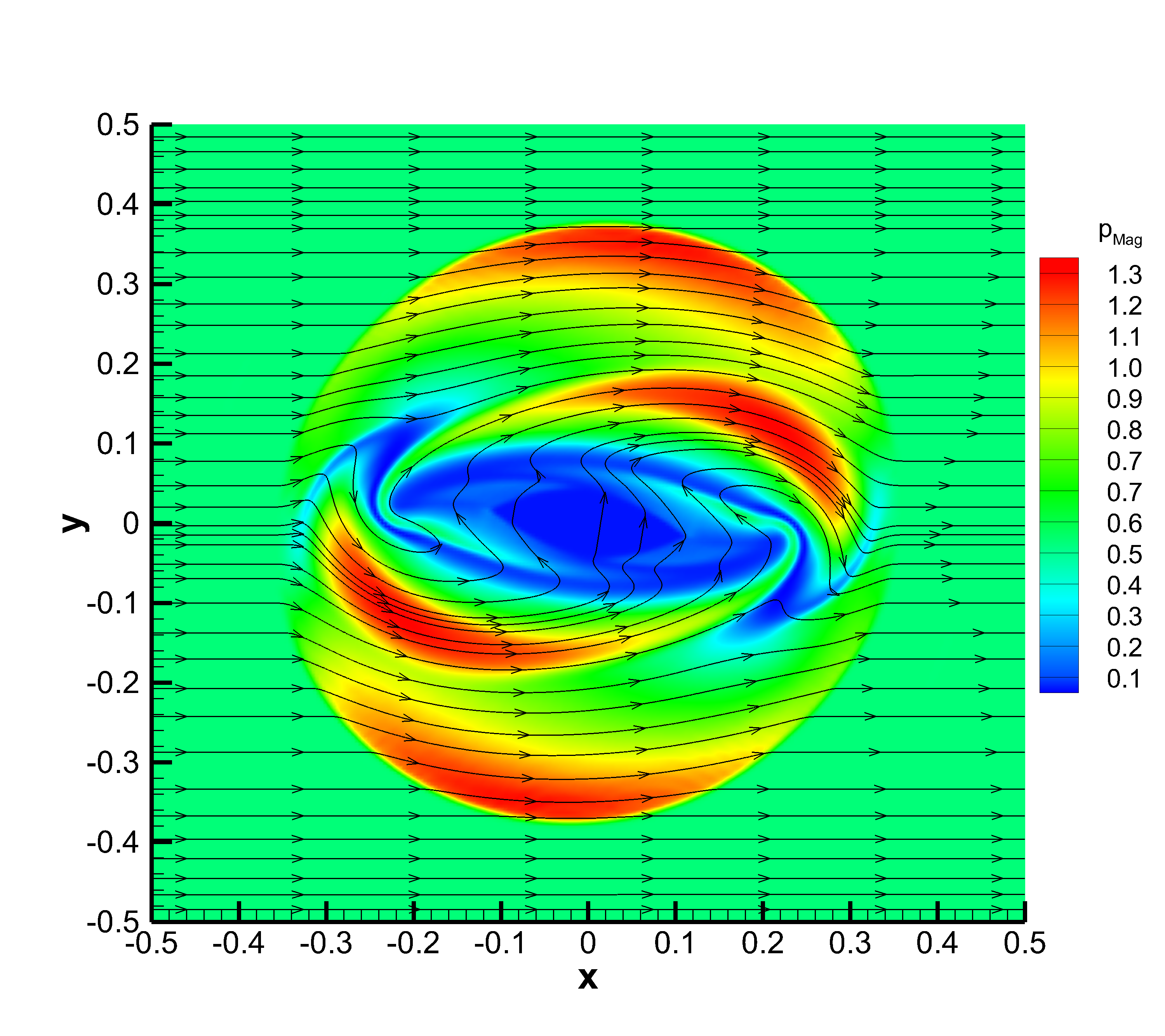}}
\end{tabular} 
\caption{Solution of the RMHD rotor problem at time $t=0.4$, obtained with the $\mathbb{P}_0\mathbb{P}_3$ scheme on a uniform grid with $300\times300$ cells. 
Top panels: rest-mass density (left) and thermal pressure (right). Bottom panels: Mach number (left) and magnetic pressure (right).
}
\label{fig:RMHD-Rotor}
\end{center}
\end{figure}

The relativistic version of the magnetic rotor problem, originally proposed 
by \cite{BalsaraSpicer1999}, has by now become a standard numerical test in RMHD.
It describes the evolution of a high density plasma which, at time $t=0$,
rotates rapidly with angular velocity $\omega$ and is
surrounded by a low density plasma at rest:
\begin{equation}
\rho=\left\{\begin{array}{cl}
10 & \text{for}\;\; 0\le r\le 0.1; \\ 1 & \text{otherwise};
\end{array}\right.,~~~
\omega=\left\{\begin{array}{cl}
9.3 & \text{for}\;\; 0\le r\le 0.1; \\ 0 & \text{otherwise};
\end{array}\right.,~~~
{\B} = \left(\begin{array}{c}
1.0  \\ 0 \\ 0
\end{array}\right),~~~
p = 1\,,\gamma=4/3.
\label{eq:MHDrotor_ic}
\end{equation}
Due to rotation, a sequence of torsional Alfv\'en waves are launched outside the cylinder, with the net effect of reducing the angular velocity of the rotor. We have solved this 
problem over a computational domain  $\Omega = [-0.6,0.6]\times[-0.6,0.6]$, discretized by $300\times300$ numerical cells and
using a fourth order finite volume scheme with the Rusanov Riemann solver.
No taper has been applied to the initial conditions, thus producing true discontinuities right at the beginning.
Fig.~\ref{fig:RMHD-Rotor} shows the rest-mass density, the thermal pressure, the relativistic Mach number and 
the magnetic pressure at time $t=0.4$. 
We obtain results which are in  good qualitative agreement with those available in the literature (see, for instance, 
\cite{DelZanna2003}, \cite{DumbserZanotti}, \cite{ADER_MOOD_14} and \cite{Kim2014}). We
emphasize that for this test the reconstruction of the primitive variables $v^i$ turns out to be more robust than that achieved through the reconstruction of the products $Wv^i$.

%

\subsection{The Baer-Nunziato equations}
\label{sec:BN}

%
As a genuinely non-conservative system of hyperbolic equations we consider the 
Baer-Nunziato model for compressible two-phase flow 
(see also \cite{BaerNunziato1986,SaurelAbgrall,AndrianovWarnecke,Schwendeman,DeledicquePapalexandris,MurroneGuillard}).  
In the rest of the paper we define the
first phase as the solid phase and the second phase as the gas phase. As a result, we will use 
the subscripts $1$ and $s$ as well as $2$ and $g$ as synonyms. 
Sticking to \cite{BaerNunziato1986}, we prescribe the interface velocity $\mathbf{v}_I$ and the pressure  $p_I$ as
$\mathbf{v}_I = \mathbf{v}_1$ and $p_I = p_2$, respectively,
although other choices are also possible \cite{SaurelAbgrall}. 
With these definitions, the system of Baer-Nunziato equations 
can be cast in the form prescribed by (\ref{NCsyst}) after defining the state vector $\u$ as 
\begin{equation}
\u=\left(
\ar, \, \ar v_1^i, \, \ar E_1, \, 
\arr, \, \arr v_2^i, \, \arr E_2, \, \phi_1 
\right)\,,
\end{equation}
where $\phi_k$ is the volume fraction of phase $k$, with the condition that $\phi_1+\phi_2=1$.
On the other hand, the fluxes ${\bf f}^i$, the sources $\bf S$ and the non-conservative matrices 
${\bf B}_i$ are expressed by
\be
{\bf f}^i=\left[\begin{array}{c}
\ar v_1^i \\ \phi_1( \rho_1 v_1^i v_1^j + p_1\delta^{ij} ) \\ \phi_1 v_1^i(\rho_1 E_1+p_1) \\ \phi_2\rho_2 v_2^i \\ \phi_2( \rho_2 v_2^i v_2^j + p_2\delta^{ij} )\\
\phi_2 v_2^i(\rho_2 E_2 + p_2) \\ 0 \\
\end{array}\right],~~~
{\bf S}=\left[\begin{array}{c}
0 \\ -\nu (v_1^i - v_2^i) \\ -\nu \v_1 \cdot (\v_1 - \v_2) \\ 0 \\ -\nu (v_2^i - v_1^i)\\
-\nu \v_1 \cdot (\v_2 - \v_1) \\ \mu(p_1-p_2) \\
\end{array}\right]\,,
\label{eq:bnsource}
\ee
\begin{table}[!t]
\begin{center}
\renewcommand{\arraystretch}{1.0}
\begin{tabular}{ccccccccc}
\hline
   & $\rho_s$ & $u_s$  & $p_s$ & $\rho_g$ & $u_g$ & $p_g$ & $\phi_s$ & $t_e$  \\
\hline 
\multicolumn{1}{l}{\textbf{BNRP1 \cite{DeledicquePapalexandris}:} } & 
\multicolumn{8}{c}{ $\gamma_s = 1.4, \quad \pi_s = 0, \quad \gamma_g = 1.4, \quad \pi_g = 0$}  \\
\hline 
L & 1.0    & 0.0   & 1.0  & 0.5 & 0.0   &  1.0 & 0.4 & 0.10 \\
R & 2.0    & 0.0   & 2.0  & 1.5 & 0.0   &  2.0 & 0.8 &      \\
\hline 
\multicolumn{1}{l}{\textbf{BNRP2 \cite{DeledicquePapalexandris}:}} & 
\multicolumn{8}{c}{ $\gamma_s = 3.0, \quad \pi_s = 100, \quad \gamma_g = 1.4, \quad \pi_g = 0$}  \\
\hline
L & 800.0   & 0.0   & 500.0  & 1.5 & 0.0   & 2.0 & 0.4 & 0.10  \\
R & 1000.0  & 0.0   & 600.0  & 1.0 & 0.0   & 1.0 & 0.3 &       \\
\hline 
\multicolumn{1}{l}{\textbf{BNRP3 \cite{DeledicquePapalexandris}:}} & 
\multicolumn{8}{c}{ $\gamma_s = 1.4, \quad \pi_s = 0, \quad \gamma_g = 1.4, \quad \pi_g = 0$}  \\ 
\hline
L & 1.0     & 0.9       & 2.5      & 1.0       & 0.0      &  1.0 & 0.9   & 0.10   \\
R & 1.0     & 0.0       & 1.0      & 1.2       & 1.0      &  2.0 & 0.2   &        \\
\hline 
\multicolumn{1}{l}{\textbf{BNRP5 \cite{Schwendeman}:}} & 
\multicolumn{8}{c}{ $\gamma_s = 1.4, \quad \pi_s = 0, \quad \gamma_g = 1.4, \quad \pi_g = 0$}  \\
\hline 
L & 1.0    & 0.0   & 1.0  & 0.2 & 0.0   &  0.3 & 0.8 & 0.20  \\
R & 1.0    & 0.0   & 1.0  & 1.0 & 0.0   &  1.0 & 0.3 &      \\
\hline 
\multicolumn{1}{l}{\textbf{BNRP6 \cite{AndrianovWarnecke}:}} & 
\multicolumn{8}{c}{ $\gamma_s = 1.4, \quad \pi_s = 0, \quad \gamma_g = 1.4, \quad \pi_g = 0$}  \\
\hline 
L & 0.2068    & 1.4166   & 0.0416  & 0.5806 & 1.5833    &  1.375 & 0.1 & 0.10  \\
R & 2.2263    & 0.9366   & 6.0     & 0.4890 & -0.70138  &  0.986 & 0.2 &       \\
\hline
\end{tabular}
\end{center}
\caption{Initial states left (L) and right (R) for the Riemann problems for the Baer-Nunziato equations.
Values for $\gamma_i$, $\pi_i$ and the final time $t_e$ are also reported.} 
\label{tab.rpbn.ic}
\end{table}
\be
{\bf{B}}_i = \left( {\begin{array}{*{20}{l}}
0&0&0&0&0&0&0&0&0&0&0\\
0&0&0&0&0&0&0&0&0&0& - p_I \mathbf{e}_i \\
0&0&0&0&0&0&0&0&0&0& - p_I v_I^i\\
0&0&0&0&0&0&0&0&0&0&0\\
0&0&0&0&0&0&0&0&0&0&p_I \mathbf{e}_i \\
0&0&0&0&0&0&0&0&0&0&p_I v_I^i\\
0&0&0&0&0&0&0&0&0&0&v_I^i
\end{array}} \right), 
\ee

%
where $\mathbf{e}_i$ is the unit vector pointing in direction $i$,  ($i \in \left\{x,y,z\right\}$)  
and $\nu$ and $\mu$ are two parameters related to the friction between the phases and to the pressure relaxation.
\footnote{In the tests below $\nu$ and $\mu$ are both set to zero.} 

\begin{figure}[!htbp]
\begin{center}
\begin{tabular}{cc} 
\includegraphics[width=0.45\textwidth]{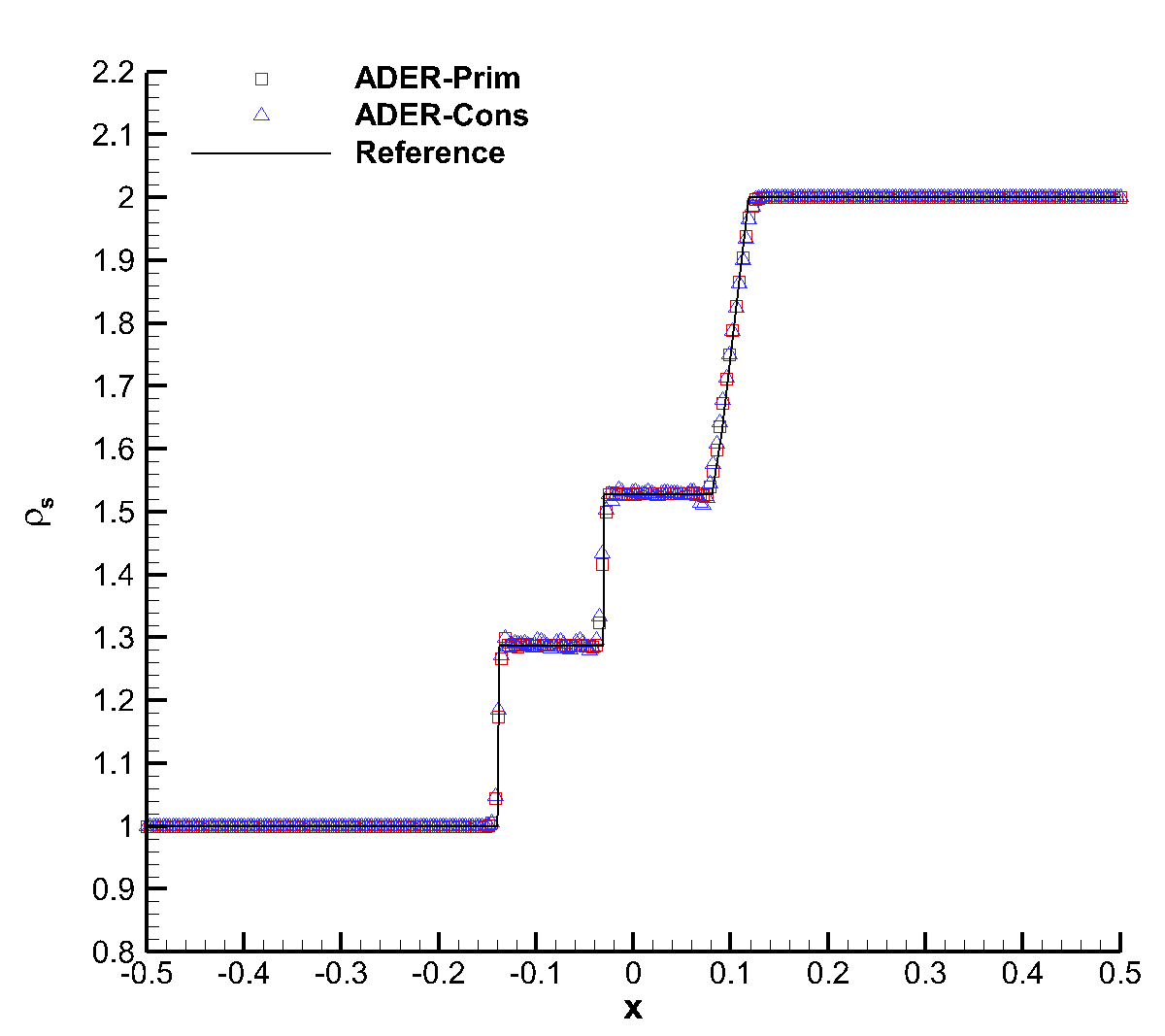}    & 
\includegraphics[width=0.45\textwidth]{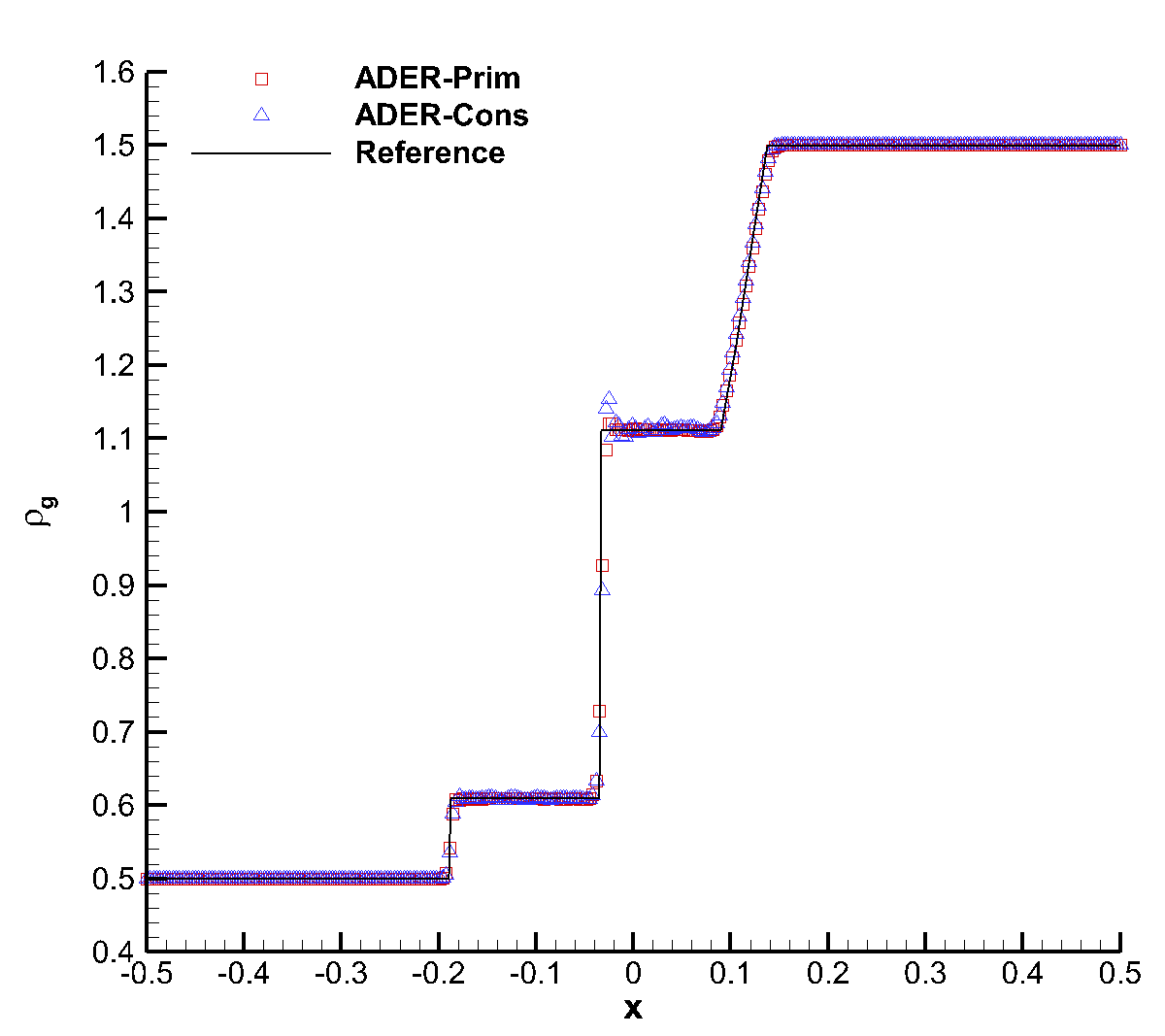}    \\ 
\includegraphics[width=0.45\textwidth]{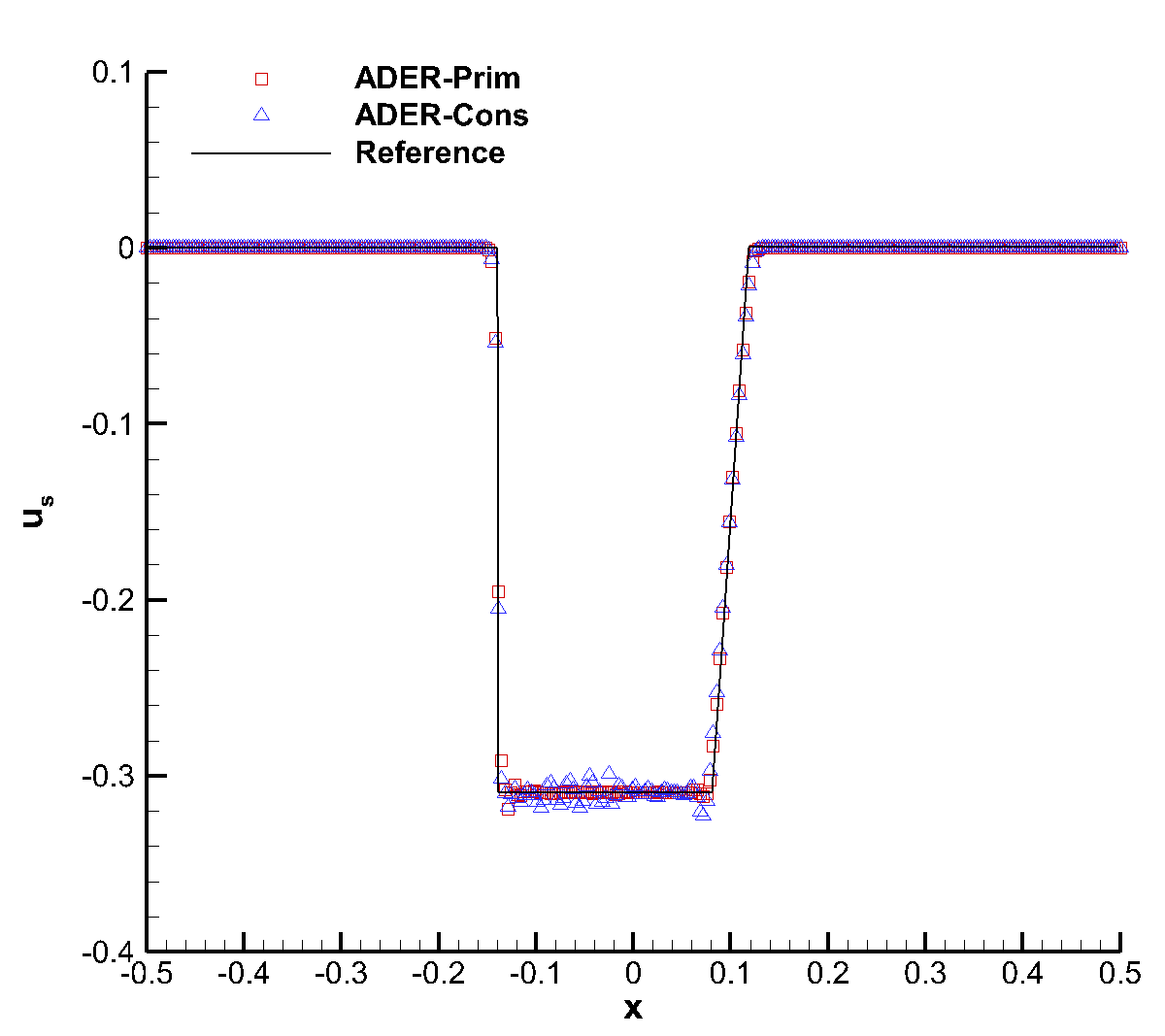}      &  
\includegraphics[width=0.45\textwidth]{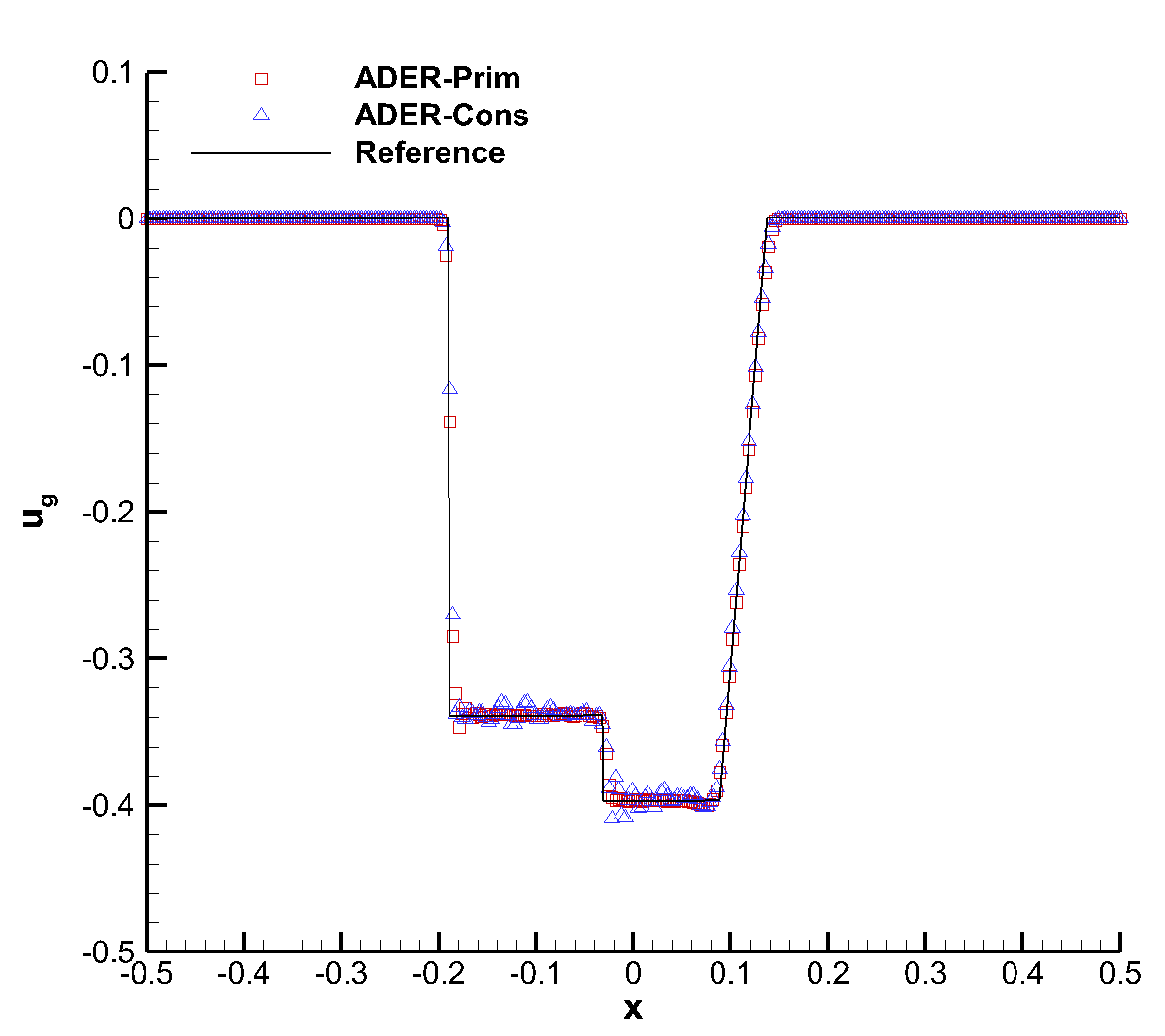}      \\ 
\includegraphics[width=0.45\textwidth]{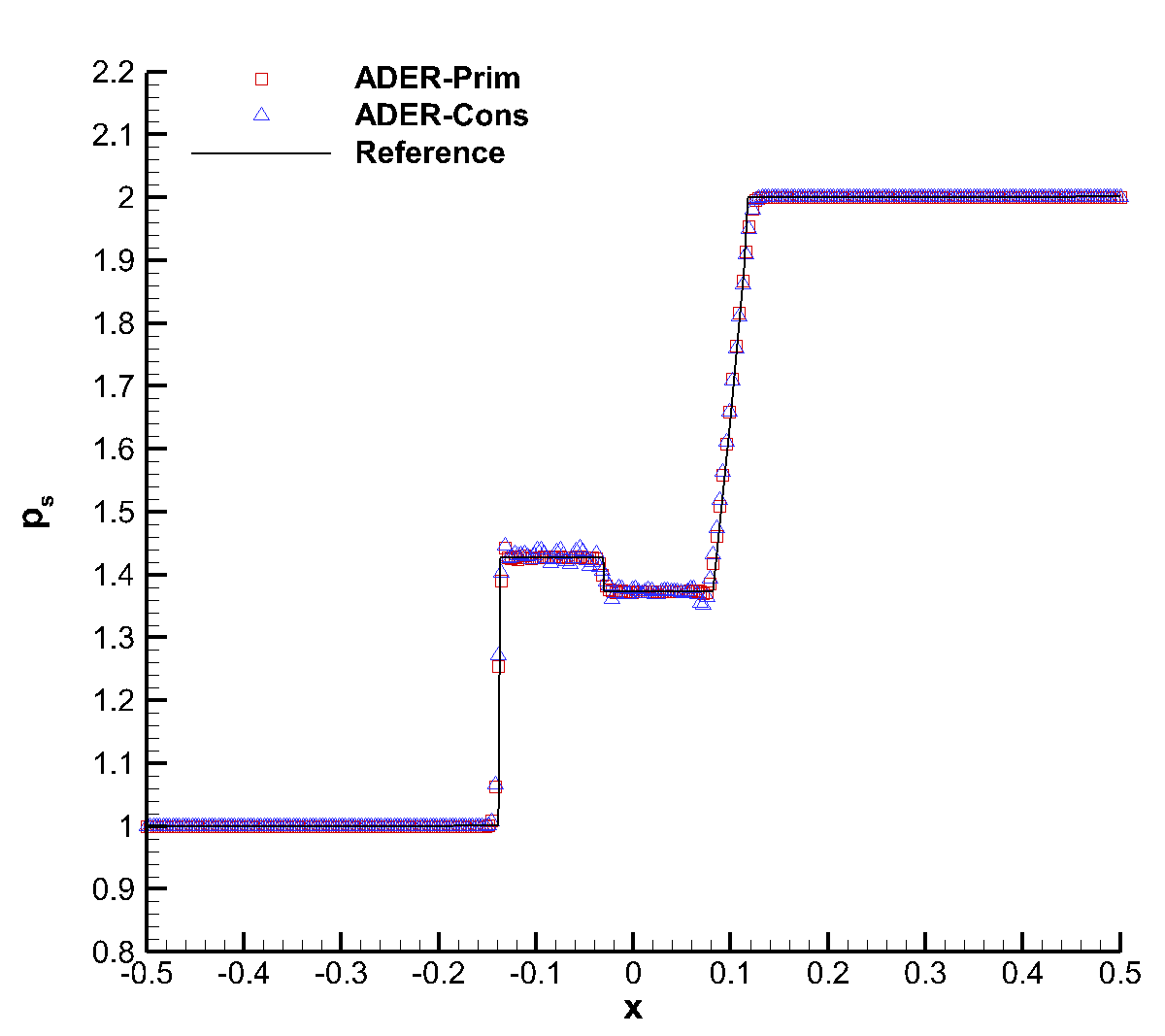}      & 
\includegraphics[width=0.45\textwidth]{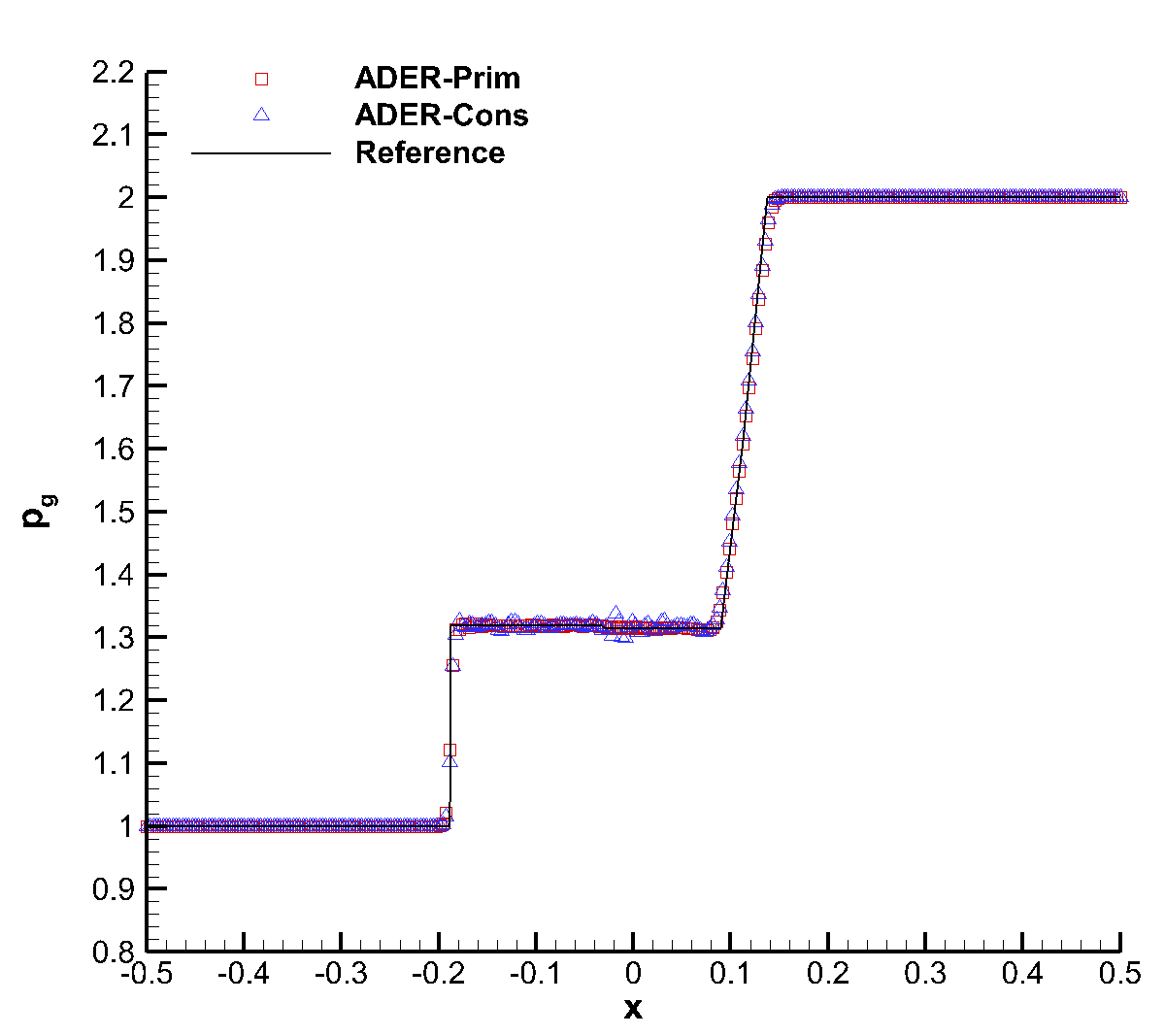}   
\end{tabular}
\caption{Results for the Baer--Nunziato Riemann problem BNRP1. The Osher Riemann solver has been used over a $300$ cells uniform grid.}
\label{fig.bn.rp1}
\end{center}
\end{figure}
\begin{figure}[!htbp]
\begin{center}
\begin{tabular}{cc} 
\includegraphics[width=0.45\textwidth]{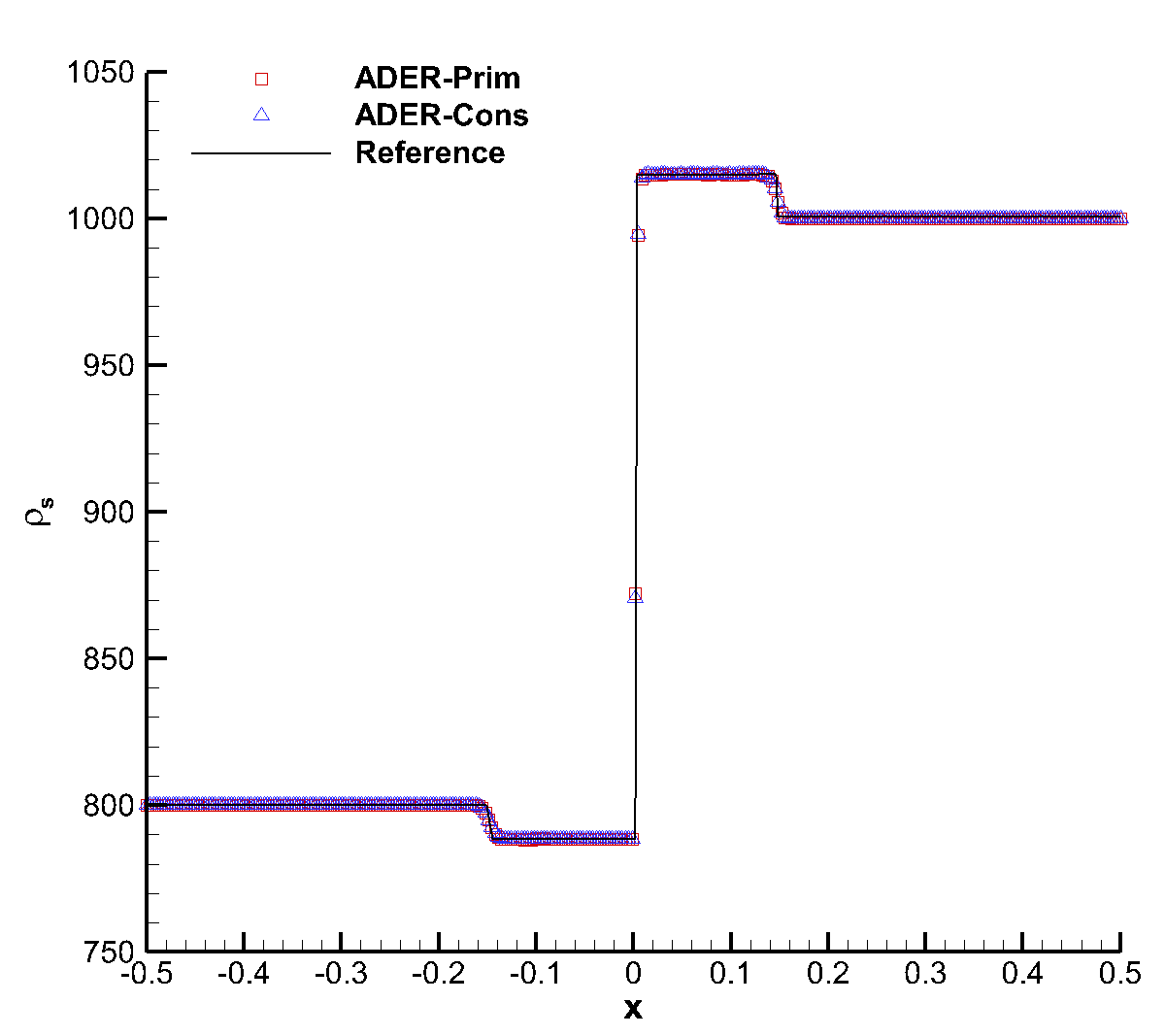}    & 
\includegraphics[width=0.45\textwidth]{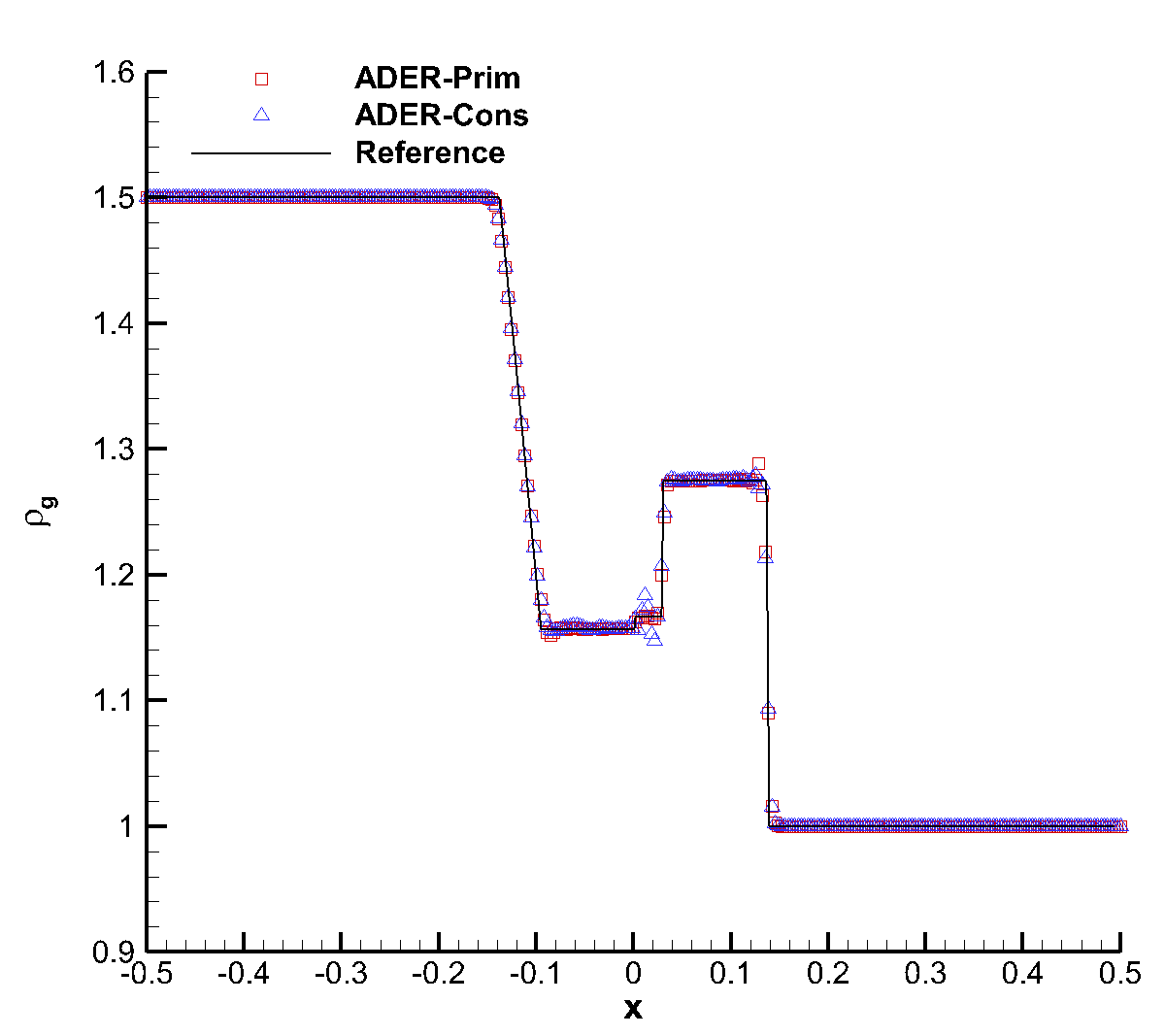}    \\ 
\includegraphics[width=0.45\textwidth]{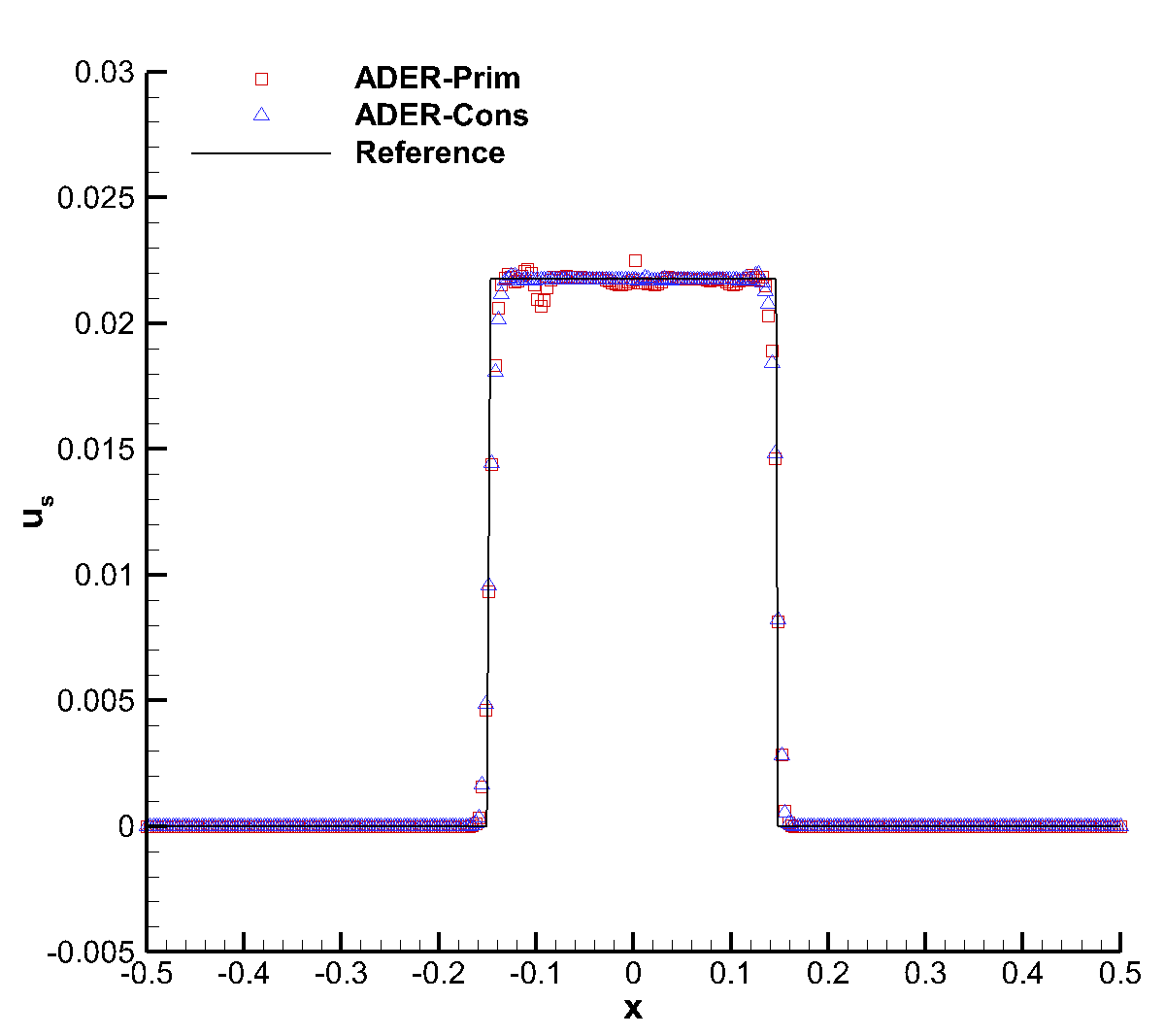}      &  
\includegraphics[width=0.45\textwidth]{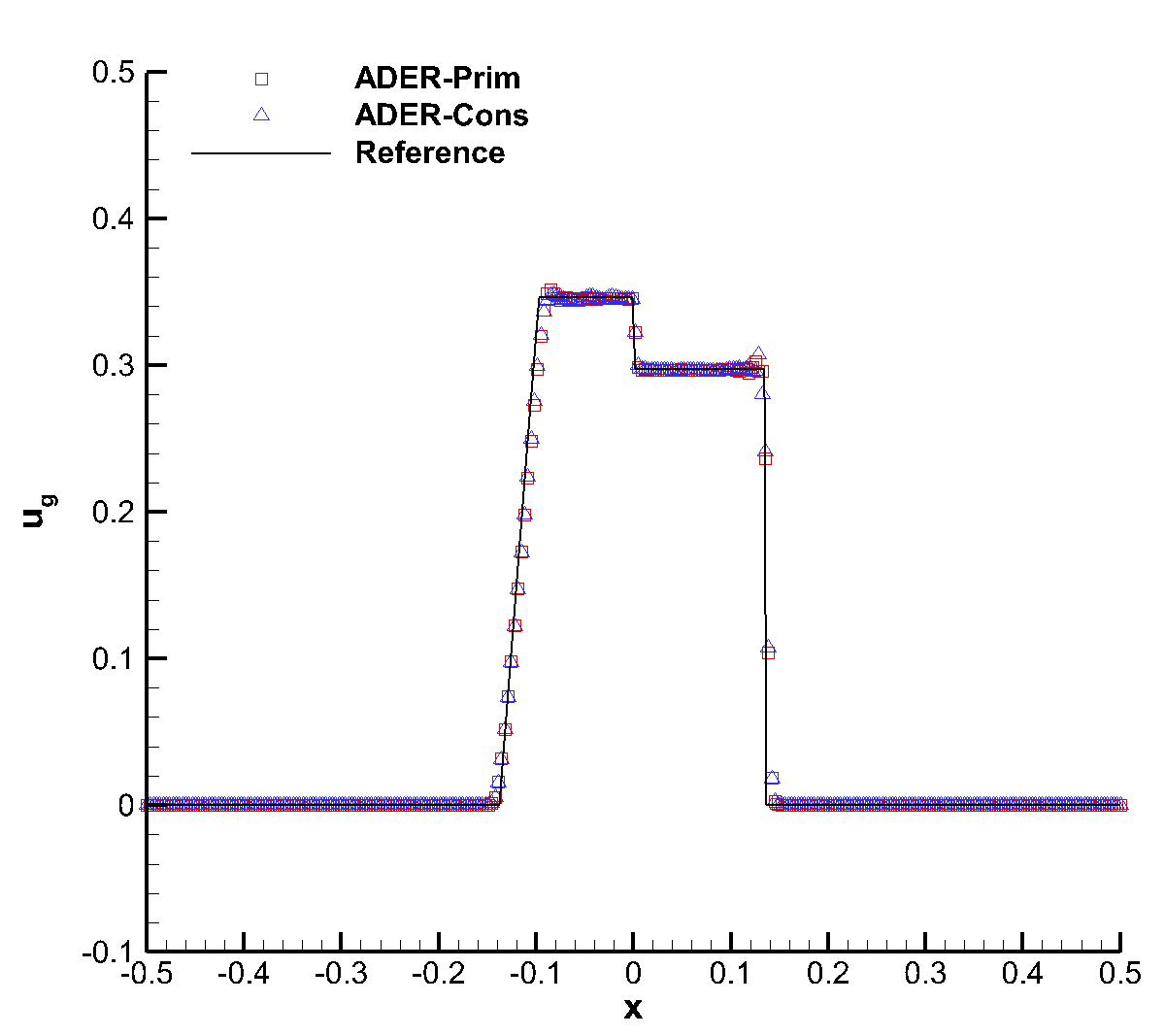}      \\ 
\includegraphics[width=0.45\textwidth]{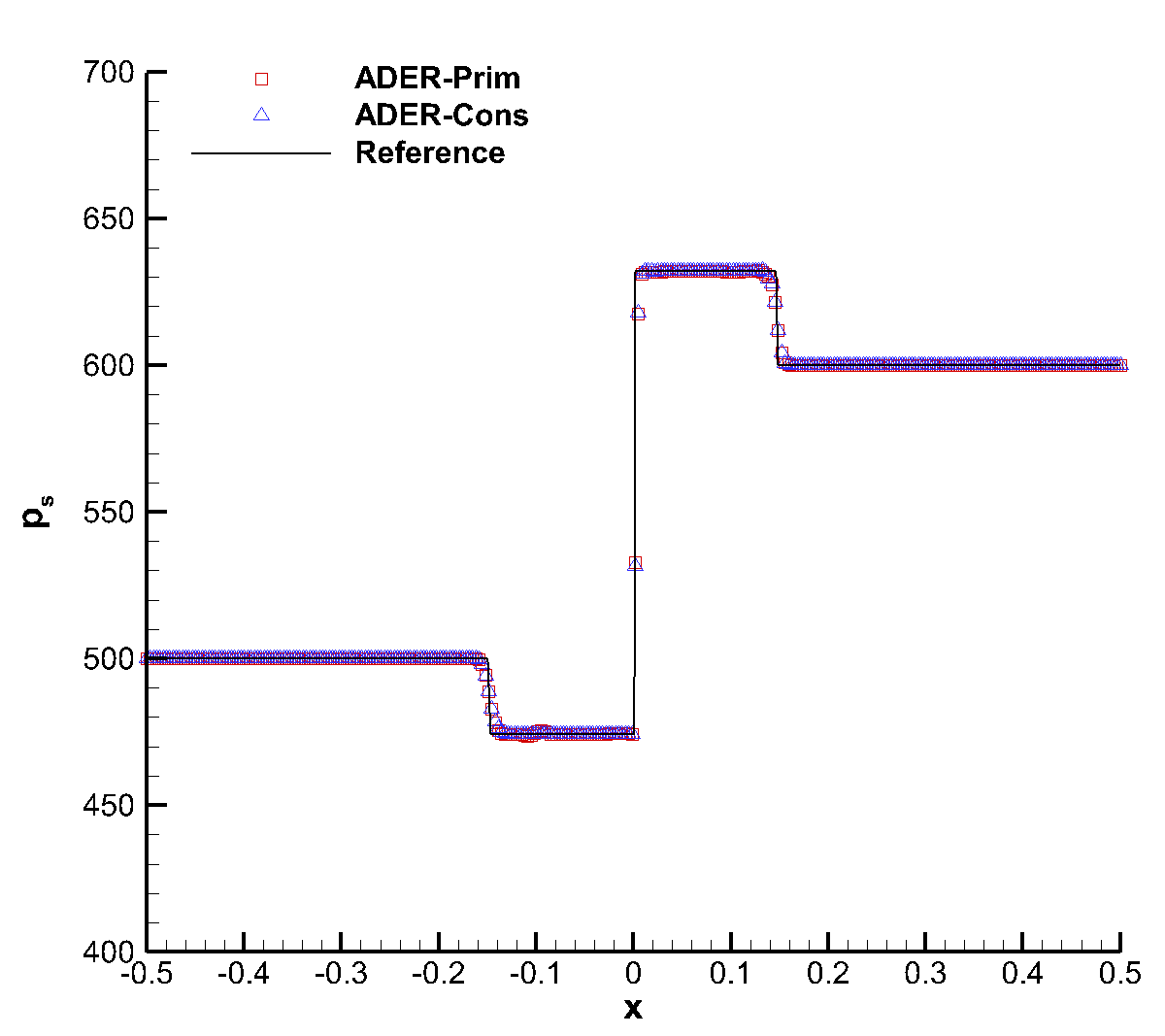}      & 
\includegraphics[width=0.45\textwidth]{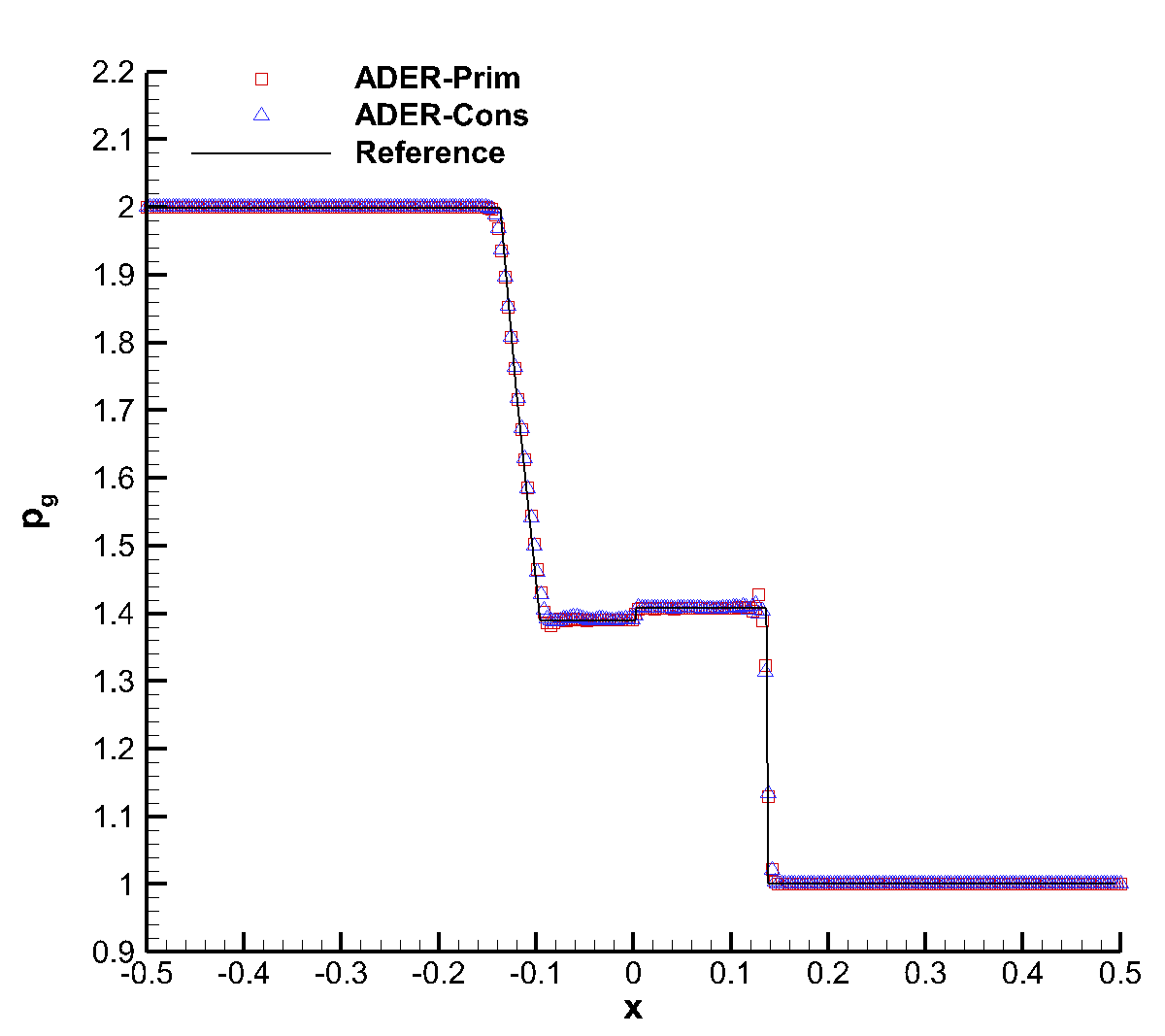}   
\end{tabular}
\caption{Results for the Baer--Nunziato Riemann problem BNRP2. The Osher Riemann solver has been used over a $300$ cells uniform grid.}
\label{fig.bn.rp2}
\end{center}
\end{figure}
\begin{figure}[!htbp]
\begin{center}
\begin{tabular}{cc} 
\includegraphics[width=0.45\textwidth]{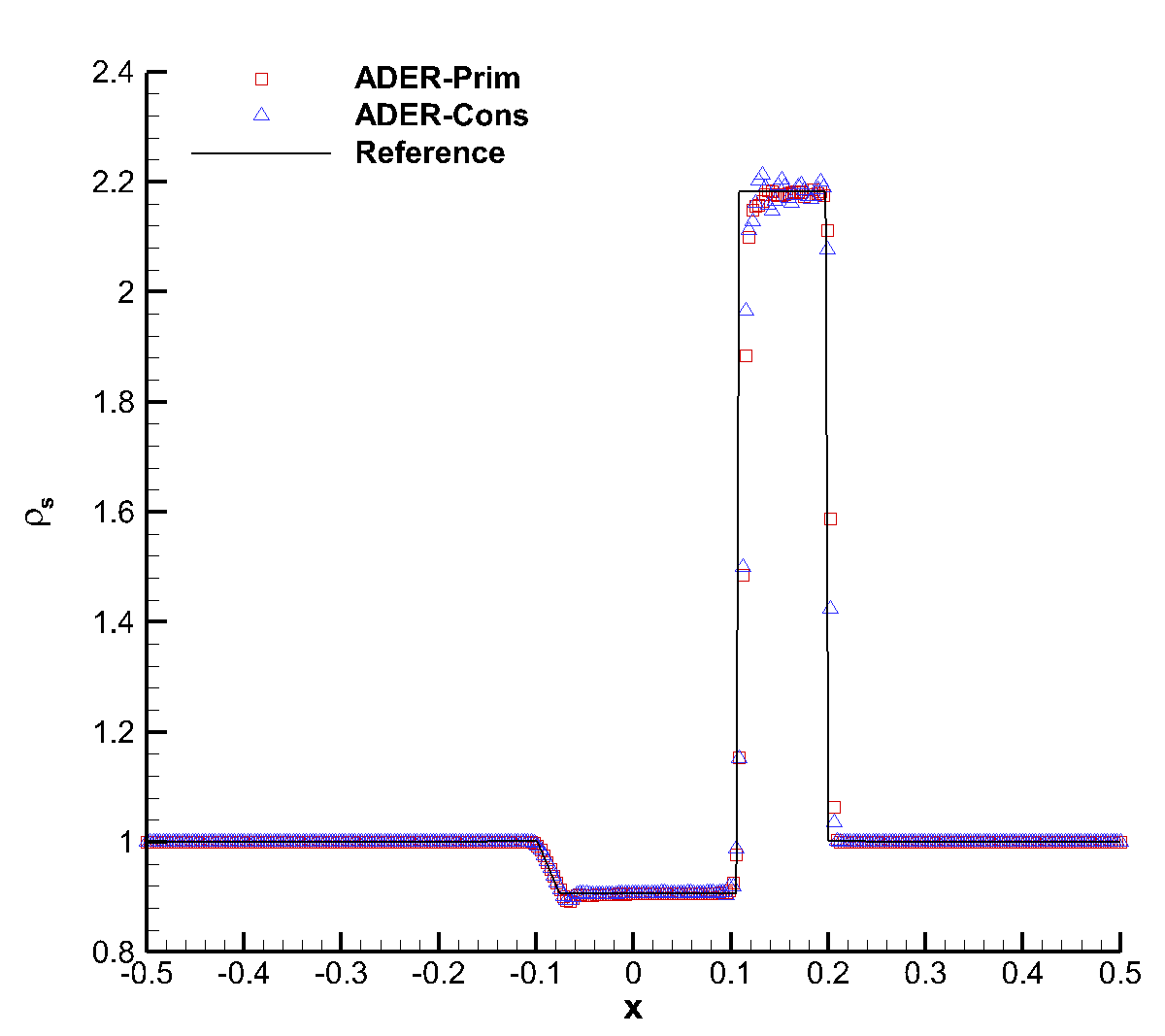}    & 
\includegraphics[width=0.45\textwidth]{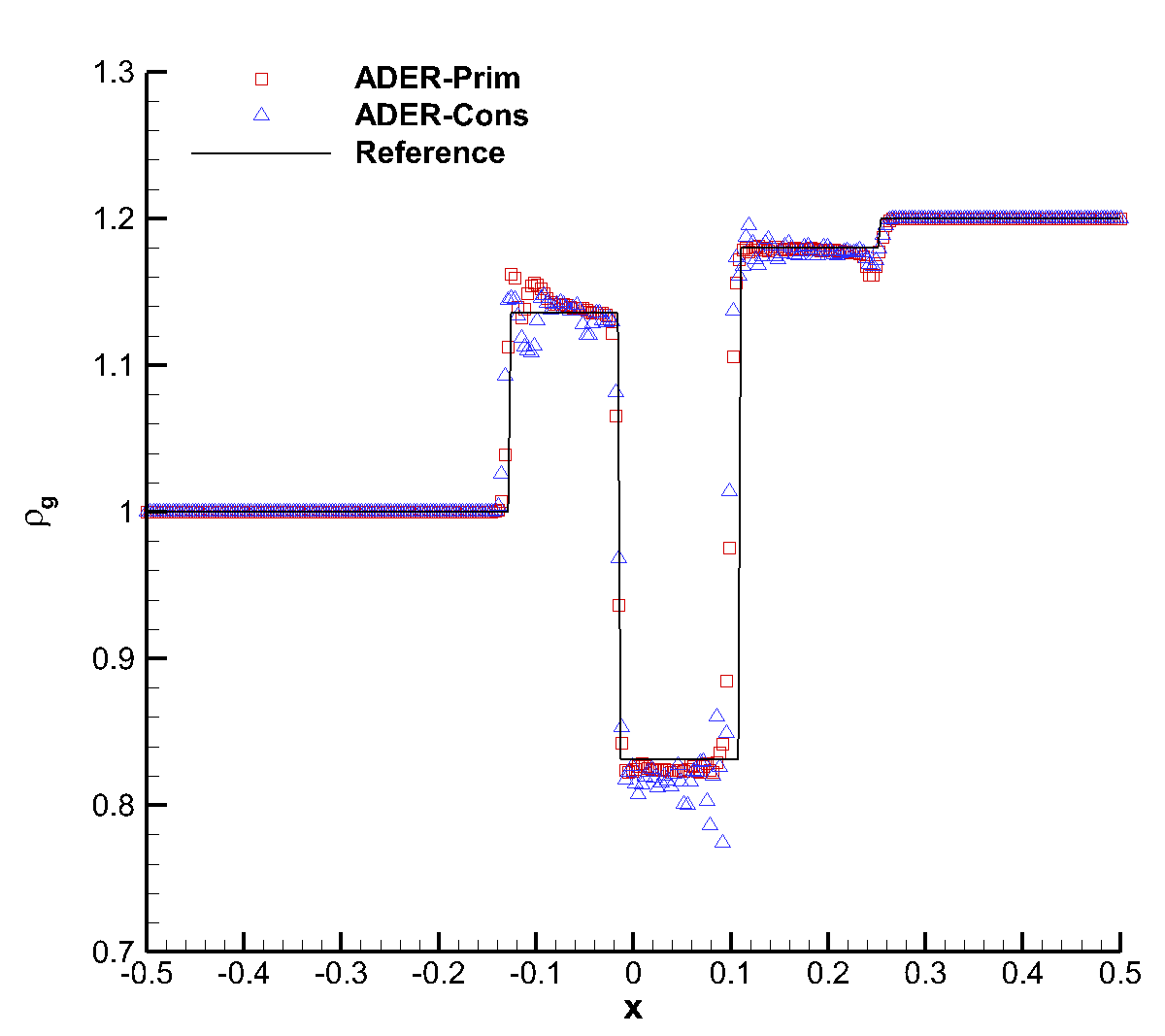}    \\ 
\includegraphics[width=0.45\textwidth]{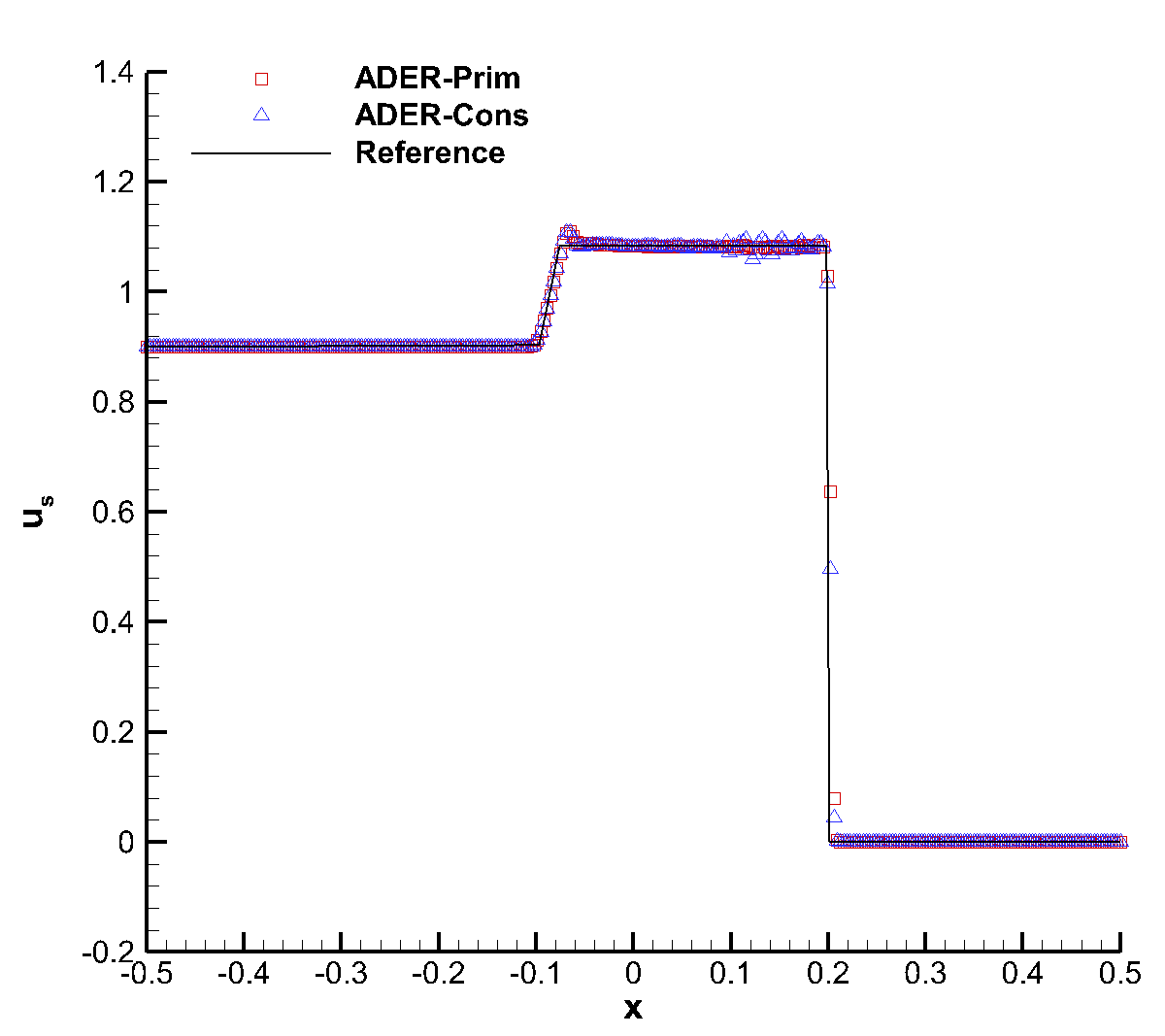}      &  
\includegraphics[width=0.45\textwidth]{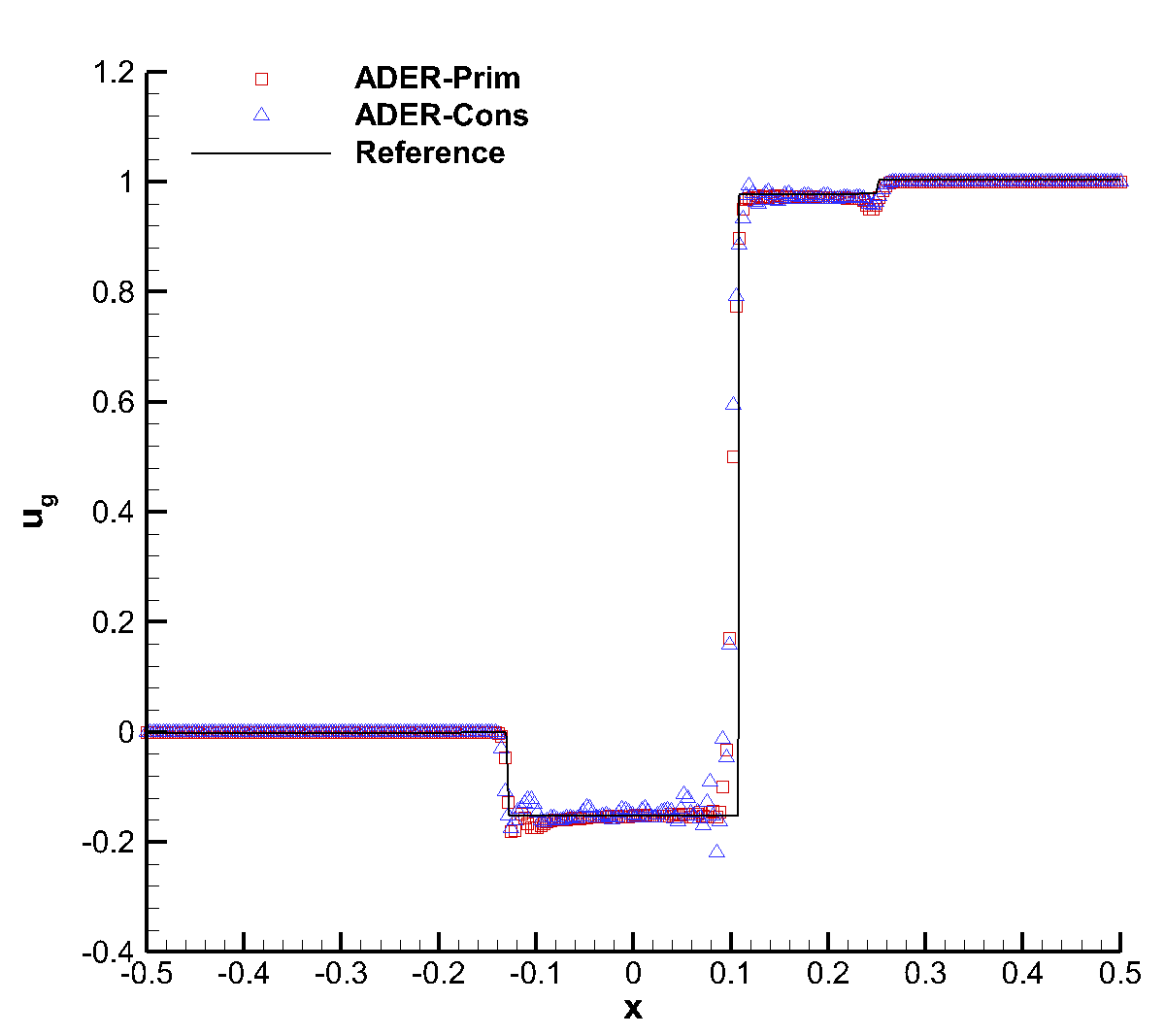}      \\ 
\includegraphics[width=0.45\textwidth]{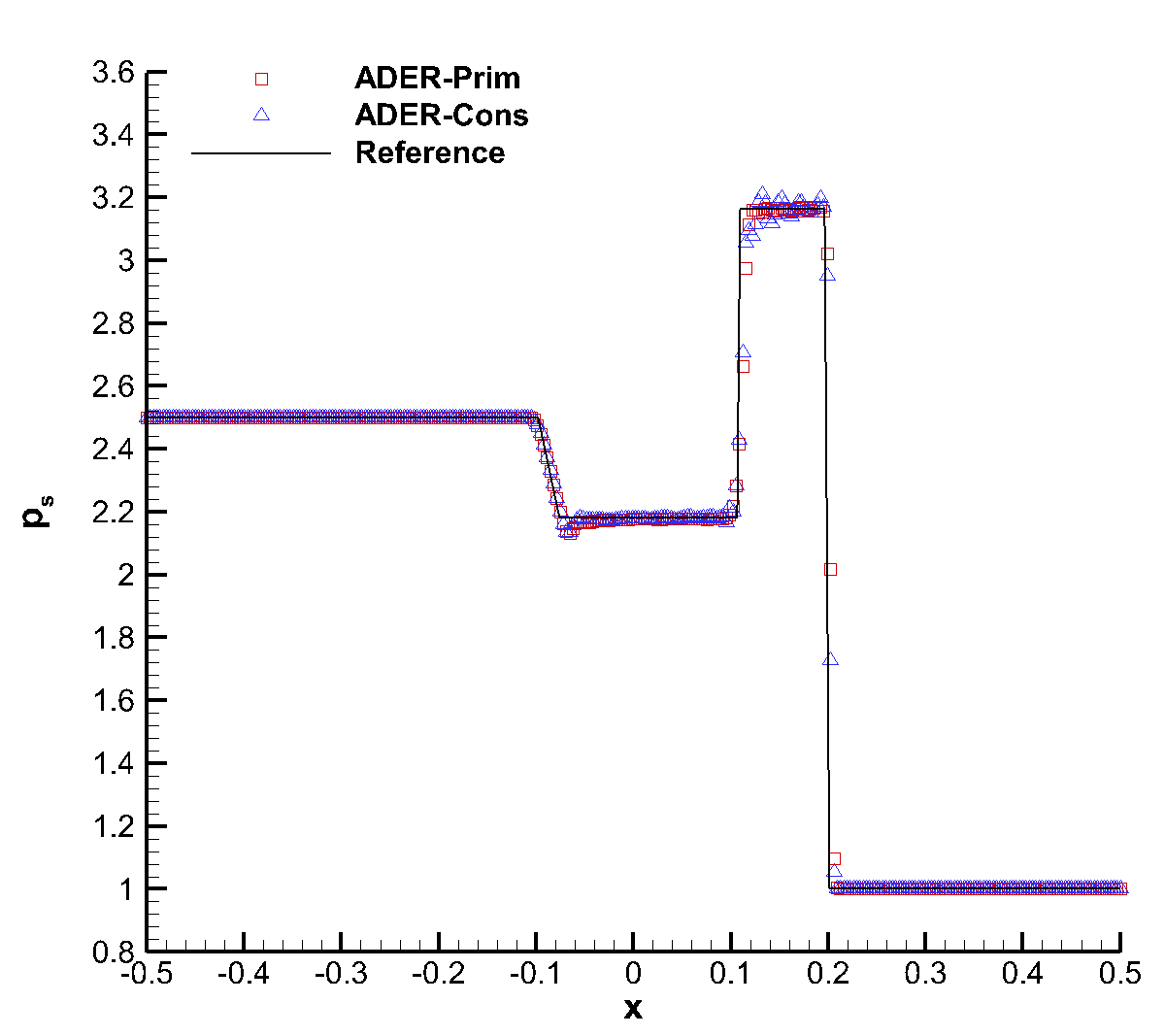}      & 
\includegraphics[width=0.45\textwidth]{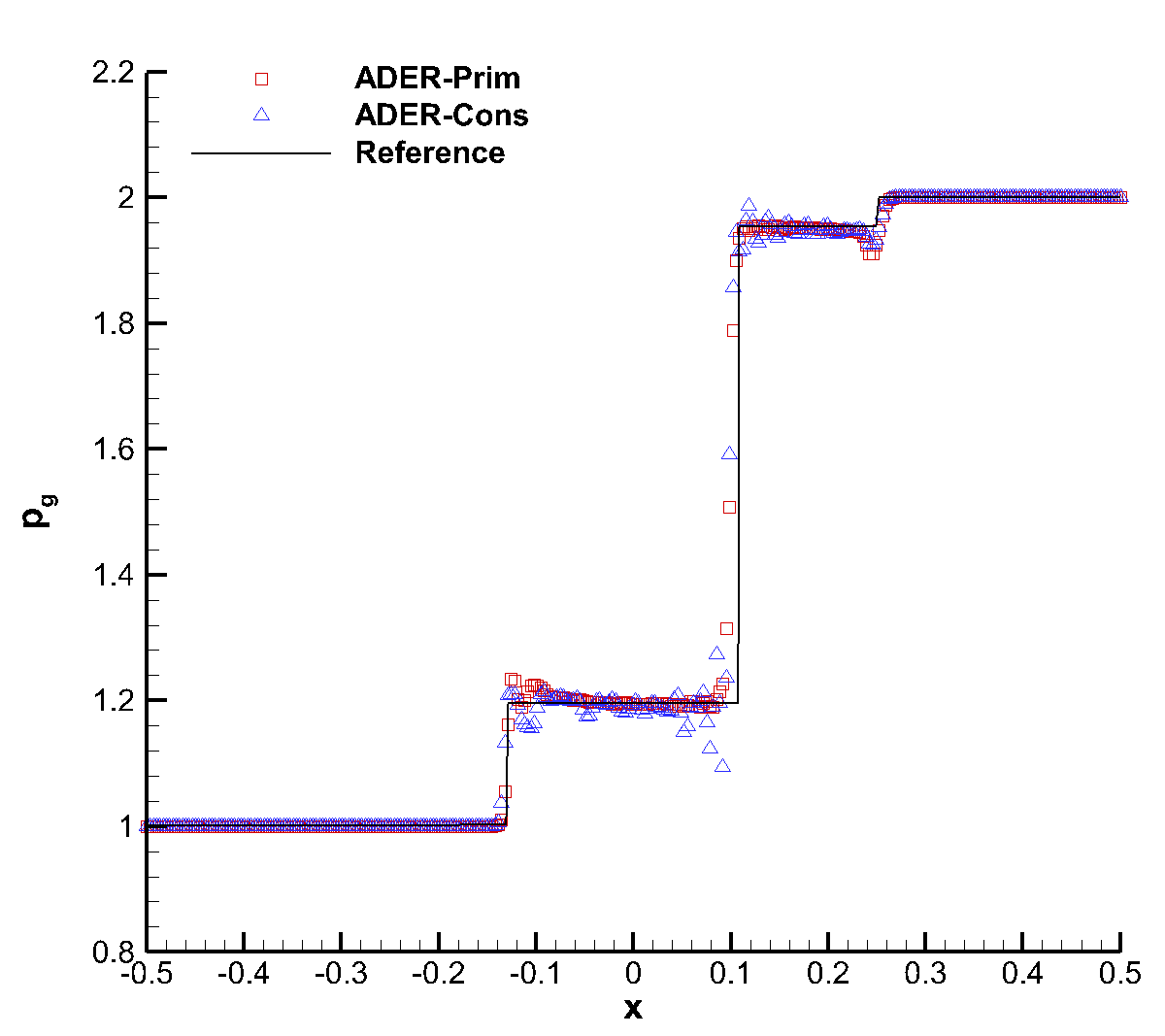}   
\end{tabular}
\caption{Results for the Baer--Nunziato Riemann problem BNRP3. The Osher Riemann solver has been used over a $300$ cells uniform grid.}
\label{fig.bn.rp3}
\end{center}
\end{figure}
\begin{figure}[!htbp]
\begin{center}
\begin{tabular}{cc} 
\includegraphics[width=0.45\textwidth]{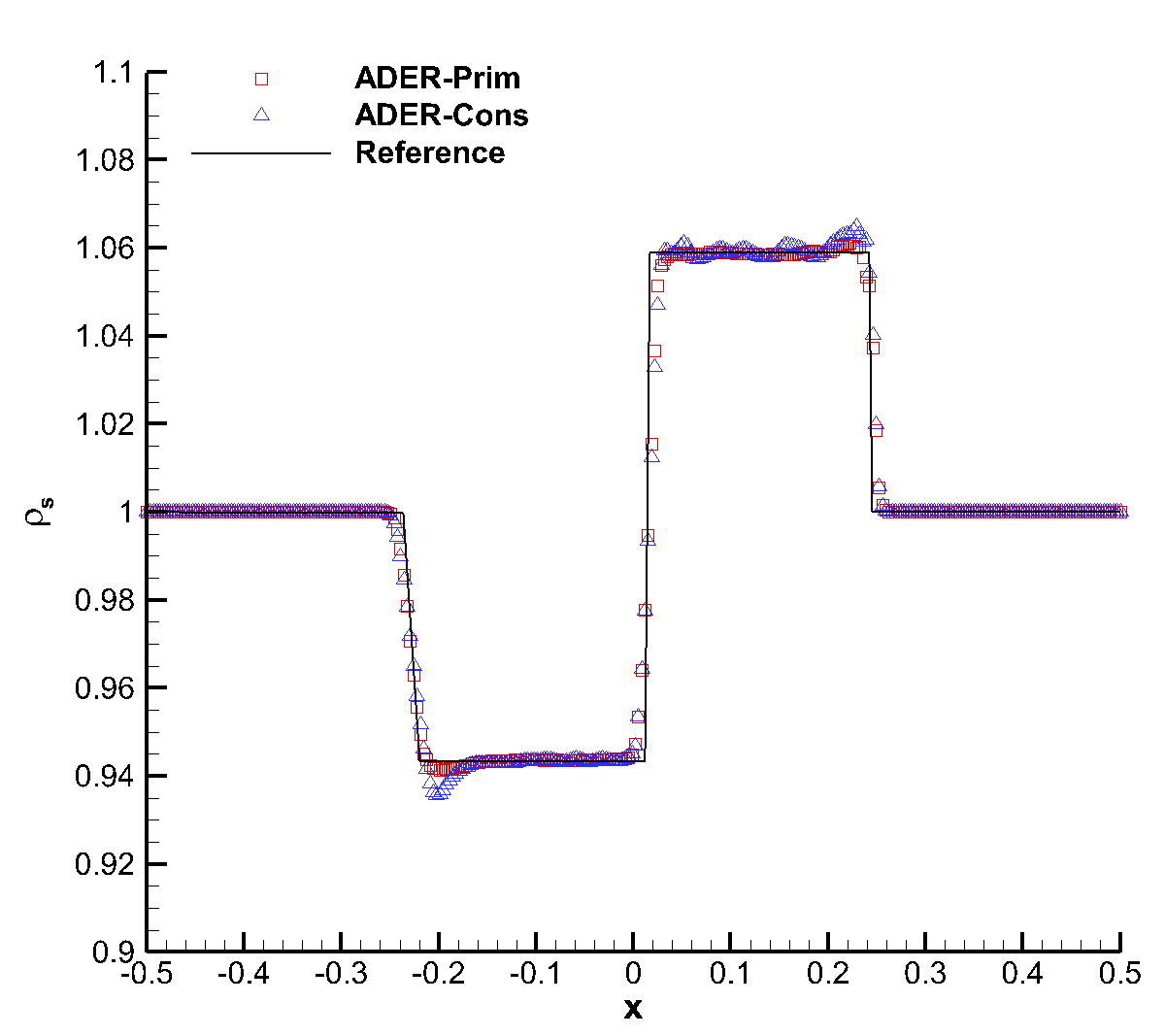}    & 
\includegraphics[width=0.45\textwidth]{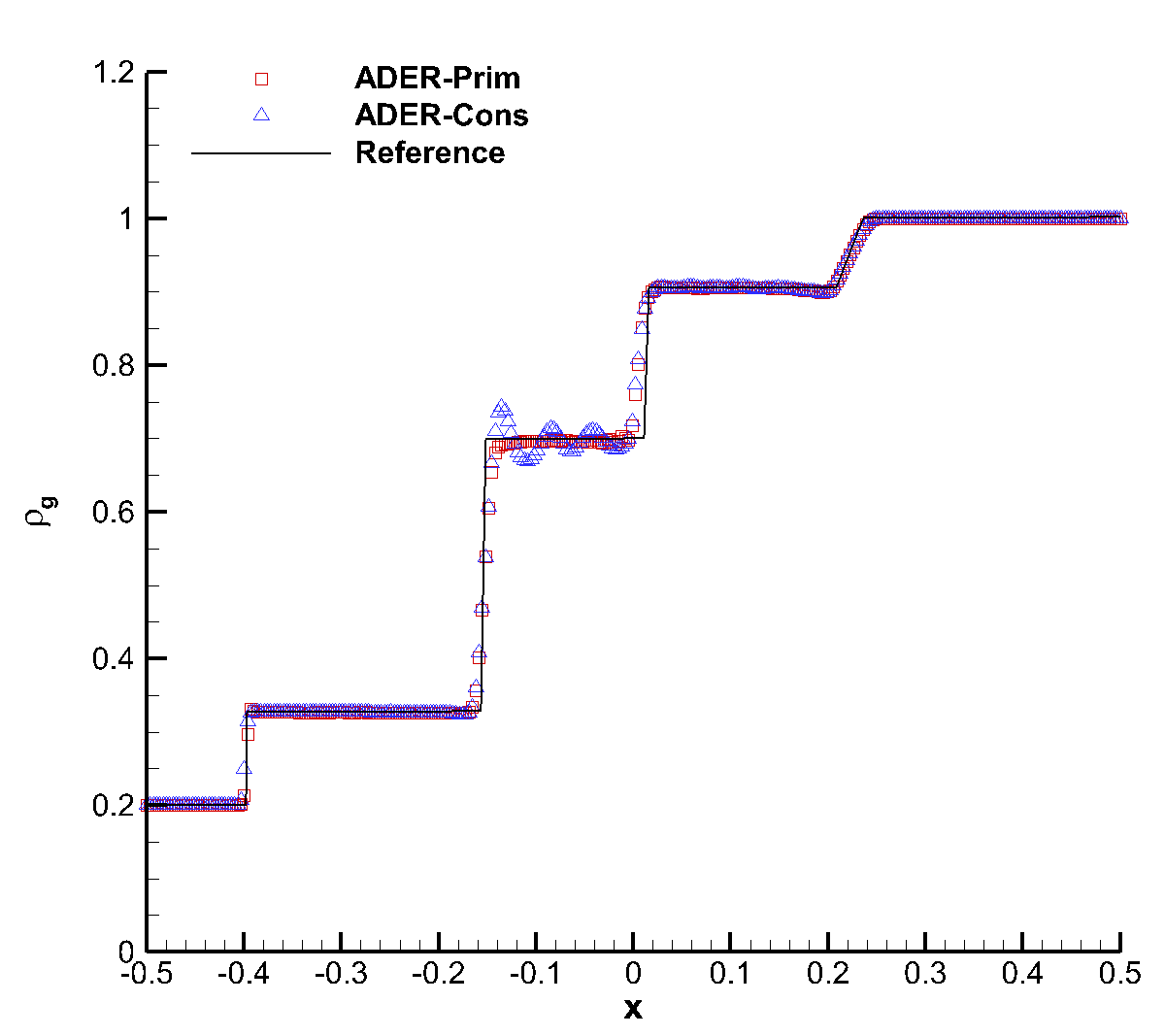}    \\ 
\includegraphics[width=0.45\textwidth]{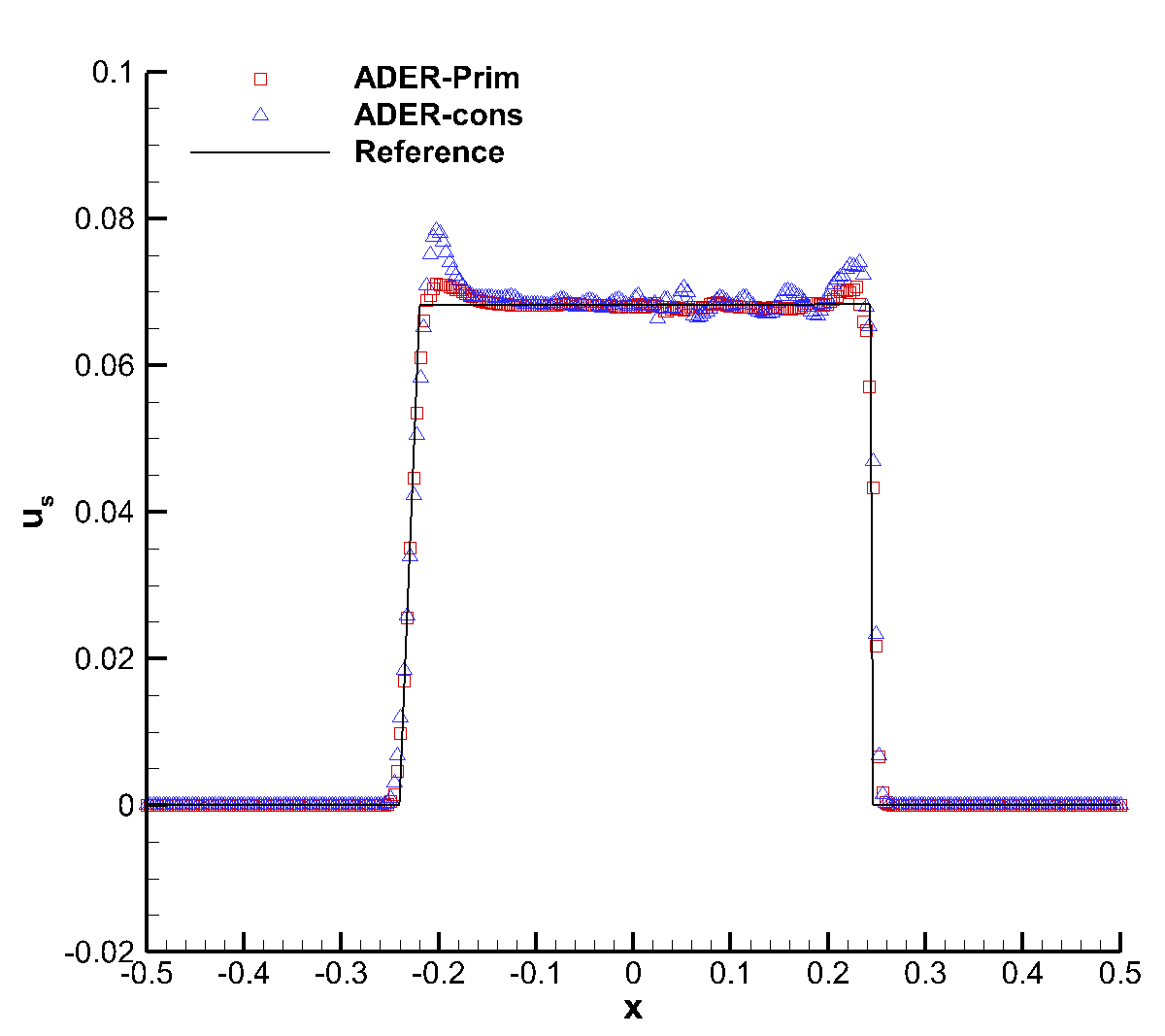}      &  
\includegraphics[width=0.45\textwidth]{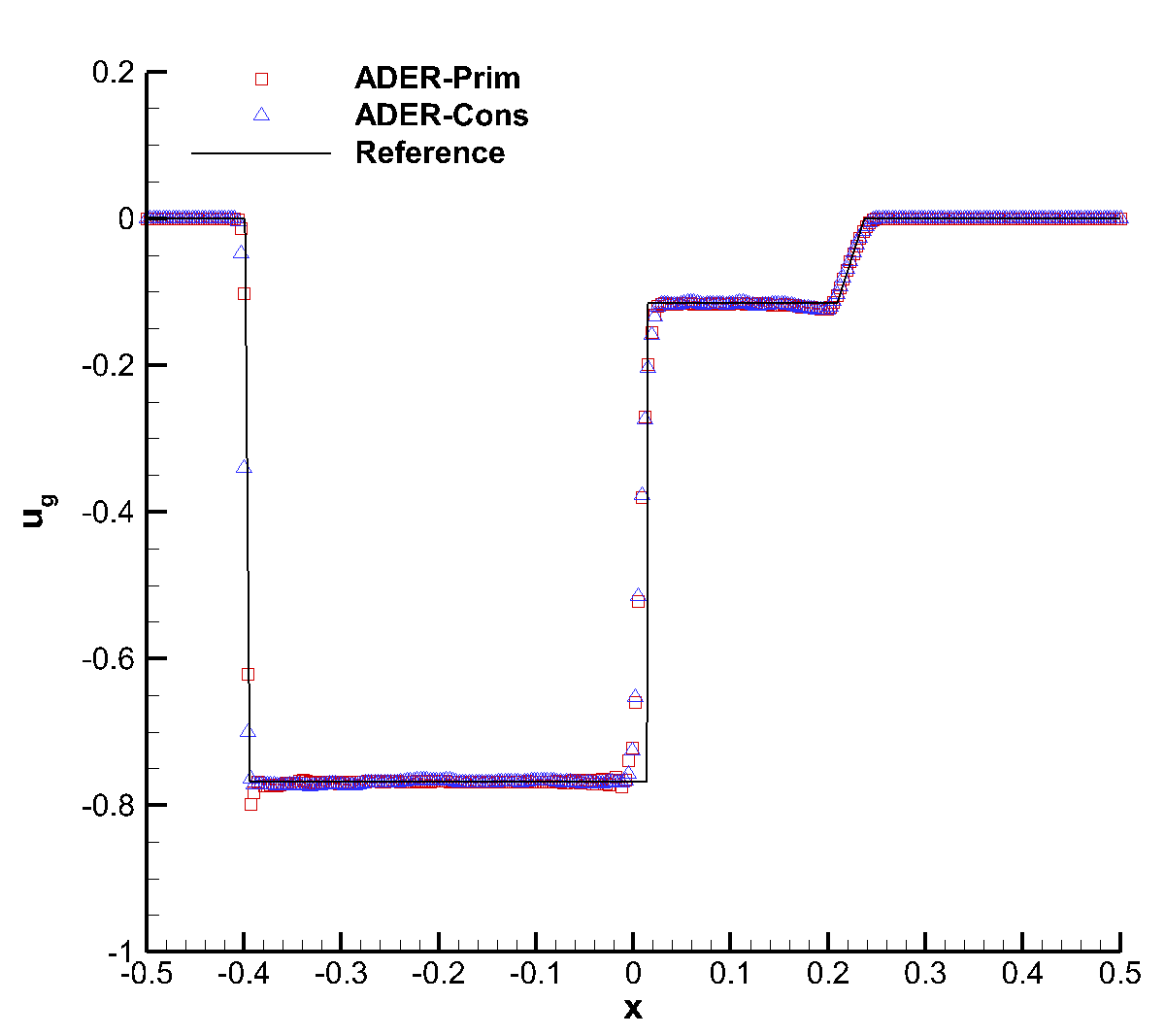}      \\ 
\includegraphics[width=0.45\textwidth]{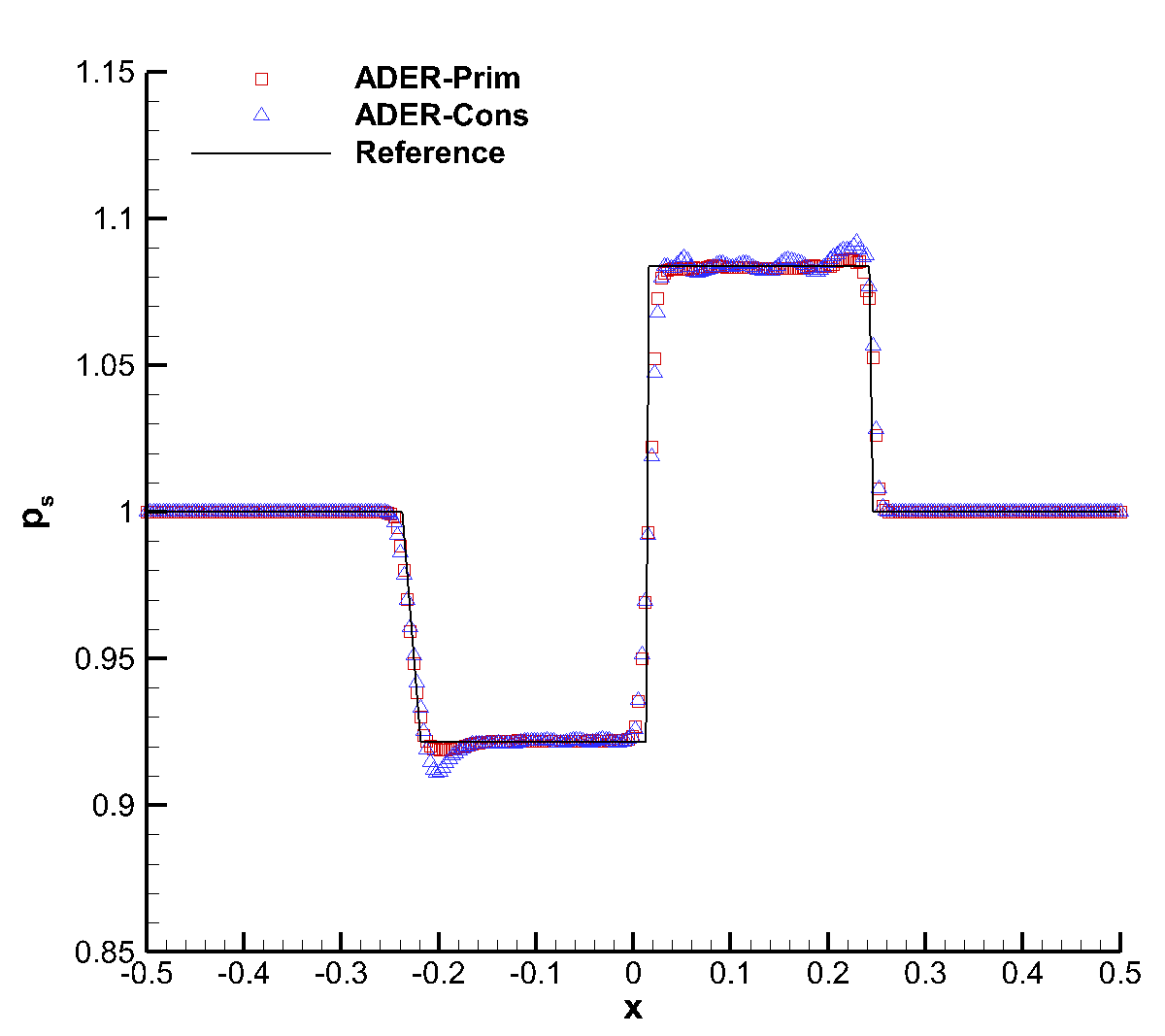}      & 
\includegraphics[width=0.45\textwidth]{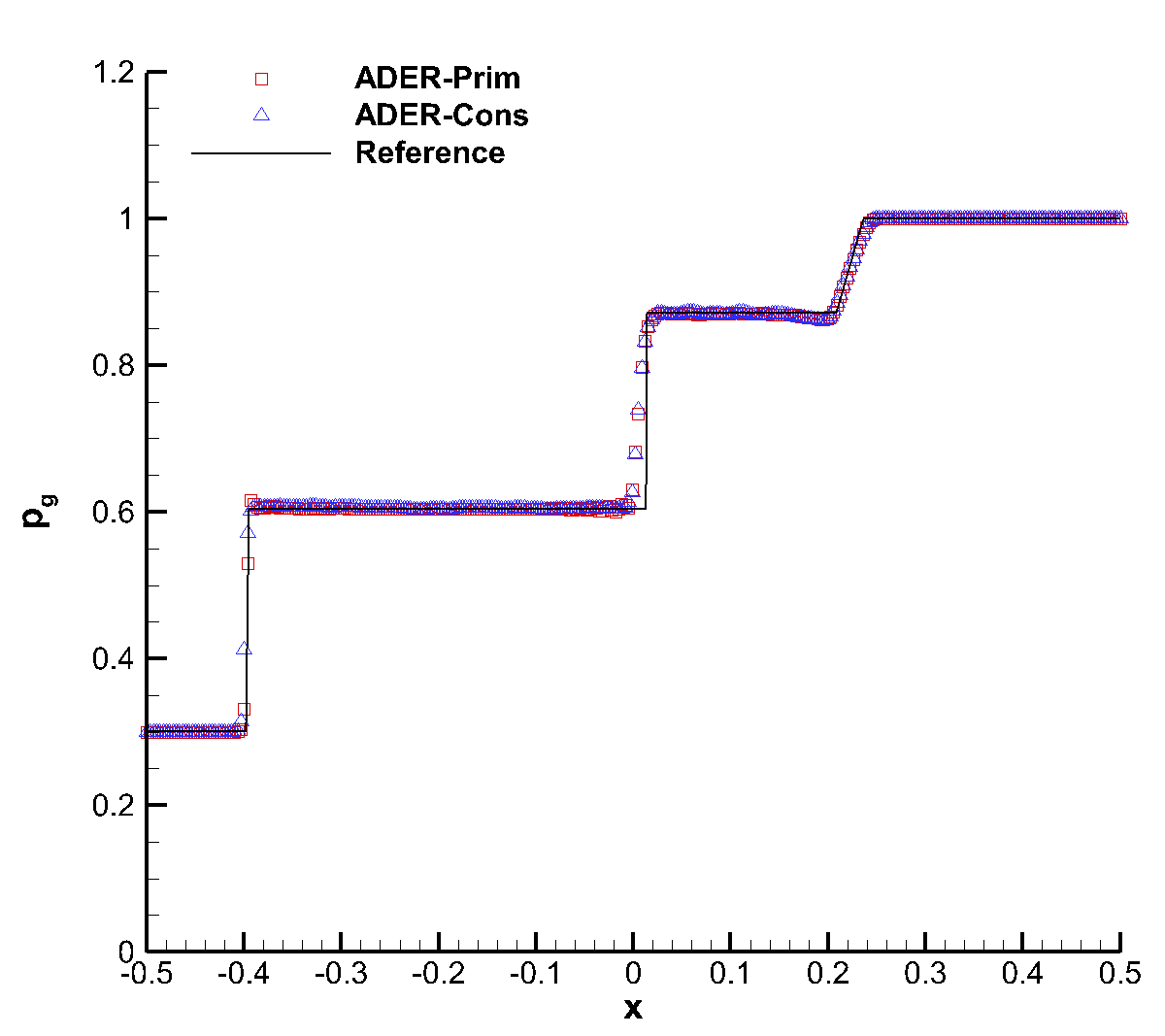}   
\end{tabular}
\caption{Results for the Baer--Nunziato Riemann problem BNRP5. The Rusanov Riemann solver has been used over a $300$ cells uniform grid.}
\label{fig.bn.rp5}
\end{center}
\end{figure}
\begin{figure}[!htbp]
\begin{center}
\begin{tabular}{cc} 
\includegraphics[width=0.45\textwidth]{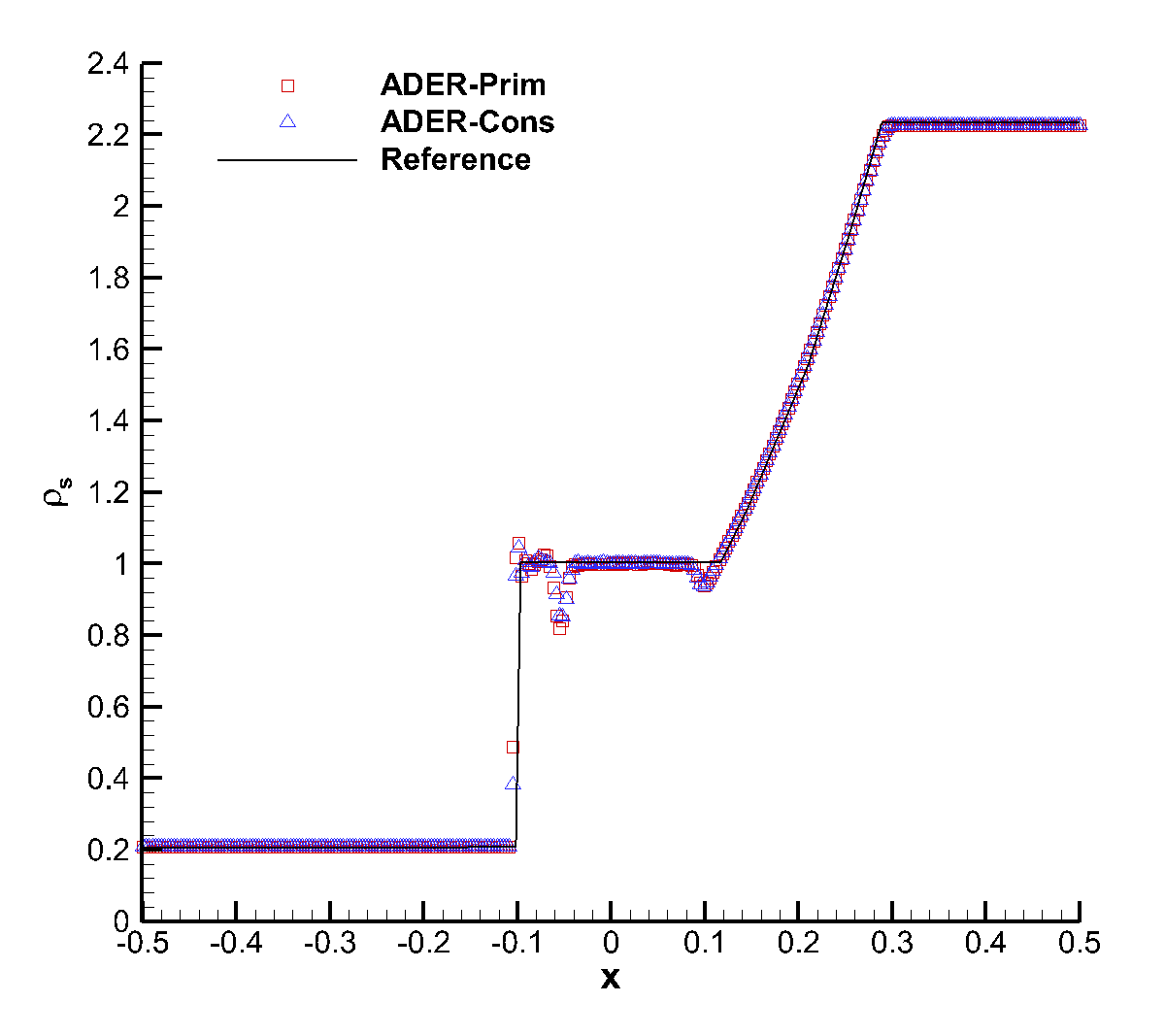}    & 
\includegraphics[width=0.45\textwidth]{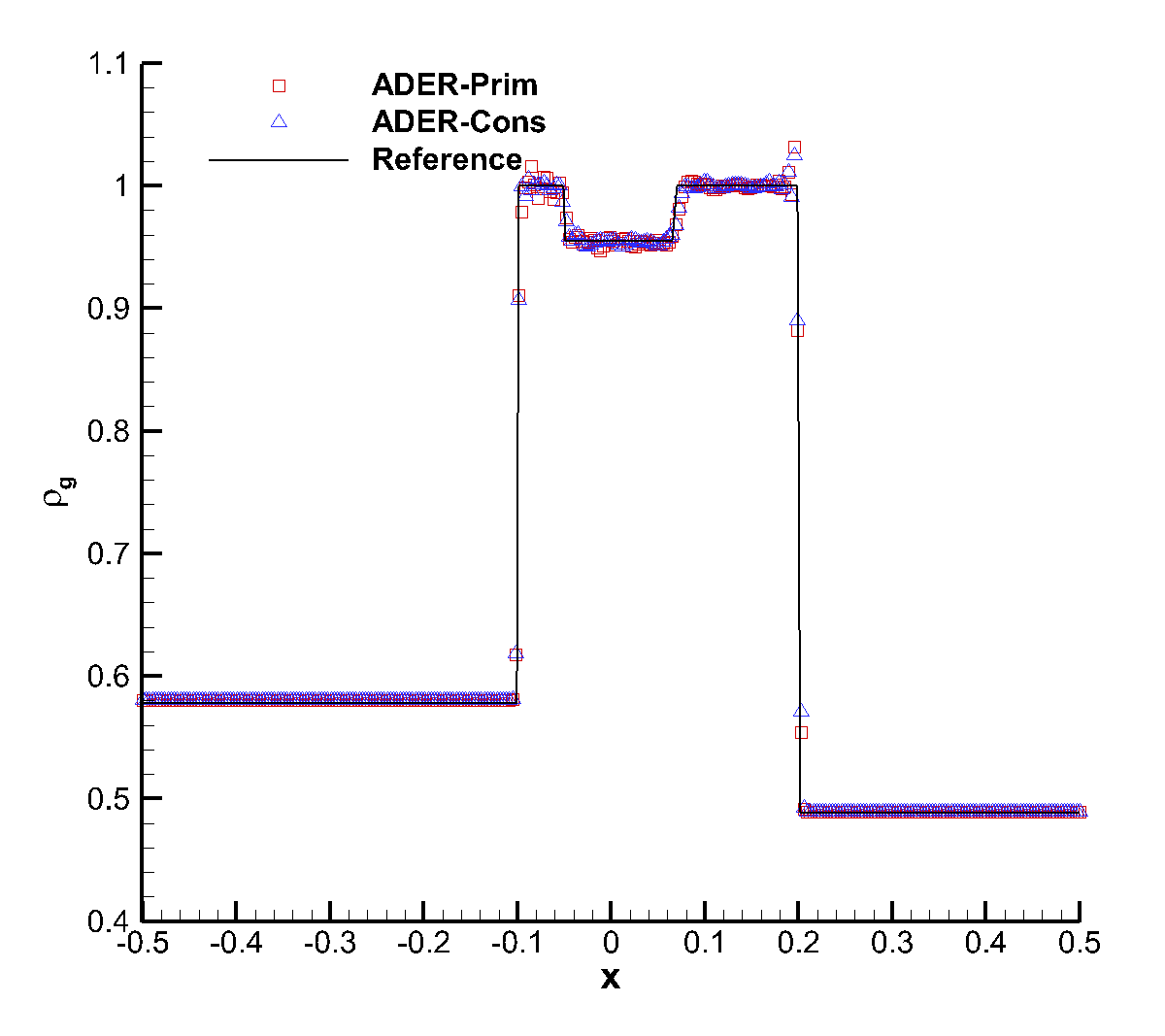}    \\ 
\includegraphics[width=0.45\textwidth]{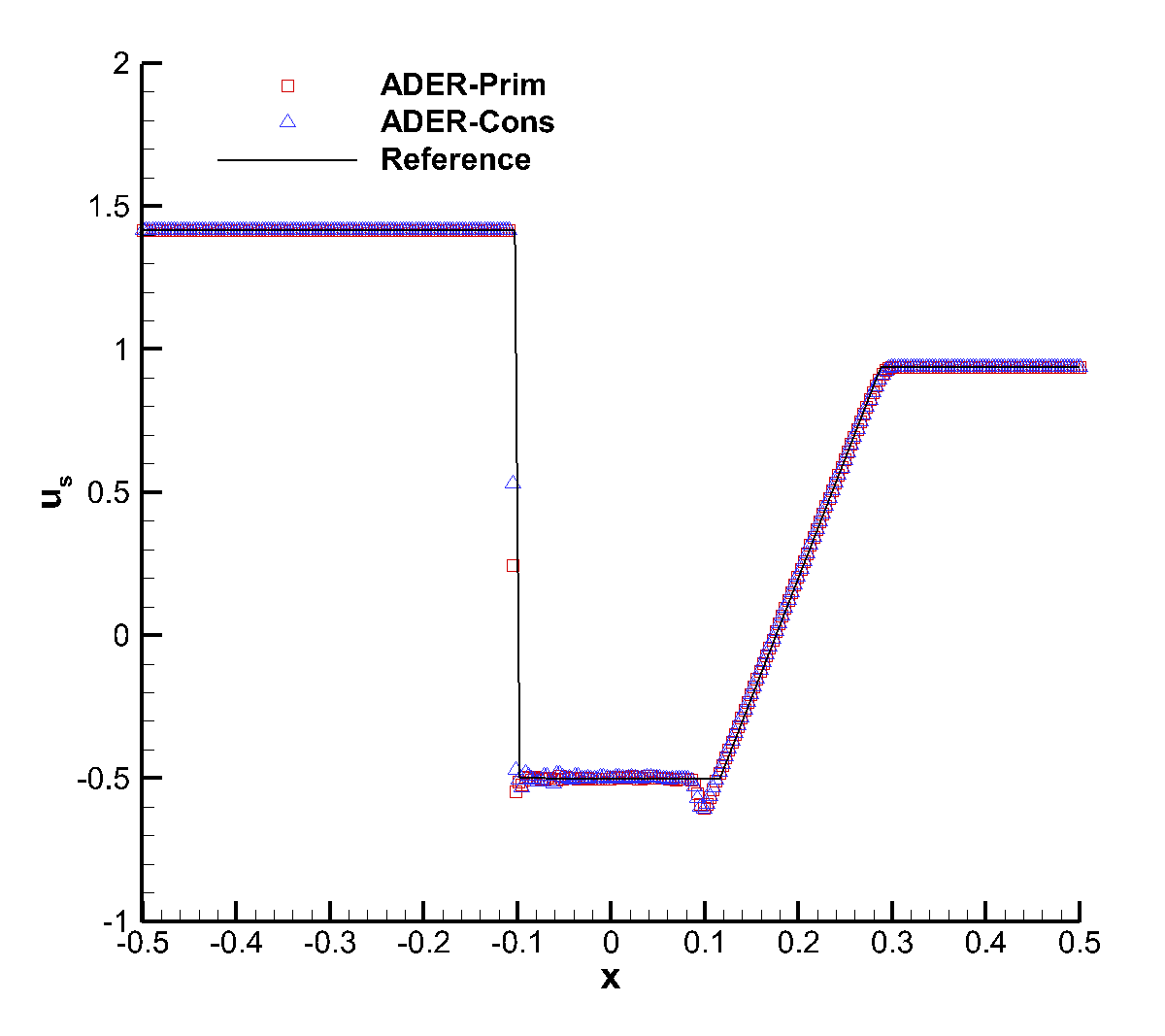}      &  
\includegraphics[width=0.45\textwidth]{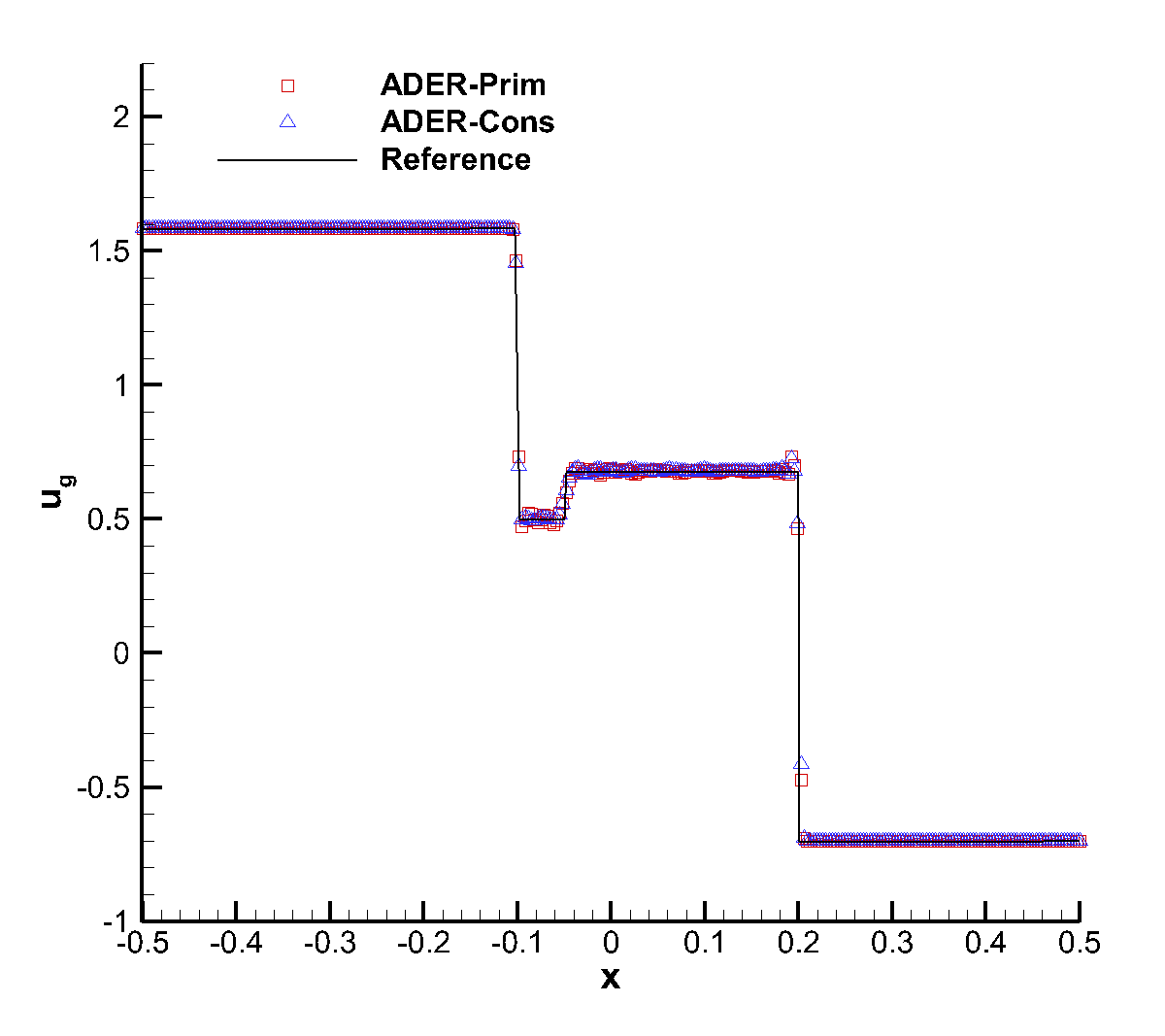}      \\ 
\includegraphics[width=0.45\textwidth]{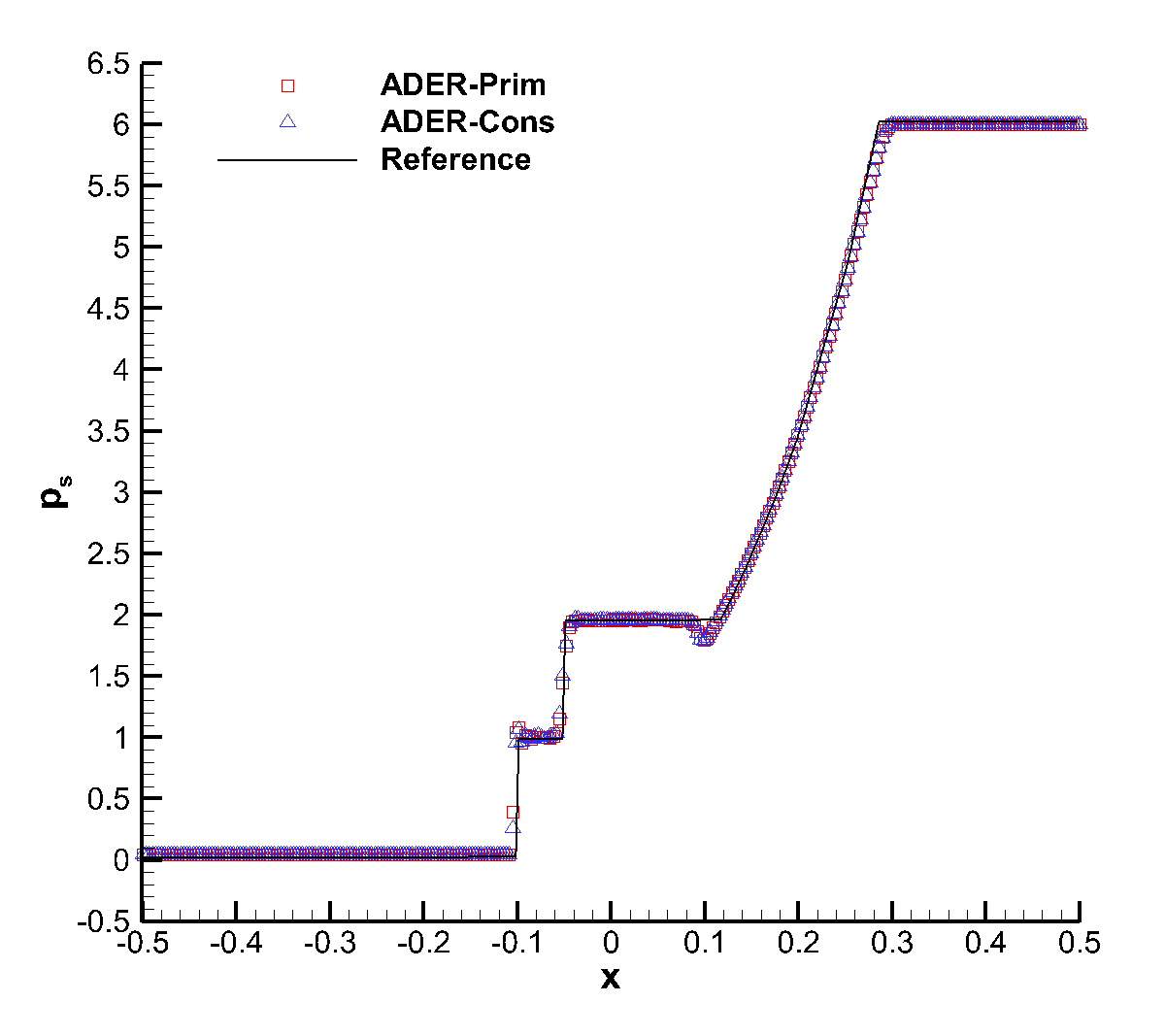}      & 
\includegraphics[width=0.45\textwidth]{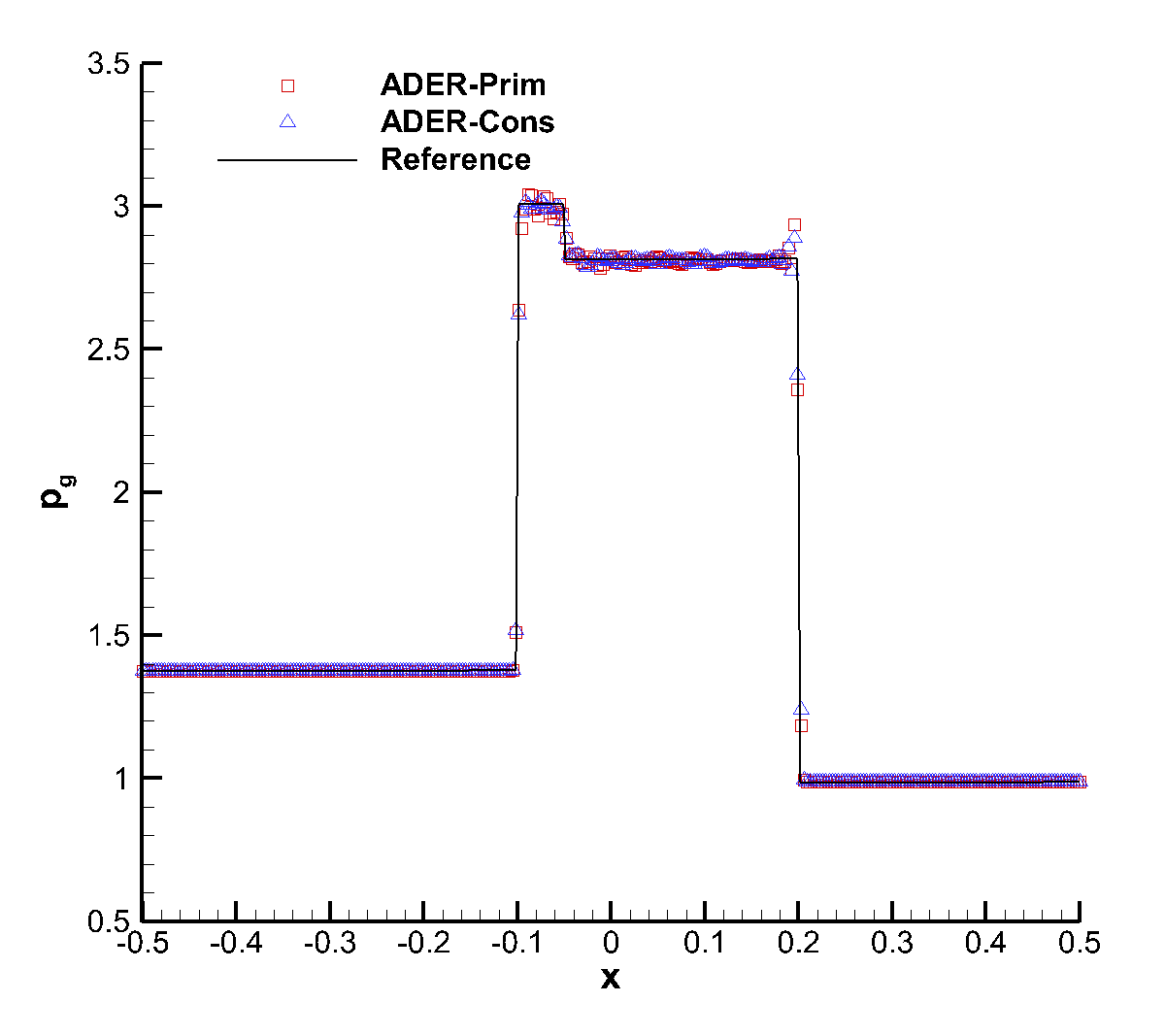}   
\end{tabular}
\caption{Results for the Baer--Nunziato Riemann problem BNRP6. The Osher Riemann solver has been used over a $300$ cells uniform grid.}
\label{fig.bn.rp6}
\end{center}
\end{figure}
The equation of state is the so-called stiffened gas equation of state, 
\begin{equation}
\label{eqn.eos} 
   \epsilon_k = \frac{p_k + \gamma_k \pi_k}{\rho_k (\gamma_k -1 )}\,,
\end{equation}
which is a simple modification of the ideal gas EOS and where $\pi_k$ expresses a reference pressure. 
For brevity, we have solved this system of equations
only for a set of one-dimensional Riemann problems,  
with initial conditions reported in Tab.~\ref{tab.rpbn.ic}.
The name of the models, BNRP1, BNRP2, etc.,  respects the numeration adopted in \cite{USFORCE2}.
A reference  solution is available for these tests, and it can be found in
\cite{AndrianovWarnecke,Schwendeman,DeledicquePapalexandris}. 
Each Riemann problem has been solved using a fourth order WENO scheme with $300$ cells uniformly distributed over the range $[-0.5;0.5]$. 
In Figs. \ref{fig.bn.rp1}-\ref{fig.bn.rp6} we have reported  the comparison 
among the solutions obtained with the ADER-Prim, with the ADER-Cons and with the exact solver.
In all the tests, with the exception of BNRP2, the ADER-Prim scheme behaves significantly better than the ADER-Cons scheme.
On several occasions, such as for $v_s$ and $v_g$ in BNRP1, or for most of the quantities in BNRP5,
the solution provided through  ADER-Cons manifest evident oscillations, which are instead strongly reduced, or even 
absent, when the ADER-Prim scheme is used. The CPU time overhead implied by ADER-Prim
is comparatively limited, and never larger than $\sim 20\%$.
%

\section{Extension to Discontinuous Galerkin and adaptive mesh refinement}
\label{sec:extension}

Although we have so far concentrated on the implementation of the new ADER-Prim scheme in the context of finite volume methods,
the same idea can be extended to Discontinuous Galerkin (DG) schemes as well. Incidentally, we note that the interest of 
computational astrophysics towards DG methods is increasing~\cite{Radice2011,Teukolsky2015}, and, especially in the relativistic context, they are expected to play  a
crucial role in the years to come.
In a sequence of papers, we have recently developed a class of 
robust DG schemes which are able to cope even with discontinuous solutions, by incorporating an aposteriori subcell limiter \cite{Dumbser2014,Zanotti2015c,Zanotti2015d}.
The whole logic can be briefly summarized as follows. 
First we assume a \emph{discrete representation} of the solution, in conserved variables, at any given time $t^n$ as
\begin{equation}
\label{eqn.ansatz.uh}
  \mathbf{u}_h(\x,t^n) = \sum_{l=0}^{N}\Phi_l(\boldsymbol{\xi}) \hat{\mathbf{u}}^n_l= \Phi_l(\boldsymbol{\xi}) \hat{\mathbf{u}}^n_l \quad \x \in T_i\,,
\end{equation}
in which the polynomials
\begin{equation}
  \Phi_l(\boldsymbol{\xi}) = \psi_p(\xi) \psi_q(\eta) \psi_r(\zeta)
\end{equation} 
are built using the spatial Lagrange interpolation polynomials already adopted for the WENO reconstruction.
The time evolution of the {\em degrees of freedom}  $\hat{\mathbf{u}}^n_l$ is then obtained after considering the 
weak form of the governing PDE, which leads to
%
\bea
\label{eqn.pde.nc.gw2}
&&\left( \int \limits_{T_i} \Phi_k \Phi_l d\x \right) \left( \hat{\mathbf{u}}_l^{n+1} -  \hat{\mathbf{u}}_l^{n} \right) +
\int\limits_{t^n}^{t^{n+1}} \int \limits_{\partial T_i} \Phi_k \, \left( {\bf \tilde f}\left(\v_h^-, \v_h^+ \right) + \frac{1}{2} \mathcal{D}\left(\v_h^-, \v_h^+ \right) 
\right) \cdot\mathbf{n} \, dS dt  \nonumber \\
&& -\int\limits_{t^n}^{t^{n+1}} \int \limits_{T_i} \nabla \Phi_k \cdot \F\left(\v_h \right) d\x dt 
+ \int\limits_{t^n}^{t^{n+1}} \int \limits_{T_i} \Phi_k{\bf B}(\v_h)  \cdot \M \nabla \v_h  \, d\x dt = \int\limits_{t^n}^{t^{n+1}} \int \limits_{T_i} \Phi_k {\bf S}\left(\v_h \right) \, d\x dt\,, 
\nonumber \\ 
\eea
%
where, just like in Eq.~(\ref{eqn.numerical.flux}), ${\bf \tilde f}$ denotes a numerical flux function and $\mathcal{D}\left(\v_h^-, \v_h^+ \right)$ a path-conservative
jump term.
Obviously, no spatial WENO reconstruction is needed within the DG framework, and the local spacetime DG predictor $\v_h(\mathbf{x},t)$ entering Eq.~(\ref{eqn.pde.nc.gw2}) 
will be computed according to the same strategy outlined in Sect.~\ref{sec:Description_of_the_predictor}. T
although acting directly over the degrees of freedom $\hat{\mathbf{p}}^n_l$ in primitive variables, which are computed from the degrees of freedom $\hat{\mathbf{u}}^n_l$ in conserved variables
simply by
\begin{equation}
   \hat{\mathbf{p}}^n_l = \V \left( \hat{\u}^n_l \right), \qquad \forall l. 
\end{equation}  
The conversion can be done in such a simple way because we use a \textit{nodal} basis $\Phi_l(\mathbf{x})$. In other words, the degrees of freedom $\hat{\mathbf{u}}^n_l$ in 
conserved variables are first converted into degrees of freedom $\hat{\p}^n_l$ in primitive variables, which are then used as initial conditions for the LSDG predictor, i.e. 
\begin{equation}
\label{LSDG-2}
{\mathbf{u}}_h(\mathbf{x},t^n)\xrightarrow{Cons2Prim} {\p}_h(\mathbf{x},t^n)\xrightarrow{LSDG} {\v}_h(\mathbf{x},t)\,,\hspace{1cm}t\in[t^n;t^{n+1}]\,.
\end{equation} 
In those cells in which the main scheme of Eq.~(\ref{eqn.pde.nc.gw2}) fails, either because unphysical values of any quantity are encountered, or
because strong oscillations appear in the solution which violate the discrete maximum principle, the computation within the troubled cell goes back to 
the time level $t^n$ and it proceeds to a complete re-calculation. In practice, a suitable subgrid is generated just within the troubled cell, and a traditional finite volume scheme 
is used on the subgrid using an alternative data representation in terms of cell averages defined for each cell of the subgrid. 
This approach and the underlying \textit{a posteriori} MOOD framework have been presented in full details in \cite{CDL1,CDL2,Dumbser2014}, to which we address 
the interested reader for a deeper understanding.

The resulting ADER-DG scheme in primitive variables can be combined with spacetime adaptive mesh refinement (AMR), in such a way to resolve the 
smallest details of the solution in highly complex flows. We refer to \cite{Zanotti2015c,Zanotti2015d} for a full account of our AMR solver in the context of ADER-DG schemes.
Here we want to show three representative test cases of the ability of the new ADER-Prim-DG scheme with adaptive mesh refinement, by considering the 
cylindrical expansion of a blast wave in a plasma with an initially uniform magnetic field (see also \cite{Komissarov1999,Leismann2005,DelZanna2007,DumbserZanotti}), 
as well as the shock problems of Leblanc, Sedov \cite{Sedov1959} and Noh \cite{noh_1987_ecs}. 

\subsection{RMHD blast wave problem} 
At time $t=0$, the rest-mass density and the pressure are $\rho=0.01$ and $p=1$, respectively, within a cylinder of radius $R=1.0$, while outside the cylinder $\rho=10^{-4}$ and $p=5\times10^{-4}$. 
Moreover, there is a constant magnetic field $B_0$ along the $x$-direction and the plasma is at rest, while a smooth ramp function between $r=0.8$ and $r=1$ modulates
the initial jump between inner and outer values, similarly to
\cite{Komissarov1999} and \cite{DelZanna2007}.

\begin{figure}
\begin{center}
\begin{tabular}{cc}  
{\includegraphics[angle=0,width=7.3cm,height=7.0cm]{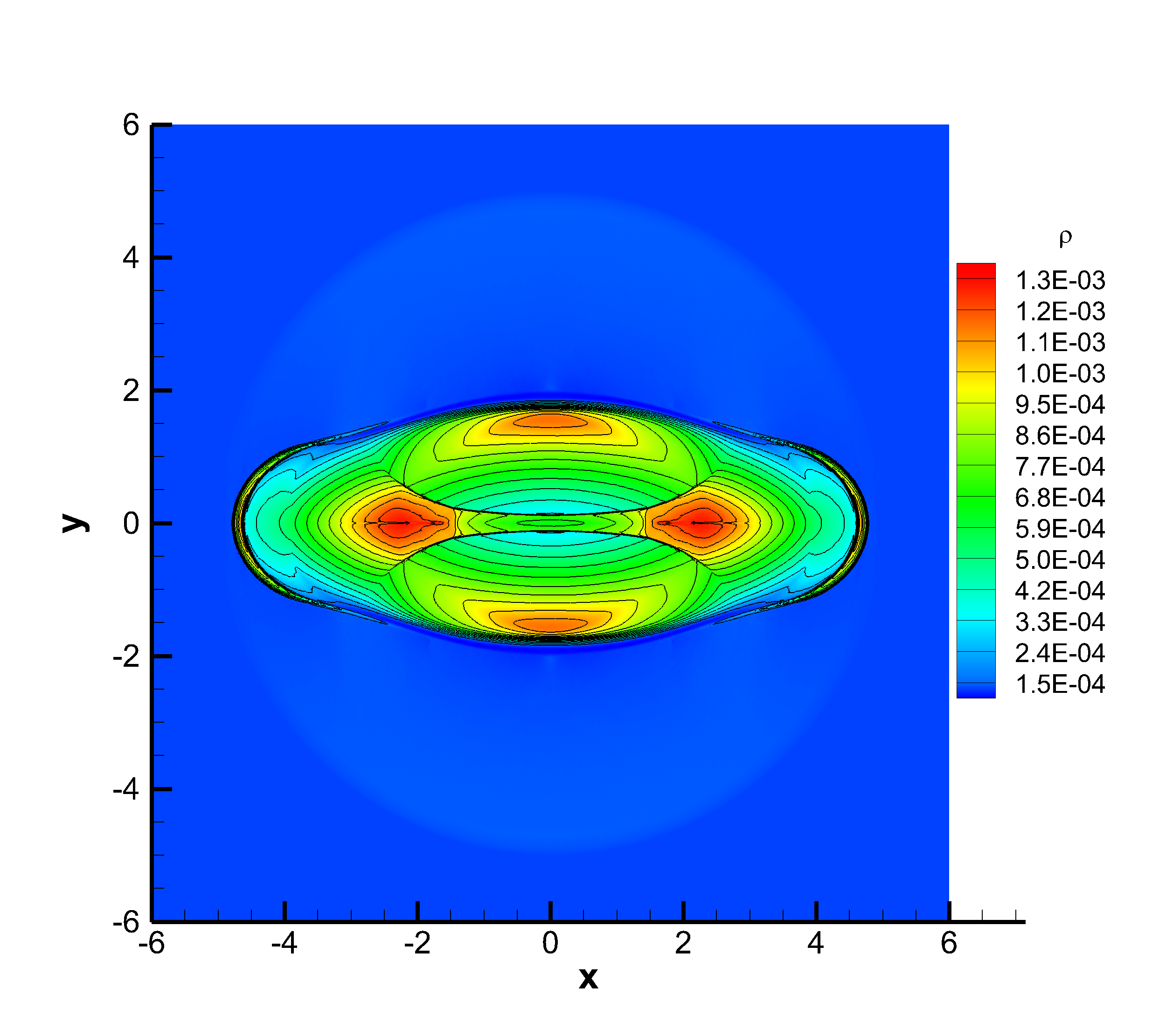}} & 
{\includegraphics[angle=0,width=7.3cm,height=7.0cm]{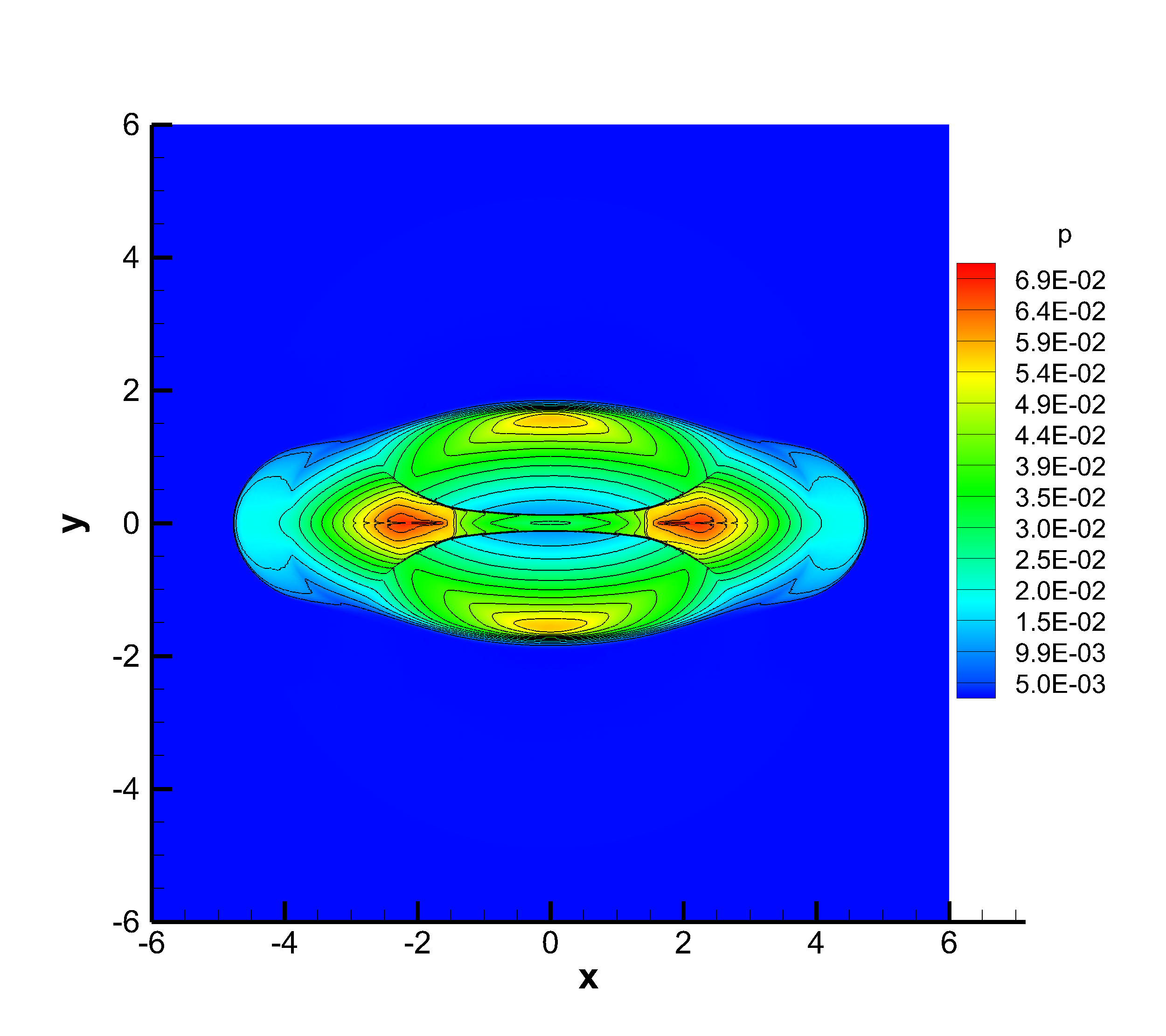}} \\
{\includegraphics[angle=0,width=7.3cm,height=7.0cm]{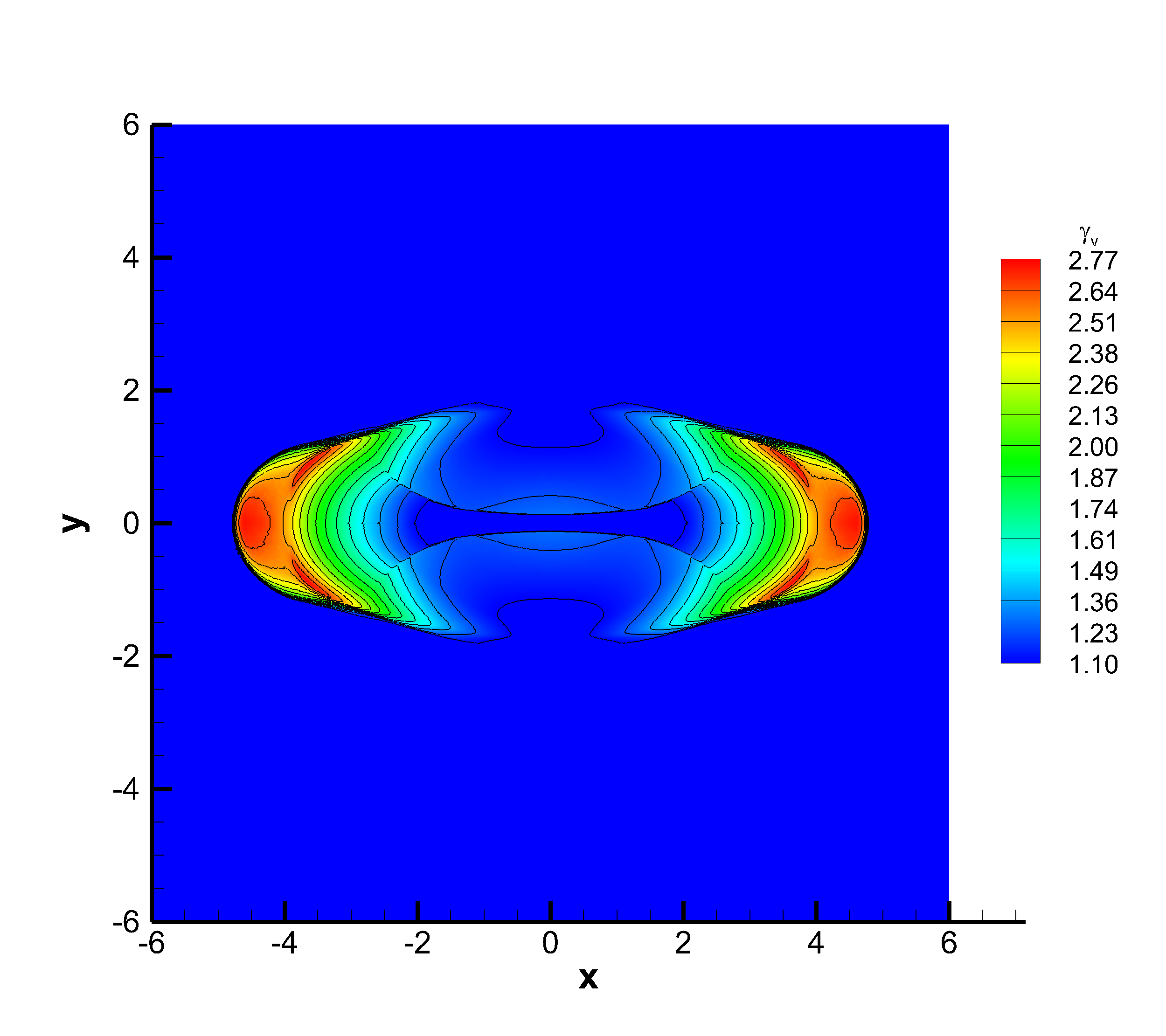}} & 
{\includegraphics[angle=0,width=7.3cm,height=7.0cm]{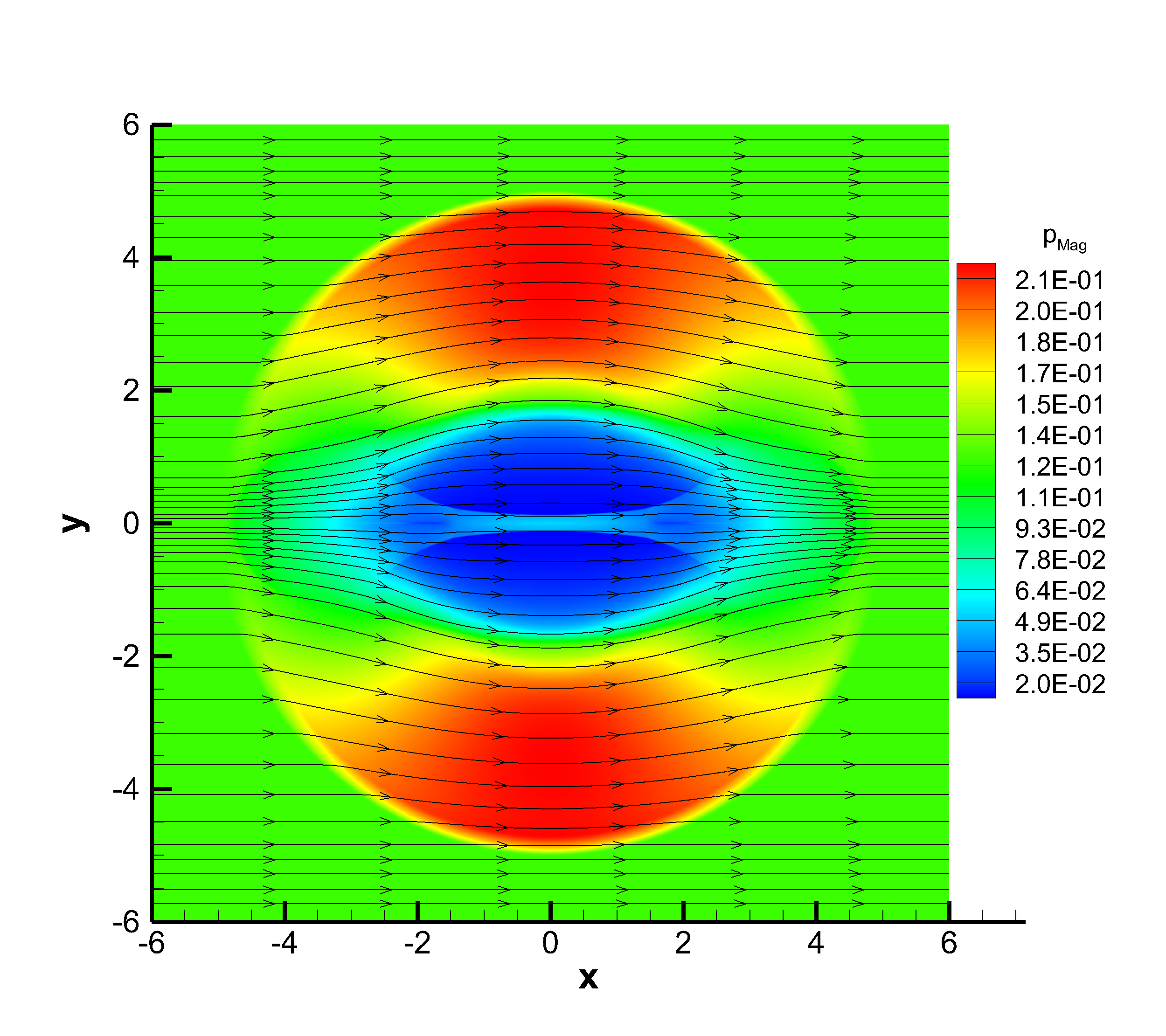}} \\ 
{\includegraphics[angle=0,width=7.3cm,height=7.0cm]{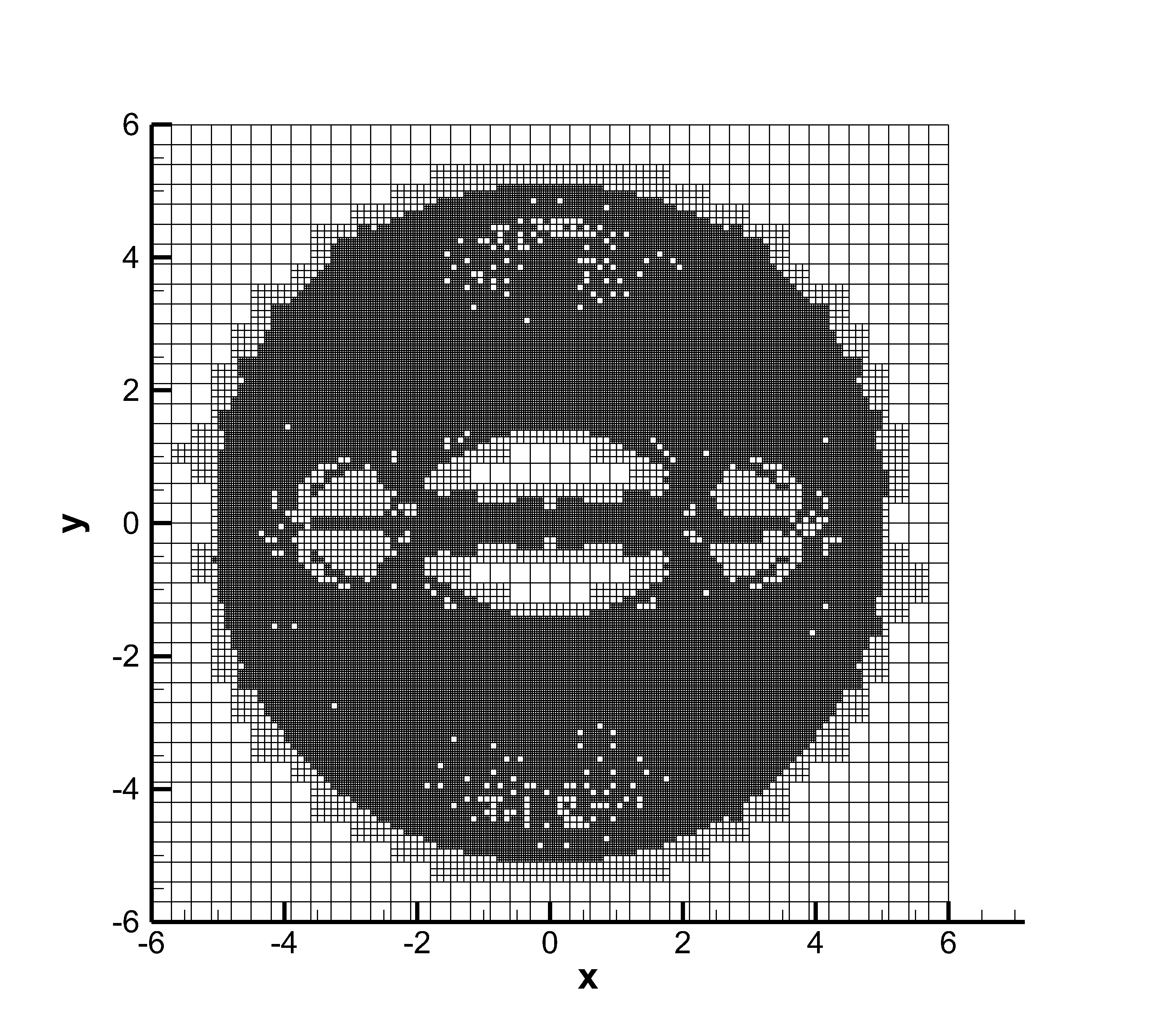}} & 
{\includegraphics[angle=0,width=7.3cm,height=7.0cm]{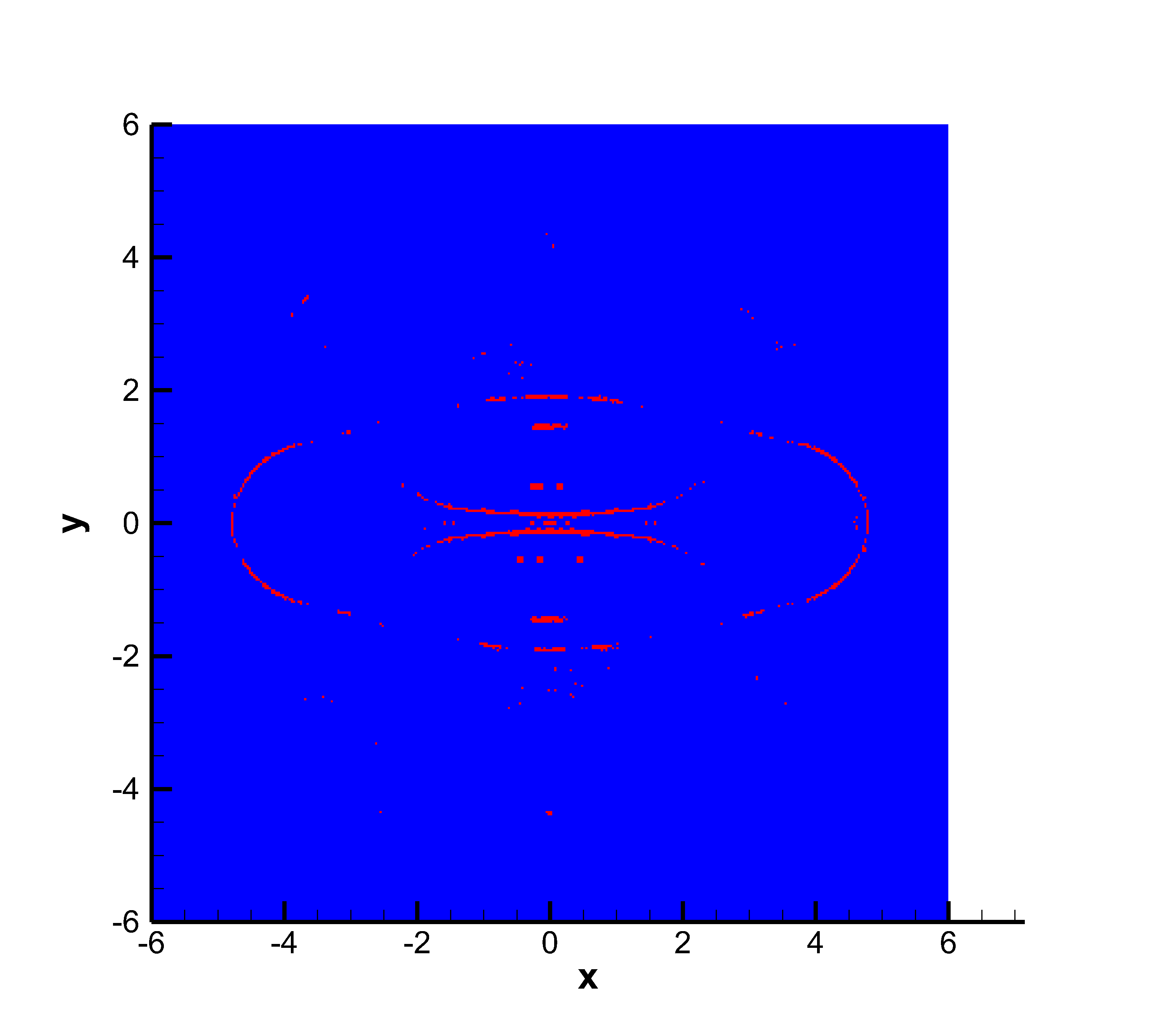}}
\end{tabular} 
\caption{Solution of the RMHD blast wave  at time $t=4.0$, obtained with the ADER-DG $\mathbb{P}_3\mathbb{P}_3$ scheme 
supplemented with the \aposteriori second order TVD subcell finite volume limiter. 
Top panels: rest-mass density (left) and thermal pressure (right). Central panels: Lorentz factor (left) and magnetic pressure (right), with magnetic field lines reported.
Bottom panels: AMR grid (left) and limiter map (right) with troubled cells marked in red and regular unlimited cells marked in blue. 
}
\label{fig:RMHD-BlastWave}
\end{center}
\end{figure}

The computational domain is
$\Omega = [-6,6]\times[-6,6]$, and the problem has been solved over an initial coarse mesh
with $40\times40$ elements. During the evolution the mesh is adaptively refined using a refinement factor along each direction $\mathfrak{r}=3$ and two levels of refinement. 
A simple Rusanov Riemann solver has been adopted, in combination with the $\mathbb{P}_3\mathbb{P}_3$ version of the ADER-DG scheme. 
On the subgrid we are free to choose any finite volume scheme that we wish, and for this specific test
we have found convenient to adopt
a  second-order TVD scheme.
The results for $B_0=0.5$ are shown in Fig.~\ref{fig:RMHD-BlastWave}, which
reports the rest-mass density, the thermal pressure, the Lorentz factor and 
the magnetic pressure at time $t=4.0$. At this time, the solution is composed by 
an external circular fast shock wave, which is hardly visible in the rest mass density, and a reverse shock wave, which is compressed along the $y$-direction.
The magnetic field is mostly confined between these two waves, as it can be appreciated from the 
contour plot of the magnetic pressure.
The two bottom panels of the figure show the AMR grid (bottom left) and the map of the limiter (bottom right).
In the latter we have used the red color to highlight those cells which required the activation of
the limiter over the subgrid, while the blue color is for the regular cells. In practice, the limiter is only needed at the inner  shock front, while the external shock front
is so weak that the limiter is only occasionally activated. 
These results confirm the ability of the new ADER-Prim scheme to work also in combination with Discontinuous Galerkin methods, and with complex
systems of equations like RMHD.

\subsection{Leblanc, Sedov and Noh problem} 

\begin{figure}
\begin{center}
\begin{tabular}{cc}  
{\includegraphics[angle=0,width=0.45\textwidth]{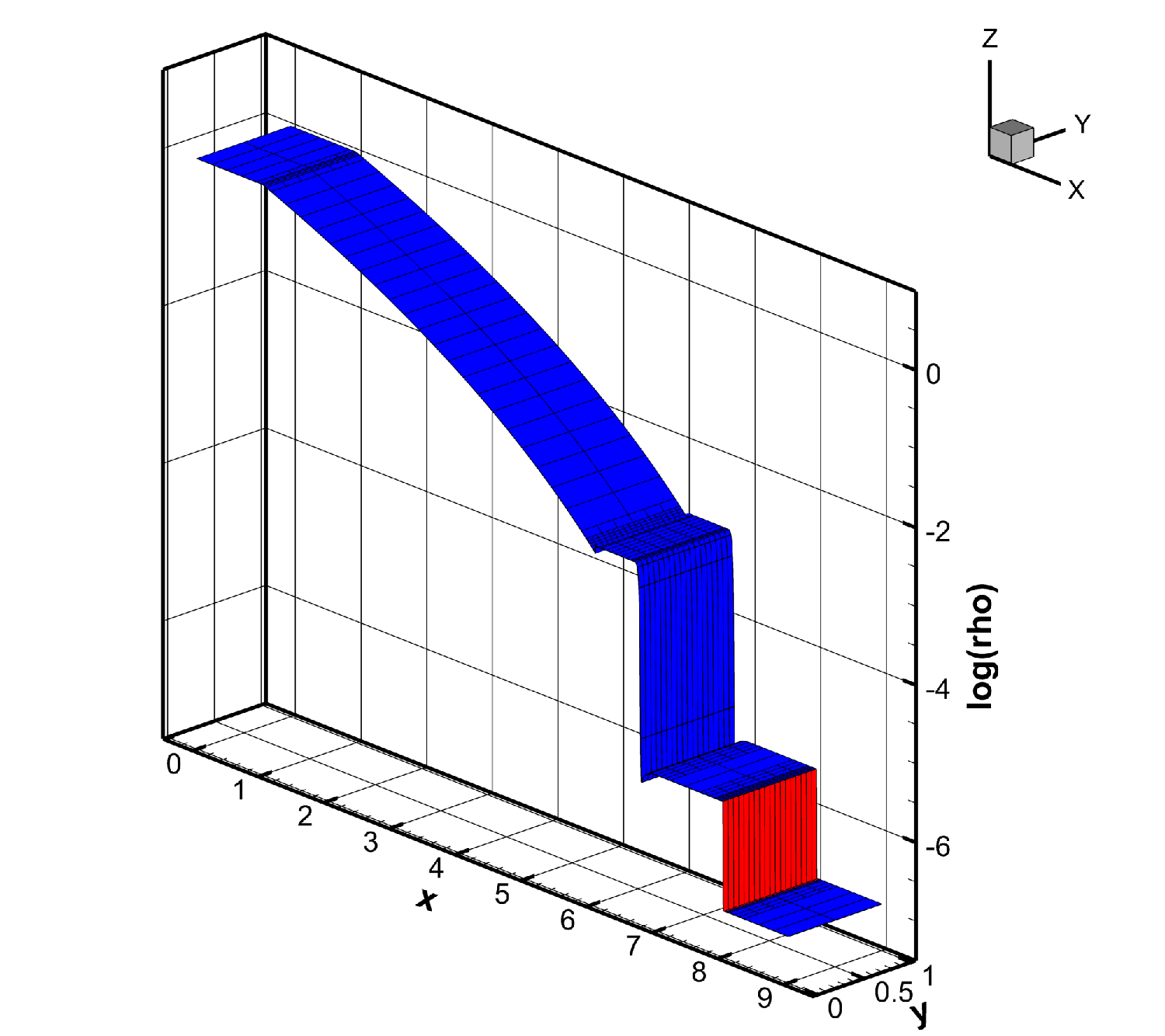}} & 
{\includegraphics[angle=0,width=0.45\textwidth]{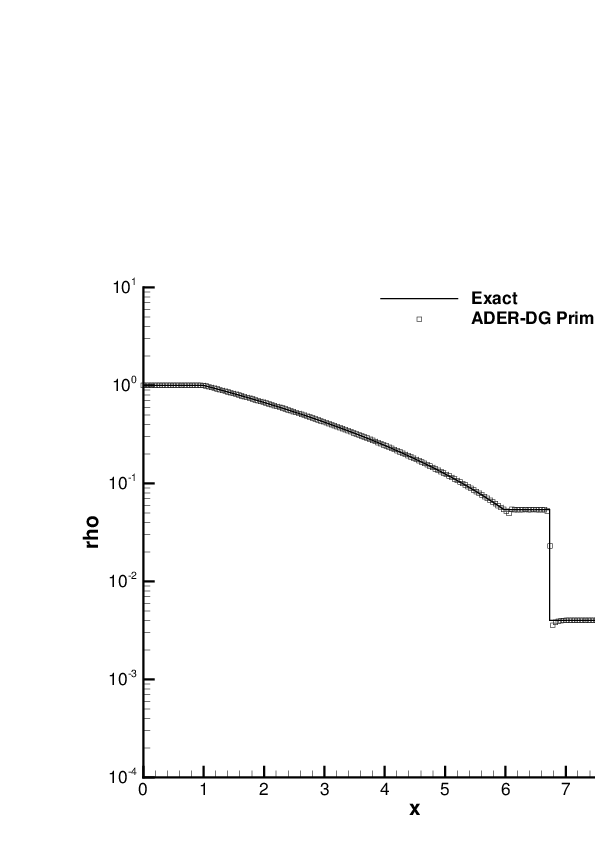}}     \\
{\includegraphics[angle=0,width=0.45\textwidth]{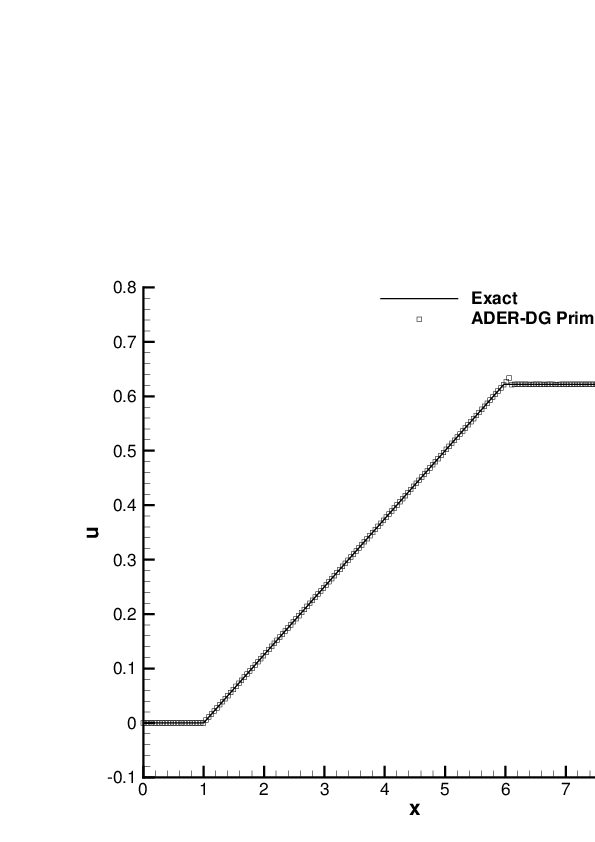}}       &
{\includegraphics[angle=0,width=0.45\textwidth]{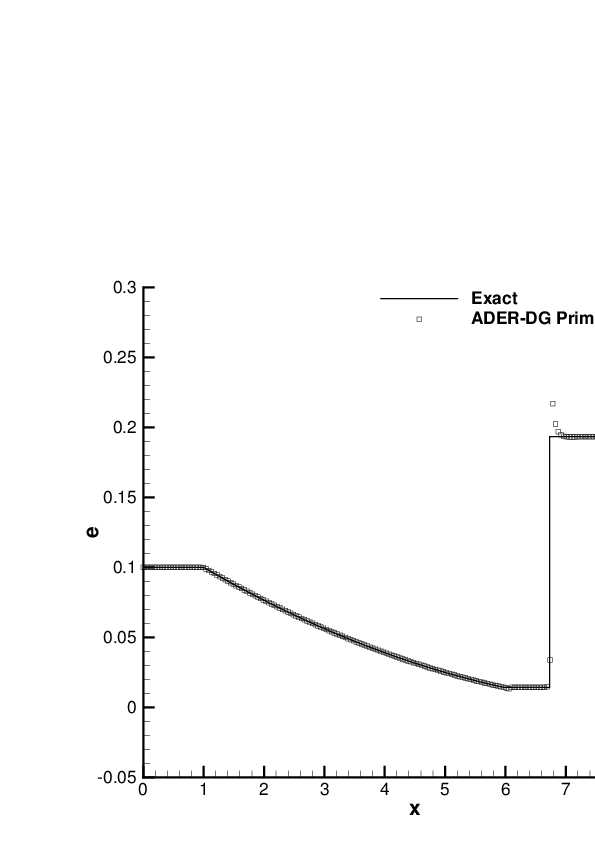}}     
\end{tabular} 
\caption{Solution of the Leblanc shock tube problem at time $t=6.0$, obtained with the ADER-DG $\mathbb{P}_9\mathbb{P}_9$ scheme 
supplemented with the \aposteriori second order TVD subcell finite volume limiter. 
Top left: Troubled cells highlighted in red and unlimited cells in blue. Top right to bottom right: Comparison with the exact solution 
using a 1D cut through the 2D solution on 200 equidistant sample points for density, velocity and internal energy.  
}
\label{fig:Leblanc}
\end{center}
\end{figure}

\begin{figure}
\begin{center}
\begin{tabular}{cc} 
{\includegraphics[angle=0,width=0.45\textwidth]{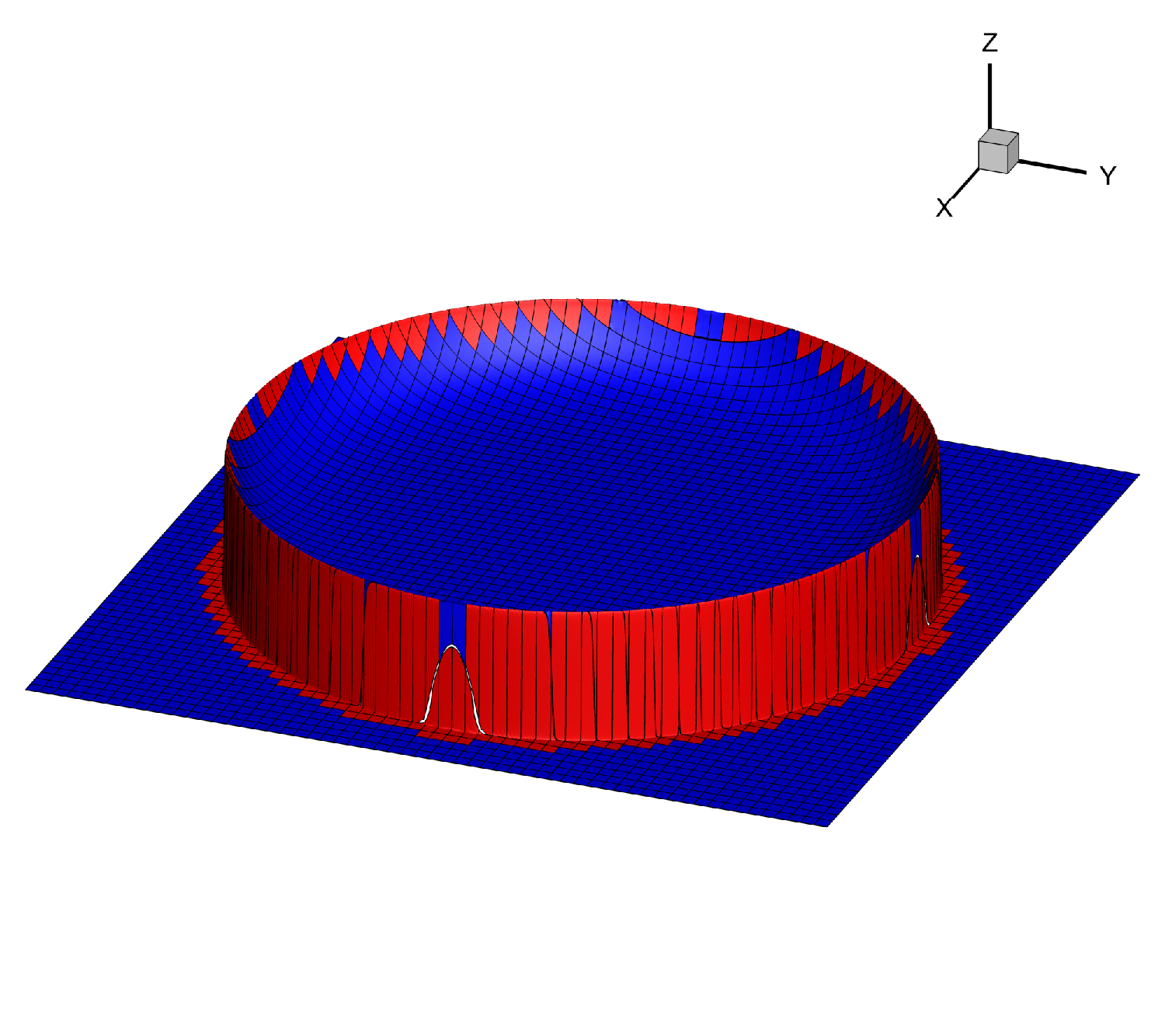}} & 
{\includegraphics[angle=0,width=0.45\textwidth]{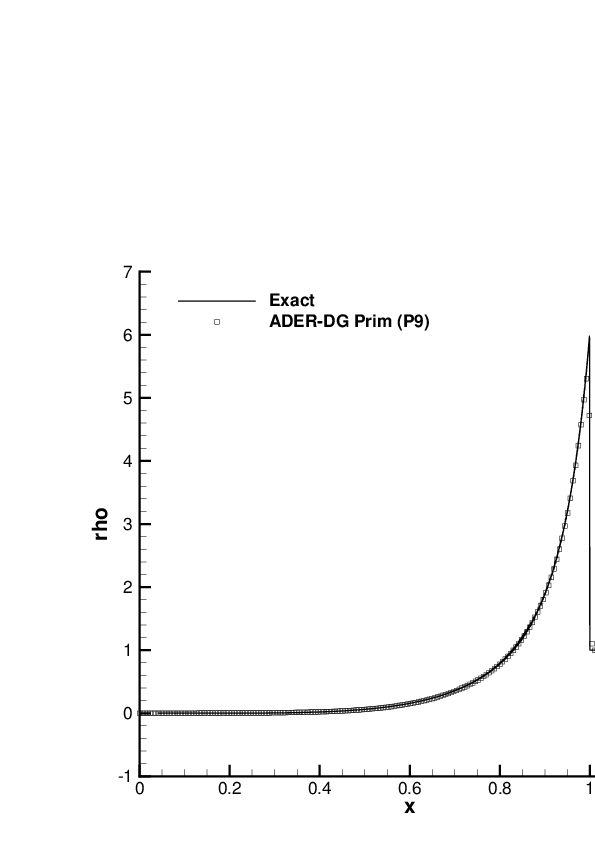}}
\end{tabular} 
\caption{Solution of the Sedov problem at time $t=1.0$, obtained with the ADER-DG $\mathbb{P}_9\mathbb{P}_9$ scheme 
supplemented with the \aposteriori second order TVD subcell finite volume limiter. 
Left: Troubled cells highlighted in red and unlimited cells in blue. Right: Comparison with the exact solution along the $x$-axis. 
}
\label{fig:Sedov}
\end{center}
\end{figure}

\begin{figure}
\begin{center}
\begin{tabular}{cc} 
{\includegraphics[angle=0,width=0.45\textwidth]{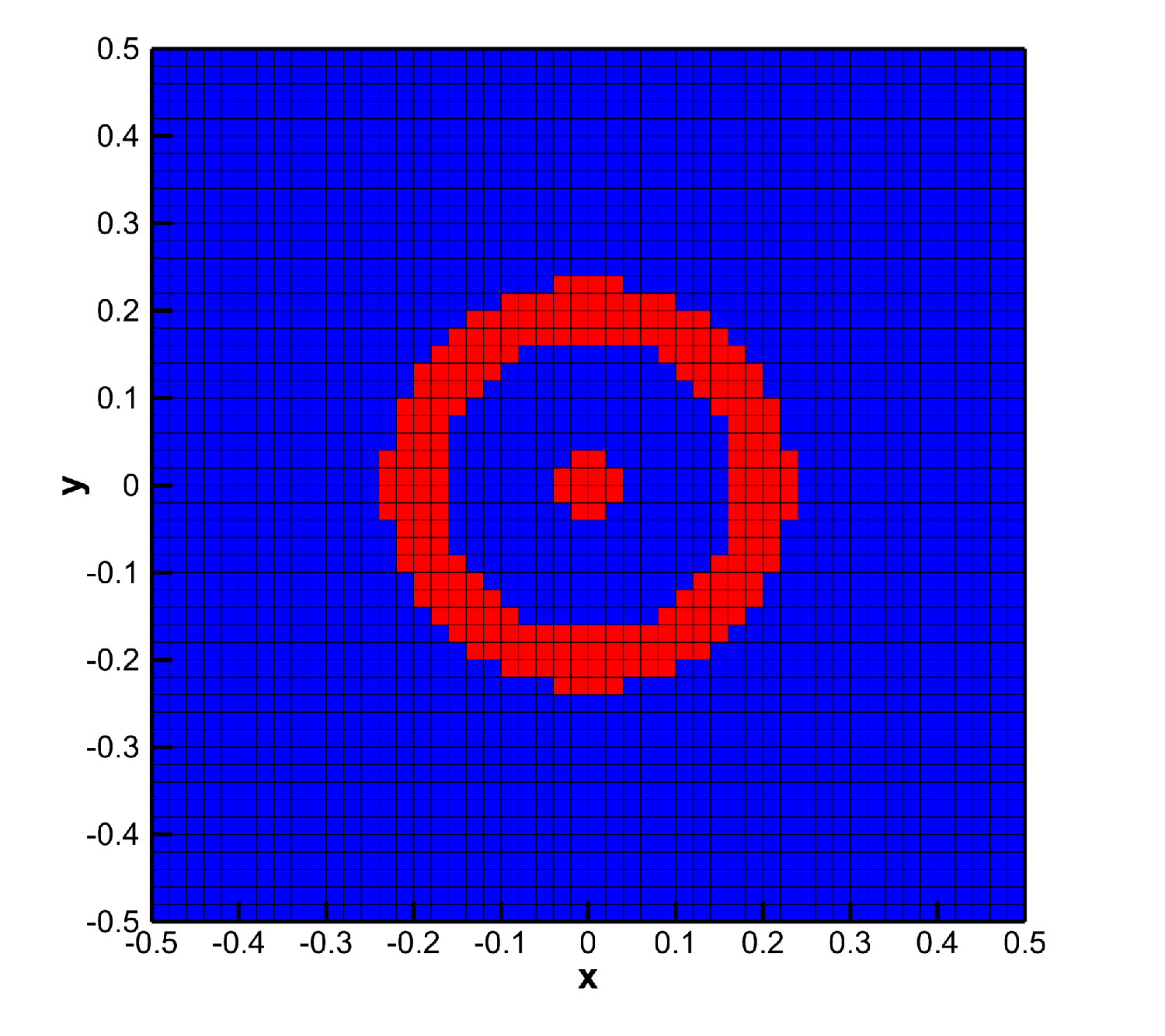}} & 
{\includegraphics[angle=0,width=0.45\textwidth]{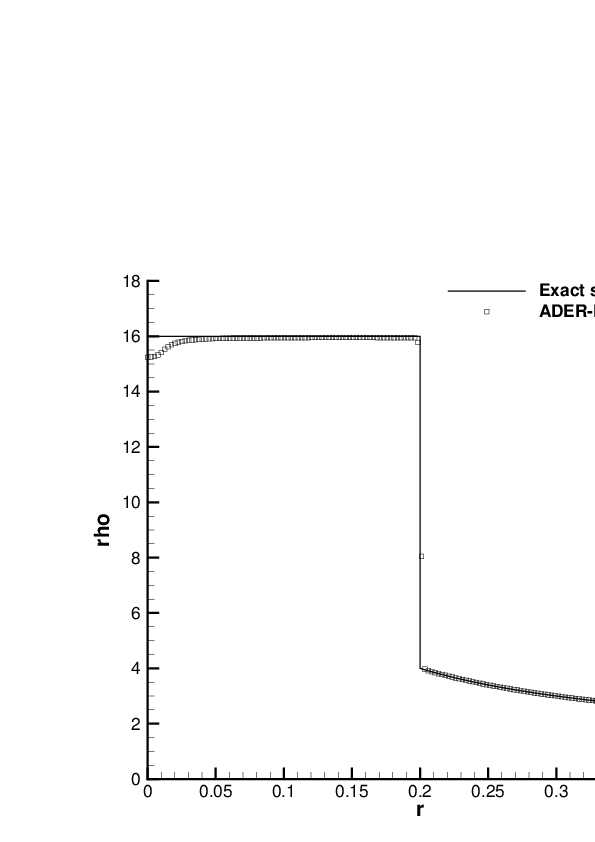}}
\end{tabular} 
\caption{Solution of the Noh problem at time $t=0.6$, obtained with the ADER-DG $\mathbb{P}_9\mathbb{P}_9$ scheme 
supplemented with the \aposteriori second order TVD subcell finite volume limiter. 
Left: Troubled cells highlighted in red and unlimited cells in blue. Right: Comparison with the exact solution along the $x$-axis. 
}
\label{fig:Noh}
\end{center}
\end{figure}

Here we solve again the classical Euler equations of compressible gas dynamics on a rectangular domain for the Leblanc problem and on a circular domain in the case of 
the shock problems of Sedov and Noh. The initial conditions are detailed in \cite{Dumbser-Uuriintsetseg2013,LagrangeMHD,Lagrange3D}. For the low pressure region that is 
present in the above test problems, we use $p=10^{-14}$ for the Leblanc and the Noh problem. The computational results obtained with very high order ADER-DG $\mathbb{P}_9\mathbb{P}_9$ 
schemes are depicted in Figures  \ref{fig:Leblanc}, \ref{fig:Sedov} and \ref{fig:Noh}, showing an excellent agreement with the exact solution in all cases, apart from the 
overshoot in the case of the Leblanc shock tube. We stress that all test problems are extremely severe and therefore clearly demonstrate the robustness of the new approach. 

\section{Conclusions}
\label{sec:conclusions}
%
The new version of ADER schemes introduced in \cite{DumbserEnauxToro} relies on a local space-time discontinuous Galerkin predictor, which is then used for the 
computation of high order accurate fluxes and sources. This approach has the advantage over classical Cauchy-Kovalewski based ADER schemes 
\cite{toro1,toro3,toro4,titarevtoro,Toro:2006a,dumbser_jsc,taube_jsc} that it is in principle 
applicable to general nonlinear systems of conservation laws. However, for hyperbolic systems in which the conversion from conservative to primitive variables 
is not analytic but only available numerically, a large number of such expensive conversions must be performed, namely one for each space-time quadrature point for 
the integration of the numerical fluxes over the element interfaces and one for each space-time degree of freedom in the local space-time DG predictor. 

Motivated by this limitation, we have designed a new version of ADER schemes, valid primarily for finite volume schemes but extendible also to the discontinuous Galerkin 
finite element framework, in which both the spatial WENO reconstruction and the subsequent local space-time DG predictor
act on the primitive variables. In the finite volume context this can be done by performing a double WENO reconstruction for each cell. In the first WENO step, piece-wise 
polynomials of the conserved variables are computed from the cell averages in the usual way. Then, these reconstruction polynomials are simply \textit{evaluated} in the 
cell centers, in order to obtain \textit{point values} of the conserved variables. After that, a single conversion from the conserved to the primitive variable is needed 
in each cell. Finally, a second WENO reconstruction acts on these point values and provides piece-wise polynomials of the primitive variables. The local space-time 
discontinuous Galerkin predictor must then be reformulated in a non-conservative fashion, supplying the time evolution of the reconstructed polynomials for the primitive 
variables. 

For all systems of equations that we have explored, classical Euler, relativistic hydrodynamics (RHD) and magnetohydrodynamics (RMHD) and the Baer--Nunziato equations, 
we have noticed a significant reduction of spurious oscillations provided by the new reconstruction in primitive variables with respect to traditional reconstruction 
in conserved variables. This effect is particularly evident for the Baer--Nunziato equations. In the relativistic regime, there is also an improvement in the ability 
of capturing the position of shock waves (see Fig.~\ref{fig:shock-tube-RS}). To a large extent, the new primitive formulation provides results that are comparable 
to reconstruction in characteristic variables. 

Moreover, for systems of equations in which the conversion from the conserved to the primitive variables cannot be obtained in closed form, such as for the 
RHD and RMHD equations, there is an advantage in terms of computational efficiency, with reductions of the CPU time around $\sim 20\%$, or more. 
We have also introduced an additional improvement, namely the implementation of a new initial guess for the LSDG predictor, which is based on 
an extrapolation in time, similar to Adams--Bashforth-type ODE integrators. This new initial guess is typically faster than those traditionally available, but 
it is also less robust in the presence of strong shocks. 

We predict that the new version of ADER based on primitive variables  will become the standard ADER scheme in the relativistic framework. 
This may become particularly advantageous for high energy astrophysics, in which both high accuracy and high computational efficiency are required.



\begin{backmatter}

\section*{Competing interests}
  The authors declare that they have no competing interests.



\section*{Acknowledgements}
\begin{tabular}{lr} 
\begin{minipage}[c]{0.8\textwidth}
The research presented in this paper was financed by i) the European Research Council 
(ERC) under the European Union's Seventh Framework Programme (FP7/2007-2013)  
with the research project \textit{STiMulUs}, ERC Grant agreement no. 278267 and 
ii) it has received funding from the European Union's Horizon 2020 Research and Innovation 
Programme under grant agreement No. 671698 (call FETHPC-1-2014, project \textit{ExaHyPE}).    
\end{minipage}  
& 
\begin{minipage}[c]{0.2\textwidth}
\includegraphics[angle=0,width=0.65\textwidth]{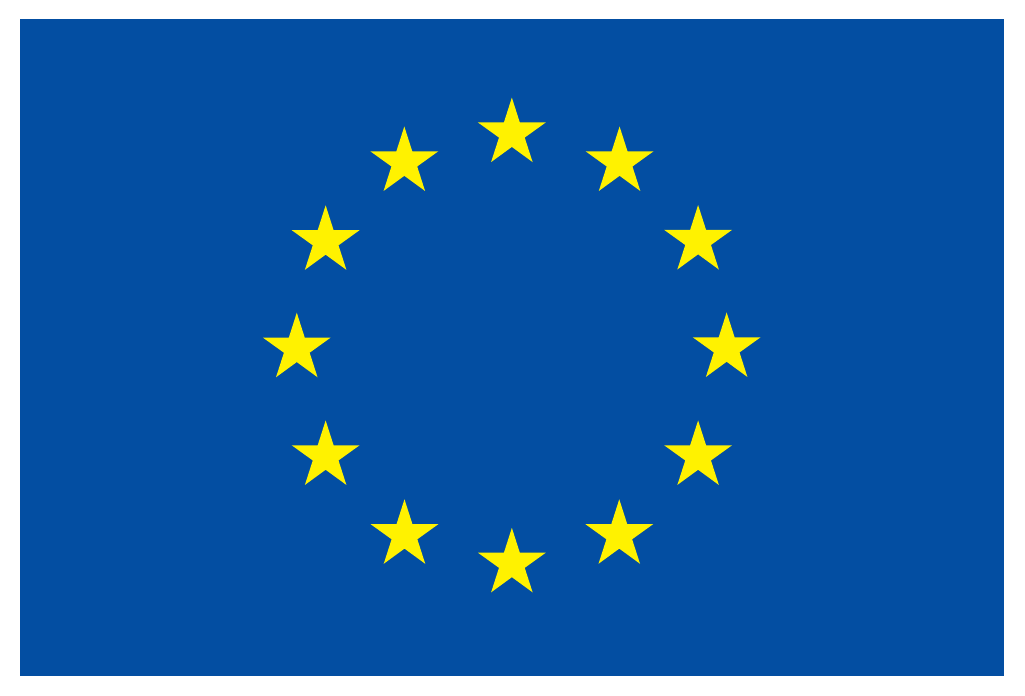}
\end{minipage}  
\end{tabular} 

We are grateful to Bruno Giacomazzo and Luciano Rezzolla for providing the numerical code for the exact 
solution of the Riemann problem in RMHD. 
We would also like to acknowledge PRACE for awarding access to 
the SuperMUC supercomputer based in Munich (Germany) at the Leibniz Rechenzentrum (LRZ), 
and ISCRA, for awarding access to the FERMI supercomputer based in Casalecchio (Italy).




\end{backmatter}
\end{document}